\documentclass{article}
\usepackage[utf8]{inputenc}

\title{Moduli space of simple polygons}
\author{Ahtziri Gonz\'alez and Manuel Sedano-Mendoza}
\date{}

\usepackage{tikzit}
\usepackage{graphicx}
\usepackage{amssymb}
\usepackage{amsthm}
\usepackage{amsfonts,amsmath,amsbsy,latexsym,mathrsfs}
\usepackage{lineno}
\usepackage{enumerate}
\usepackage[section]{placeins}
\usepackage{multicol}

\usepackage{pgfplots}
\pgfplotsset{compat=1.15}
\usepackage{mathrsfs}
\usetikzlibrary{arrows}
\usepackage[normalem]{ulem}
\usepackage{tikz}
\usepackage{pgf,tikz}
\usetikzlibrary{arrows}
\usepackage{pgfplots}
\usepackage{arcs}
\usetikzlibrary{matrix}

\DeclareFontFamily{OMX}{yhex}{}
\DeclareFontShape{OMX}{yhex}{m}{n}{<->yhcmex10}{}
\DeclareSymbolFont{yhlargesymbols}{OMX}{yhex}{m}{n}
\DeclareMathAccent{\wideparen}{\mathord}{yhlargesymbols}{"F3}

\makeatletter
\DeclareFontFamily{U}{tipa}{}
\DeclareFontShape{U}{tipa}{m}{n}{<->tipa10}{}
\newcommand{\arc@char}{{\usefont{U}{tipa}{m}{n}\symbol{62}}}%

\newcommand{\arc}[1]{\mathpalette\arc@arc{#1}}

\newcommand{\arc@arc}[2]{%
  \sbox0{$\m@th#1#2$}%
  \vbox{
    \hbox{\resizebox{\wd0}{\height}{\arc@char}}
    \nointerlineskip
    \box0
  }%
}
\makeatother

\usepackage[colorinlistoftodos]{todonotes} 

\usepackage{caption}
\usepackage{subcaption}
\usepackage{float}
\usepackage{hyperref}
\usepackage{pb-diagram}

\newtheorem{thm}{Theorem}[section]
\newtheorem{cor}[thm]{Corolary}
\newtheorem{lem}[thm]{Lemma}
\newtheorem{prop}[thm]{Proposition}
\newtheorem{defi}[thm]{Definition}
\newtheorem{conj}[thm]{Conjecture}
\newtheorem{obs}[thm]{\rm Remark}
\newtheorem{ejem}[thm]{Example}
\newtheorem{con}[thm]{Construction}

\graphicspath{{Figures/}}

\renewcommand{\Re}{\textrm{Re}}
\renewcommand{\Im}{\textrm{Im}}

\def \cA{{\mathcal A}}
\def \cB{{\mathcal B}}
\def \cC{{\mathcal C}}
\def \cD{{\mathcal D}}

\def \cG{{\mathcal G}}

\def \cK{{\mathcal K}}
\def \cL{{\mathcal L}}
\def \cM{{\mathcal M}}
\def \cN{{\mathcal N}}

\def \cR{{\mathcal R}}
\def \cS{{\mathcal S}}

\def \cU{{\mathcal U}}
\def \cV{{\mathcal V}}
\def \cW{{\mathcal W}}

\def \cZ{{\mathcal Z}}

\def \B{\mathbb{B}}
\def \C{\mathbb{C}}
\def \D{\mathbb{D}}
\def \H{\mathbb{H}}

\def \N{\mathbb{N}}
\def \P{\mathbb{P}}

\def \R{\mathbb{R}}
\def \S{\mathbb{S}}

\begin{document}
\maketitle

\begin{abstract}
In this work we study the space $S(n)$ of positively oriented, simple $n$-gons with labeled vertices up to oriented similarity and $\cM(n)$ the moduli space of simple $n$-gons as a quotient of $S(n)$. We give a local description of $S(n)$ and $\cM(n)$ around the regular $n$-gon and describe completely the two spaces when $n = 4$ using a topological Morse function. We also provide an asymptotic description of the compactifications $\overline{S(4)}$ and $\overline{\cM(4)}$ up to a complicated set at their boundaries.
\end{abstract}


\section{Introduction}

Spaces of polygons have been widely studied since they capture much of the technical difficulties in different areas of mathematics. Examples of this are the difference between fine and coarse moduli spaces and its generalization to moduli stacks in algebraic geometry \cite{BEH}, the existence or non-existence of global regularizing flows in differential geometry \cite{ChowGli}, \cite{SmBrouFra}, \cite{CoDeRo}, geometric constructions on spaces of polyhedrons \cite{THU} and \cite{fill}, compactifications of other Moduli spaces \cite{HauKnu}, \cite{KaMil2}. In addition to this, the study of spaces of polygons have consecuences on physics and engineering, with applications such as robotic motion planing and pattern recognition \cite{Loz}, \cite{Gott}. The study of the geometry and topology of spaces of polygons with fixed properties is interesting by itself. Examples of this are the study of spaces of polygons with fixed side lengths \cite{KaMil1}, \cite{KaMil2}, fixed angles \cite{BaGhys} and star shaped polygons \cite{KoYam}. 

\noindent
In this work we study the space of simple $n$-gons, with labeled vertices, modulo oriented similarity. More precisely, the set of simple $n$-gons is naturally embedded in the complex projective space $\C\P^{n-2}$, having two homeomorphic connected components, corresponding to positively and negatively oriented $n$-gons respectively, see Section \ref{preliminaries}. The connected component of positively oriented polygons is denoted here as $S(n)$ and as an open subset of $\C\P^{n-2}$, $S(n)$ has a structure of an analytic manifold. We propose the following conjecture

\begin{conj}\label{conj_1}
$S(n)$ is homeomorphic to $\R^{2n-4}$.
\end{conj}

\noindent
There is a parallelism between the space $S(n)$ and the concept of Teichm\"uller space of Riemann surfaces, and between its symmetry group $D_n$ (isomorphic to the dihedral group) and the modular group of the corresponding surface, see Remark \ref{remark_teich_parallelism}. Motivated by this, we refer to the quotient $\cM(n) = S(n)/D_n$ as the moduli space of simple polygons. If $C(X) = \big(X \times (0,1]\big) / ( X \times \{1\})$ denotes the open cone over a topological space $X$, then by linearizing the action of the dihedral group around the regular polygon (i.e. the polygon whose sides and angles are all equal), we can recover the local result

\begin{thm}\label{thm_principal_local}
There is a neighborhood $U \subset \cM(n)$ around the regular $n$-gon such that $U$ is homeomorphic to $C(\S^{2n-5}/\langle \Sigma,T\rangle)$, where $$(z_1,z_2,\dots,z_{n-2})\overset{\Sigma}\longrightarrow(e^{4\pi i/n}z_1,e^{6\pi i/n}z_2,\dots,e^{2 \pi i(n-1)/n}z_{n-2}),\quad\text{and}$$ $$(z_1,z_2,\dots,z_{n-2})\overset{T}\longrightarrow-(e^{-4\pi i/n}\overline{z}_1,e^{-6\pi i/n}\overline{z}_2,\dots,e^{-2 \pi i(n-1)/n}\overline{z}_{n-2}).$$
\end{thm}

\noindent
We conjecture that such result holds globally.

\begin{conj}\label{conj_2}
$\cM(n)$ is homeomorphic to $C(\S^{2n-5}/\langle \Sigma,T\rangle)$.
\end{conj}

\noindent
\textbf{General strategy.} Conjectures \ref{conj_1} and \ref{conj_2} can be proven independently, but a strategy to prove both simultaneously is to find a differentiable Morse function in $S(n)$ with a unique singular point at the regular polygon and whose level sets are invariant under the action of the symmetry group $D_n$. In fact, we can observe that proving the cone structure on $S(n)$ and constructing a Morse function with a unique singular point are equivalent problems in the differentiable category. Moreover, the existence of such Morse function implies the existence of a regularizing flow on $S(n)$ (its gradient flow), that is, a flow taking each polygon onto the regular one withouth loosing the properties of being simple and withouth collapsing vertices. To the knowledge of the authors, the existence of such regularizing flow remains open besides numerous attempts.

\noindent
In this work, we propose a height function $h : \overline{S(4)} \subset \C\P^2 \rightarrow [0,1]$, which is a Morse function in the topological sense (see \cite{Morse}, \cite{Kui} and \cite{Cant}), which has a unique maximum $h = 1$ at the regular quadrilateral (square) and whose level sets are homeomorphic to $\S^3$ and are invariant under the action of the dihedral group, leading to the following

\begin{thm}\label{thm_principal_one}
There is a homeomorphism between $S(4)$ and the open cone  $C(\S^3) \cong \R^4$, such that the height function $h$ corresponds to the function $C(\S^3) \rightarrow (0,1]$, induced from the projection $\S^3 \times (0,1] \rightarrow (0,1]$.
\end{thm}

\noindent
In \cite{AHT-JO} it was already proven that $S(4)$ is homeomorphic to $\R^4$ and a description of the boundary $\partial S(4) \subset \C\P^2$ was given by gluing a polyhedron, making evident that the closure $\overline{S(4)}$ is complicated, in particular it is not a manifold. The advantage of the approach given here, is that the cone structure on $S(4)$ is invariant under the action of the dihedral group, leading to the following result which was out of reach with the techniques of \cite{AHT-JO}.

\begin{cor}
Conjectures \ref{conj_1} and \ref{conj_2} are true for $n=4$.
\end{cor}

\noindent
\textbf{Difficulties.} Here, we are calling $h$ a ``Morse function'', however, it is not easy to prove this statement. In fact, it seems to the authors that it is equally difficult to prove that the height function $h$ is a Morse function than to prove Theorem \ref{thm_principal_one} directly and then conclude that $h$ is a Morse function from there. The most problematic regions are the neighborhoods of the square, where there are a lot of quadrilaterals with different symmetries and where the height function loses differentiability (a phenomenon present in other constructions were symmetries cause problems, such as in \cite{BEH}). To prove Theorem \ref{thm_principal_one}, we may attempt to give a family of homeomorphisms $f_s : h^{-1}(s) \rightarrow \S^3$ which vary continuously on the parameter $s$, and then construct the homeomorphism $S(4) \rightarrow C(\S^3)$ as $Z \mapsto [f_{h(Z)}(Z), h(Z)]$. This strategy is hard to carry out and it seems to the authors that the construction of such a continuous family of homeomorphisms equals the difficulty of the construction of a regularizing flow, which theoretically exists precisely by Theorem \ref{thm_principal_one}. We overcome this difficulty by splitting $S(4)$ into subsets where the height function $h$ behaves well (see Remark \ref{remark_differentiability_height} for instance) and where we are able to construct the continuous family of homeomorphisms, thus obtaining a cone structure on every such set. With the previous partition, the problem translates into understanding the combinatorics of the intersections of the sets, which we give explicitly. The last main difficulty we have to face is on the compactifications $\overline{S(4)}$ and $\overline{\cM(4)}$, where there are some bad convergence phenomena. We face this difficulty by identifying the problematic region which gives bad behaviour at the boundary and prove that we have a very good asymptotic behaviour outside of it

\begin{thm}\label{thm_principal_two}
There exists a subspace $\cD \subset \overline{S(4)}$ such that $\overline{S(4)} \smallsetminus \cD \cong \H^4 \sqcup \H^4$, where $\H^4$ denotes the fourth-dimensional, closed upper half-space. The space $\cD$ is an open cone over the two sphere $\cD \cap S(4) \cong C(\S^2)$, compactified with a space $\cD \cap \partial S(4)$, which is obtained as a union of a wedge of spheres $\S^2 \vee \S^2$ and a closed disc $\D^2$.
\end{thm}

\noindent
As these constructions are invariant under the action of the dihedral group, we obtain the asymptotic description of the Moduli space

\begin{thm}\label{thm_principal_three}
There is a subspace $\cD' \subset \overline{\cM(4)} = \overline{S(4)} / D_4$, that is homeomorphic to an open cone over the $2$-disc $\cD' \cap \cM(4) \cong C(\D^2)$, compactified with a topological disc $\cD' \cap \partial \cM(4) \cong \D^2$, so that the complement $\overline{\cM(4)} \smallsetminus \cD'$ is homeomorphic to $\H^4$.
\end{thm}

\noindent
As the subset $\cD'$ in Theorem \ref{thm_principal_three} is capturing the bad asymptotic behaviour, the compactification of the open cone $\cD' \cap \cM(4)$ with the topological disc $\cD' \cap \partial \cM(4)$ is complicated, see Remark \ref{remark_compactification_cone_over_disc}.

\noindent
It is worth noting that there have been previous uses of Morse functions on spaces of polygons, by considering parameters such as side lengths \cite{Haus} or diagonals \cite{KaMil1}. Those functions don't work here, because they don't capture the asymptotic behaviour to the boundary $\partial S(n)$. There is however a Morse function on Teichm\"uller spaces of Riemann surfaces whose behaviour resembles much more to the one presented here, which is the systole \cite{Sch}.

\noindent
\textbf{Structure of the paper.} In Section \ref{preliminaries}, we define precisely the spaces $S(n)$ and $\cM(n)$, give the examples for $n = 3$, and conclude with an index of notations. In Section \ref{features}, we give the local description of $\cM(n)$, prove Theorem \ref{thm_principal_local}, and in Example \ref{eje4} we give a local, linearized version of Theorem \ref{thm_principal_one}. In Section \ref{height}, we introduce the height function $h$ for quadrilaterals and give some basic properties, we also give the example of an analogous height function on $S(3)$ and explain how it determines its topology. Section \ref{levelcurves} is devoted to give a good decomposition of $S(4)$ into subsets on which the height function behaves well and describe the combinatorics of their intersections; we finish this section by proving that the level sets of $h$ in $S(4)$ are homeomorphic to the sphere $\S^3$. In Section \ref{special}, we prove Theorem \ref{thm_principal_one} and in the process, we prove Conjectures \ref{conj_1} and \ref{conj_2}. In Section \ref{section_asymptotic_desc}, we study the compactifications $\overline{S(4)}$ and $\overline{\cM(4)}$, and as it is the most technical section of the paper, we split it into four subsections. In Subsection \ref{subsection_bad_asymptotic} we describe the set $\cD$ mentioned in Theorem \ref{thm_principal_two}, which is the set of quadrilaterals presenting bad asymptotic phenomena at the boundary. In Subsection \ref{subsection_good_behaviour}, we give a decomposition of the boundary $\partial S(4) \smallsetminus \cD$, with equivalent combinatorics of the decompositions of every other level set $h^{-1}(s) \smallsetminus \cD$. In Subsection \ref{subsection_good_asymptotic}, we use all the previous constructions to prove Theorem \ref{thm_principal_two}. In Subsection \ref{subsection_compactification_moduli} we study the decomposition in $\overline{\cM(4)}$ induced by the previous construction and prove Theorem \ref{thm_principal_three}. We present in Section \ref{section:open_questions} a list of concluding remarks and open questions. Finally, we add two appendices: Appendix \ref{section_gluing_balls}, were we give two technical theorems on gluing two topological closed balls and two topological closed half-spaces, and Appendix \ref{section_intermediate_quotients}, were we give a homeomorphism between two quotients, one with a group and the other as a double quotient by a normal subgroup.

\subsection*{Acknowledgements} We would like to thank Gerardo Arizmendi Echegaray for fruitful conversations on an early stage of the project, as well as for suggesting a modification to the height function which decreased the complexity of the combinatorics in the decomposition of its level sets, c.f. Definition \ref{conos}. Manuel Sedano would like to thank CONACyT grant CB2016-283960 as well as DGAPA-UNAM Postdoctoral Scholarship. 

\section{Preliminaries}
\label{preliminaries}

A polygon with $n$ sides (or simply an $n$-gon) in the complex plane with consecutive vertices at $z_1,z_2,\dots,z_n \in \C$, is completely determined by the point $(z_1,z_2,\dots,z_n)\in\C^n$. In this way, the set of $n$-gons in $\C$ with labeled vertices is equipped with the topology of $\C^n$. For all $Z=(z_1,\dots,z_n)\in\C^n$, $\mathfrak c(Z)$ denotes the polygonal curve $\overline{z_1z_2}\cup\overline{z_2z_3}\cup\cdots\cup\overline{z_{n-1}z_n}\cup\overline{z_nz_1}$, where $\overline{ab}\subset\C$ is the segment from $a$ to $b$. 

\begin{defi}
\label{simple}
Let $Z\in\C^n$ be a $n$-gon.
\begin{enumerate}[(a)]
\item We say that $Z$ is simple if its vertices are different and $\mathfrak c(Z)$ is a Jordan's curve. $\S(n)\subset\C^n$ denotes the set of simple $n$-gons.
\item If $Z$ is simple, then we say that $Z$ is positively (negatively) oriented if the increasing order in the vertices of $\mathfrak{c}(Z)$ is in counter-clockwise (clockwise). 
\end{enumerate}
\end{defi}
In \cite{AHT} it is proved that $\S(n)$ is an open subset and has two homeomorphic connected components, corresponding to positively and negatively oriented $n$-gons. The group of complex affine transformations on $\C$, denoted here as $\cA_{\C}$, acts diagonally on $\C^n$ as
\[  (az+b,(z_1,\dots,z_n))\mapsto(az_1+b,\dots,az_n+b), \]
preserving angles, labeled vertices and proportional side lengths in the $n$-gons. The quotient $P(n):=\left(\C^n\smallsetminus\{(z,z,\dots,z)\}\right)/\cA_{\C}$ is the \textit{space of shapes} of $n$-gons with labeled vertices up to oriented similarity as considered in \cite{THU}. Notice that two equivalent $n$-gons in the linear subspace $V_n=\{(0,z_2,\dots,z_n)\in\C^n\}$ differ by just a complex factor, so the space $P(n)$ is a complex projective space 
    \[  P(n) \cong \P_{\C}V_n\cong\C\P^{n-2}.  \]
The action of $\cA_{\C}$ preserves orientation and the property of being simple in $n$-gons, so if $\eta\colon\C^n\smallsetminus\{(z,z,\dots,z)\}\to P(n)$ denotes the natural projection, then $\eta\left(\S(n)\right)\subset P(n)$ is open and has two connected components, corresponding to positively and negatively oriented $n$-gons. Denote by $S(n) \subset \eta(\S(n))$ the connected component of positively oriented shapes of $n$-gons and denote by $\eta(z_1,z_2,\dots,z_n) = [z_1,z_2,\dots,z_n] \in S(n)$ the corresponding class of an $n$-gon.

The linear isomorphism $\widehat{\sigma}(z_1,z_2,\dots,z_n)=(z_2,\dots,z_n,z_1)$, re-enumerates the vertices of an $n$-gon in a cyclic way. Notice that two $n$-gons $Z,W\in\C^n$ have the same shape if and only if $\widehat{\sigma}(Z)$ and $\widehat{\sigma}(W)$ have the same shape, so there exists a function $\sigma\colon P(n)\to P(n)$ such that $\sigma\circ\eta=\eta\circ\widehat{\sigma}$.

\begin{obs}
\label{chart}
We may observe that, as the diagonal $\{(z,z,\dots,z)\}$ is not considered in $P(n)$, at least two vertices of a polygon in $P(n)$ must be different. So if we consider the set $U_{12}=\{[z_1,z_2,\dots,z_n]\in P(n)\colon z_1\neq z_2\}$ and the coordinate chart $\phi_{12}\colon U_{12}\to\C^{n-2}$, given by
    \[  \phi_{12}\left([z_1,z_2,\dots,z_n]\right)=\left(\frac{z_3-z_1}{z_2-z_1},\dots,\frac{z_n-z_1}{z_2-z_1}\right),   \]
then for every $Z \in P(n)$, $\sigma^k(Z) \in U_{12}$, for some $k \in \N$. That is, up to re-enumeration of the vertices, $P(n)$ is contained in $U_{12}$. Through the paper, we will mostly work with representatives $[0,1,z_3,\dots,z_n] \in U_{12}$. 
\end{obs}

\begin{ejem}[Triangles]
\label{eje1} 
The space $P(3)$ is biholomorphic to the Riemann sphere. There are three options for a triangle $[0,1,x+iy]\in U_{12}$:
    \[  (a)~~y>0,\qquad\quad (b)~~y<0,\qquad\text{~~and~~} \qquad ~(c)~~y=0.  \]
In the plane $\phi_{12}(U_{12})$, the case $(a)$ is the upper half plane and corresponds to the set $S(3)$ of positively oriented triangles; the case $(b)$ is the lower half plane and corresponds to the set of negatively oriented triangles; and the case $(c)$, which corresponds to triangles degenerated to a segment, compactifies to a circle in $P(3)$ with the degenerated triangle $[0,0,1]$. Thus, the set of degenerated triangles is a circle that splits the sphere $P(3)$ into two discs, corresponding to positively and negatively oriented triangles.
\end{ejem}

The involution $\tau\colon U_{12}\to U_{12}$, which reflects the last $(n-2)$-vertices trough the line $\Re(z)=\frac{1}{2}$, changes the order in which the interior angles are presented but without changing the orientation and can be written as 
    \[  \tau([0,1,z_3,\dots,z_n])=[0,1,1-\bar{z}_n,\dots,1-\bar{z}_3].   \]
Since $U_{12}\subset P(n)$ is dense (Remark \ref{chart}), then $\tau$ can be continuously extended to an involution of $P(n)$.

\begin{prop}
\label{diffeo}
The functions $\sigma,\tau\colon P(n)\to P(n)$ are diffeomorphisms and the group generated by them $\langle\sigma,\tau\rangle$ is isomorphic to the dihedral group of order $2n$.
\end{prop}
\noindent
\begin{proof}
As the funciton $\sigma$ is defined by the property $\sigma\circ\eta=\eta\circ\widehat{\sigma}$, with $\widehat{\sigma}(z_1,z_2,\dots,z_n)=(z_2,\dots,z_n,z_1)$ and $\eta\colon\C^n\smallsetminus\{(z,z,\dots,z)\}\to P(n)$ the natural projection, the differentiability of $\sigma$ follows from the differentiability of $\widehat{\sigma}$. Analogously, if we define $\widehat{\tau}\colon\C^n\to\C^n$ as $\widehat{\tau}(z_1,\dots,z_n) = -(\overline{z_2},\overline{z_1},\overline{z_n},\dots,\overline{z_3})$, then $\sigma\circ\tau=\eta\circ\widehat{\tau}$ and the differentiability of $\tau$ follows from the differentiability of $\widehat{\tau}$. Finally, it is straightforward to check that $\sigma^n=\tau^2=(\tau\sigma)^2=1$. Then $\langle\sigma,\tau\rangle$ is isomorphic to the dihedral group of order $2n$.
\end{proof}

\noindent
Just as the action of the group $\cA_{\C}$ in $\C^n$ eliminates translations, rotations and homotheties in $n$-gons, the action of the dihedral group $\langle\sigma,\tau\rangle$ in $P(n)$, removes the labels at the vertices and identifies $n$-gons equivalent under reflections.

\begin{defi}
\label{moduli}
Let $\overline{S(n)} \subset P(n)$ denote the closure of $S(n)$. Then
\begin{enumerate}[(a)]
    \item the quotient $\cM(n)=S(n)/\langle\sigma,\tau\rangle$ is the Moduli Space of $n$-gons.
    \item the quotient $\overline{\cM(n)} = \overline{S(n)} /\langle\sigma,\tau\rangle$ is a compactification of $\cM(n)$.
\end{enumerate}
\end{defi}

\begin{obs}\label{remark_teich_parallelism}
A different description of the space $S(n)$, is to consider it as the set of graph isomorphisms $\mathfrak{R}_n \rightarrow Z$ (where $\mathfrak{R}_n$ is the regular $n$-gon and $Z$ is an arbitrary simple polygon), with the action of the dihedral group acting via precomposition. This approach creates a parallelism between the concept of Teichm\"uller space and $S(n)$, and between the modular group of a surface and the dihedral group as a symmetry group of $S(n)$ (see \cite{bestvina} for this and a similar construction for the group of outer automorphisms $Out(F_n)$). Moreover, by construction, the action of the dihedral group extends to the closure $\overline{S(n)} \subset \C\P^{n-2}$, which is a well behaved property for this parallelism and a good justification to consider $\overline{S(n)} / \langle \sigma, \tau \rangle$ as the compactification of $\cM(n)$, see \cite{Okoshita}.
\end{obs}

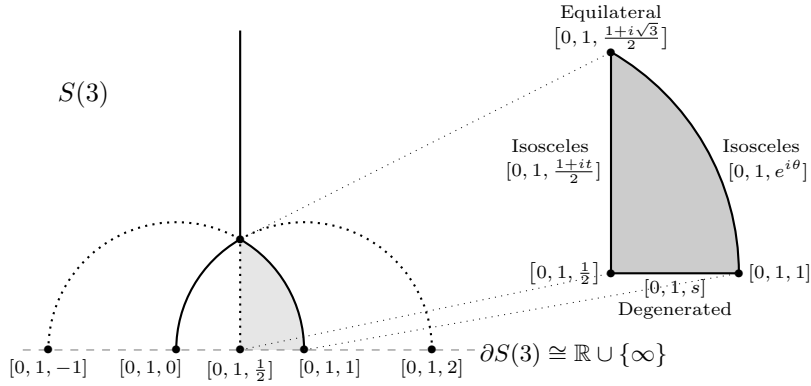
\begin{figure}[h]
\begin{center}
\begin{tikzpicture}[scale=1.7]
\fill [gray!20] (-1.5,0)--(-1,0)--(-1.5,.866)--(-1,0) arc (0:60:1);
\draw [black!50, dashed] (-3.2,0)--(.35,0);
\draw [thick=0.4] (-1,0) arc (0:60:1);
\draw [thick=0.4, dotted] (-1.5,.866) arc (60:180:1);

\draw [thick=0.4, dotted] (0,0) arc (0:120:1);
\draw [thick=0.4] (-1.5,.866) arc (120:180:1);

\draw [thick=0.4] (-1.5,.866)--(-1.5,2.5);
\draw [thick=0.4, dotted] (-1.5,0)--(-1.5,.866);

\draw (.3,-0.04) node [right] {\small{$\partial S(3)\cong\R\cup\{\infty\}$}};
\draw (-3,2) node [right] {$S(3)$};

\draw (-3,0) node [below] {\scriptsize{$[0,1,-1]$}};
\draw (-3,0) node {\scriptsize{$\bullet$}};
\draw (-2.2,0) node [below] {\scriptsize{$[0,1,0]$}};
\draw (-2,0) node {\scriptsize{$\bullet$}};
\draw (-.8,0) node [below] {\scriptsize{$[0,1,1]$}};
\draw (-1,0) node {\scriptsize{$\bullet$}};
\draw (0,0) node [below] {\scriptsize{$[0,1,2]$}};
\draw (0,0) node {\scriptsize{$\bullet$}};

\draw (-1.5,0) node {\scriptsize{$\bullet$}};
\draw (-1.5,0) node [below] {\scriptsize{$[0,1,\frac{1}{2}]$}};

\draw (-1.5,.866) node {\scriptsize{$\bullet$}};

\fill [gray!40] (1.4,.6)--(2.4,.6)--(1.4,2.33)--(2.4,.6) arc (0:60:2);
\draw [thick=0.3] (1.4,.6)--(2.4,.6);
\draw [thick=0.3] (1.4,.6)--(1.4,2.33);
\draw [thick=0.3] (2.4,.6) arc (0:60:2);

\draw (1.4,.6) node {\scriptsize{$\bullet$}};
\draw (1.4,0.6) node [left] {\scriptsize{$\big[0,1,\frac{1}{2}\big]$}};
\draw (2.4,0.6) node {\scriptsize{$\bullet$}};
\draw (2.4,0.6) node [right] {\scriptsize{$[0,1,1]$}};

\draw (1.4,2.33) node {\scriptsize{$\bullet$}};
\draw (1.4,2.5) node [above] {\scriptsize{Equilateral}};
\draw (1.4,2.28) node[above] {\scriptsize{$\big[0,1,\frac{1+i\sqrt{3}}{2}\big]$}};

\draw (1.3,1.6) node [left] {\scriptsize{Isosceles}};
\draw (1.4,1.4) node [left] {\scriptsize{$[0,1,\frac{1+it}{2}]$}};

\draw (2.2,1.6) node [right] {\scriptsize{Isosceles}};
\draw (2.24,1.4) node [right] {\scriptsize{$[0,1,e^{i\theta}]$}};

\draw (1.9,0.63) node [below] {\scriptsize{$[0,1,s]$}};
\draw (1.9,0.45) node [below] {\scriptsize{Degenerated}};

\draw [black!90, dotted] (-1,0)--(2.4,0.6);
\draw [black!90, dotted] (-1.5,0)--(1.4,0.6);
\draw [black!90, dotted] (-1.5,.866)--(1.4,2.33);

\end{tikzpicture}
\caption{Different types of triangles in the compactification $\overline{\cM(3)}$.}
\label{figure1}
\end{center}
\end{figure}

\begin{ejem}
\label{eje2}
From Example \ref{eje1} we know that $S(3)\subset U_{12}$ is a copy of the upper half plane. In the representatives $[0,1,z_3]$, the functions $\sigma([0,1,z_3])=[0,1,\frac{1}{1-z_3}]$ and  $\tau([0,1,z_3])=[0,1,1-\bar{z}_3]$ are a hyperbolic rotation of angle $\frac{2\pi}{3}$, centered at the equilateral triangle $\mathfrak{R}_3=[0,1,\frac{1+i\sqrt{3}}{2}]$, and a reflection through the line $[0,1,\frac{1}{2}+iy]$, respectively. A fundamental domain of the $\langle\sigma,\tau\rangle$-action in $S(3)$ is$$\cM(3)=\left\{[0,1,z_3]\in U_{12}\colon |z_3|\leq 1,Im(z_3)>0,Re(z_3)\geq\frac{1}{2}\right\}.$$The compactification $\overline{\cM(3)}$ is obtained by adding to $\cM(3)$ the degenerated triangles $[0,1,s]$ with $s\in [\frac{1}{2},1]$, see Figure \ref{figure1}.
\end{ejem}

\noindent
From Example \ref{eje2} it follows that $\cM(3)$ is an hyperbolic orbifold of dimension two, with boundary $\partial\overline{\cM(3)}$ at infinite distance with respect to this metric. This is the only case where the moduli space $\cM(n)$ is a geometric quotient.

\subsection*{Index of notations}
We collect here some of the notations that will be used throughout the paper.
\begin{enumerate}
    \item Open cone over a topological space $X$: $C(X) = \big(X \times (0,1]\big) / \{\sim\}$, where $(x,1) \sim (x',1)$, for all $x \in X$.

    \item Suspension over a topological space $X$: $\Sigma(A) = \big(X \times [-1,1]\big) / {\sim}$, where $(x,1) \sim (x',1)$ or $(x,-1) \sim (x',-1)$, for all $x \in X$.
    
    \item Closed half spaces $\H^n = \{(x_1,\cdots , x_n) \in \R^n : x_n \geq 0\}$

    \item Closed ball $\D^n = \{(x_1,\cdots , x_n) \in \R^n : x_1^2 + \cdots + x_n^2 = 1\}$,

    \item The distance $d(z,A) = \min\{|z-a| : z \in A\}$, from $z$ to a subset $A \subset \C$.

    \item Segment $\overline{z_1 z_2} = \{(1-t)z_1 + tz_2 : 0 \leq t \leq 1\} \subset \C$, between two points $z_1,z_2 \in \C$. Not to be confused with conjugation $\overline{x + i y} = x - iy$, as we won't ever use complex conjugation on a multiplication of two complex numbers.

    \item Polygonal curve $\mathfrak{c}(Z) = \overline{z_1z_2}\cup\overline{z_2z_3}\cup\cdots\cup\overline{z_{n-1}z_n}\cup\overline{z_nz_1}$ defined by a polygon $Z = (z_1, \cdots, z_n) \in \C^n$.

    \item Polygonal curve $\mathcal{L}_j(Z) = \overline{z_{j+1}z_{j+2}}\cup\overline{z_{j+2}z_{j+3}} \cup \cdots \overline{z_{j-2}z_{j-1}}$, opposite to a vertex $z_j$ of a polygon $Z = (z_1, \cdots, z_n) \in \C^n$. We take the cyclic convention on indices mod $n$.

    \item The shape $[z_1, \cdots ,z_n] \in S(n)$ as the class of the polygon $(z_1, \cdots ,z_n) \in \C^n$.
    
    \item The regular polygon $\mathfrak{R}_n=[1,e^{2i\pi/n},\dots,e^{2i\pi(n-1)/n}]\in S(n)$. We call $\mathfrak{R}_4$ ``the square'' as well.

    \item The closure $\overline{A} \subset \C\P^{n-2}$ of a subset $A \subset S(n) \subset \C\P^{n-2}$.
    
    \item The restricted closure $\textrm{cl}(A) = \overline{A} \cap S(n)$ of a subset $A \subset S(n) \subset \C\P^{n-2}$.

    \item  We will use the word ``side'' por the edges of a polygon, including bigons (c.f. Definici\'on \ref{bigonization-def}). On the other hand we will use the word ``edge'' for the one-dimensional sides of polyhedrons the corresponding graphs (such as in Figure \ref{Diagrama_pegados} and Figure \ref{pegado_frontera}).
\end{enumerate}

\noindent
We clarify on two more concepts widely used in the text:
\begin{obs}\label{remark_clarifying_I_bundle}
If we consider the trivial $I$-bundle $\D^2 \times \D^1$, which is homeomorphic to $\D^3$ and collapse some fibers on a closed interval of the boundary $I \subset \partial \D^2 \cong \S^1$: $(x,t) \sim (x,t')$ for all $t \in \D^1$ and $x \in I$, then the quotient 
    \[     (\D^2 \times \D^1) / {\sim} \]
is still homeomorphic to $\D^3$. We will call this kind of quotients an ``$I$-bundle with collapsed fibers''. We can further observe that the we can collapse fibers in a series of disjoint intervals $I_1 \sqcup \cdots \sqcup I_k \subset \partial \D^2$ with analogous behaviour.
\end{obs}

\begin{obs}\label{remark_different_omega_notations}
We will mostly use Remark \ref{remark_clarifying_I_bundle} on $I$-bundles with collapsed fibers, where the bases space is a closed subspace of $\{z \in \C : |z| \leq 1, \ \Im(z) \geq 0\}$ homeomorphic to $\D^2$. We will be using mostly the notation $\Omega, \Omega', \Omega'',...$ with some added sub-indices. We don't have a unified notation for this as we will be using different subsets for different constructions and the precise description of the domain is subtle and necessary. We will define each domain where needed.
\end{obs}

\section{Local description of $\cM(n)$}
\label{features}

In this section we describe the topology of $\cM(n)$ around the regular $n$-gon and explain in detail the construction for triangles and quadrilaterals.

Let $\mathfrak{R}_n=[1,e^{2i\pi/n},\dots,e^{2i\pi(n-1)/n}]\in S(n)$ be the shape of the regular $n$-gon. It is straightforward to check that $\mathfrak{R}_n=\tau(\mathfrak{R}_n)=\sigma(\mathfrak{R}_n)$, and in fact, $\mathfrak{R}_n$ is the only element in $S(n)$ that is fixed by $\sigma$.

\smallskip
\begin{lem}
\label{diferencial}
Let $D_{\mathfrak{R}_n}\sigma$ and $D_{\mathfrak{R}_n} \tau$ be the differentials of $\sigma$ and $\tau$ at $\mathfrak{R}_n$. Then there is a basis $\{\zeta_2,\dots, \zeta_{n-1} \}$ of the complex vector space $T_{\mathfrak{R}_n} S(n)$ such that$$D_{\mathfrak{R}_n}\sigma\Big(\!\sum_{k=2}^{n-1} s_k\zeta_k\!\Big)\!=\!\sum_{k=2}^{n-1} e^{2\pi ik/n}s_k\zeta_k,~\text{and}~~D_{\mathfrak{R}_n}\tau\Big(\!\sum_{k=2}^{n-1}s_k\zeta_k\!\Big)\!=\!\sum_{k=2}^{n-1}-e^{-2 \pi ik/n}\overline{s_k} \zeta_k\,.$$
\end{lem}

\begin{proof}
Let $B(Z,W)=\sum_j z_j\overline{w}_j$ denotes the standard Hermitian product in $\C^n$. Consider the basis $\beta=\big\{\widehat{\zeta}_k=(1,e^{2\pi ik/n},(e^{2\pi ik/n})^2,\dots,(e^{2\pi ik/n})^{n-1})\big\}$ with $k=0,1,\dots,n-1$, of $\C^n$. The products between these vectors are
    \[  B(\widehat{\zeta}_k,\widehat{\zeta}_l) = \sum_{m=0}^{n-1}\left(e^{\frac{2\pi i(k-l)}{n}}\right)^m=\left\{\begin{array}{lcr} n & \text{if} & k=l \\ 0 & \text{if} & k\neq l \end{array}\right.\]
and therefore, $\beta$ is an orthogonal basis. Let $V=\{(z_1,z_2,\dots,z_n)\in\C^n\colon z_1+z_2+\cdots+z_n=0\}$ be the complex vector subspace of $n$-gons that have center of mass at the origin. Note that $V=\langle\widehat{\zeta}_1,\widehat{\zeta}_2,\dots,\widehat{\zeta}_{n-1}\rangle$, $\C^n=\langle\widehat{\zeta}_0\rangle\oplus V$, $\eta(V)=\P_{\C}V=P(n)$ and $\eta(\widehat{\zeta}_1)=\mathfrak{R}_n$. Since $\langle\widehat{\zeta}_1\rangle$ is the kernel of the differential $D_{\widehat{\zeta}_1}\eta\colon V\to T_{\mathfrak{R}_n}P(n)$, then the restriction$$D_{\widehat{\zeta}_1}\eta\colon\langle\widehat{\zeta_2},\cdots,\widehat{\zeta}_{n-1}\rangle\to T_{\mathfrak{R}_n}S(n)$$is a linear isomorphism.

Let $\widehat{\tau}\colon\C^n\to\C^n$ be the function $(z_1,\dots,z_n)\mapsto-(\overline{z_2},\overline{z_1},\overline{z_n},\dots,\overline{z_3})$. Since $\widehat{\sigma}(\widehat{\zeta}_k)=e^{2ik\pi/n}\widehat{\zeta}_k$ and $\widehat{\tau}(\widehat{\zeta}_k)=-e^{-2ik\pi/n}\widehat{\zeta}_k$, then $V$ and $\langle\widehat{\zeta_2}, \dots,\widehat{\zeta}_{n-1}\rangle$ are invariant under $\widehat{\sigma}$ and $\widehat{\tau}$, moreover, $\widehat{\sigma}$ and $\widehat{\tau}$ satisfy $\eta\circ\widehat{\sigma}=\sigma\circ\eta$ and $\eta\circ\widehat{\tau}=\tau\circ\eta$. Differentiating in $\widehat{\zeta}_1$, we obtain$$D_{\widehat{\zeta}_1}\eta\circ\widehat{\sigma}=D_{\mathfrak{R}_n}\sigma\circ D_{\widehat{\zeta}_1}\eta\,,\quad\text{and}\quad D_{\widehat{\zeta}_1}\eta\circ\widehat{\tau}=D_{\mathfrak{R}_n}\tau\circ D_{\widehat{\zeta}_1}\eta\,.$$The proof is finished by taking $\zeta_k=D_{\widehat{\zeta}_1}\eta(\widehat{\zeta_k})$, for $k=2,3,\dots,n-1$.
\end{proof}

\noindent
Let $\Sigma,T\colon\C^{n-2}\to\C^{n-2}$ denote the orthogonal transformations:
    \[  (z_1,z_2,\dots,z_{n-2})\overset{\Sigma}\longrightarrow(e^{4\pi i/n}z_1,e^{6\pi i/n}z_2,\dots,e^{2 \pi i(n-1)/n}z_{n-2}),\quad\text{and} \]
    \[  (z_1,z_2,\dots,z_{n-2})\overset{T}\longrightarrow-(e^{-4\pi i/n}\overline{z}_1,e^{-6\pi i/n}\overline{z}_2,\dots,e^{-2 \pi i(n-1)/n}\overline{z}_{n-2}).    \]
The projection of the regular $n$-gon into $\cM(n)$ is also denoted by $\mathfrak{R}_n$. The following result provides a local topological description of $\cM(n)$ around $\mathfrak{R}_n$.
\begin{thm}
\label{local}
Local neighbourhoods of $\mathfrak{R}_n\in\cM(n)$ are homeomorphic to the cone over the quotient $\S^{2n-5}/\langle \Sigma,T\rangle$. 
\end{thm}

\begin{proof}
Since the maps $\widehat{\sigma},\widehat{\tau}\colon V\to V$ are orthogonal with respect to the Euclidean product, then the projections $\sigma,\tau\colon P(n)\to P(n)$ are Riemannian isometries of the Fubini-Study metric in $P(n)$. The exponential map around $\mathfrak{R}_n$, locally conjugates the actions of $\sigma$ and $\tau$ in $S(n)$ with their differential maps at $T_{\mathfrak{R}_n} S(n)$. Then, around $\mathfrak{R}_n$, the moduli space $\cM(n)=S(n)/\langle\sigma,\tau\rangle$ is homeomorphic to the cone over the quotient of the unit sphere in $T_{\mathfrak{R}_n} S(n)$ with the action of the dihedral group $\langle D_{\mathfrak{R}_n}\sigma,D_{\mathfrak{R}_n}\tau\rangle$. From Lemma \ref{diferencial}, it follows that the actions of $\langle D_{\mathfrak{R}_n}\sigma,D_{\mathfrak{R}_n}\tau\rangle$ in $T_{\mathfrak{R}_n} S(n)$ and $\langle\Sigma,T\rangle$ in $\C^{n-2}$ are conjugated.
\end{proof}

\begin{ejem}
\label{eje3}
For $n=3$, $\Sigma(z)=e^{4\pi i/3}z$ and $T(z)=-e^{-4\pi i/3}\overline{z}$. Since $T$ is a reflection through the line with argument $-\frac{\pi}{6}$, then a fundamental domain in the circle $\S^1$ for the action of $\langle \Sigma,T\rangle$, is the arc $\{e^{i\theta}\colon -\frac{\pi}{6}\leq\theta\leq 0\}$. In Example \ref{eje2} we showed that the compactification $\overline{\cM(3)}$ is homeomorphic to the closed cone over a copy of this arc, this is an exceptional case, as the higher dimensional compactifications $\overline{\cM(n)}$ will be much more involved than a closed cone over the quotient $\S^{2n-5}/\langle \Sigma,T\rangle$.
\end{ejem}

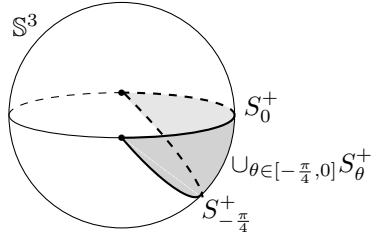
\begin{figure}[h]
\begin{center}
\begin{tikzpicture}[scale=0.75]
\fill [gray!20] (-6,-0.4) arc (-90:0:2 and 0.4) -- (-4,0) arc (0:90:2 and 0.4) -- (-6,0.4) [rotate=-49] arc (-90:0:1.8 and 0.4);

\fill [gray!35] (-6,-0.4) arc (-90:0:1.98 and 0.4) -- (-4,0) arc (0:-46:2) -- (-6,-0.4) [rotate=-50] arc (-90:0:1.7 and 0.4);

\draw (-6,0) circle (2);

\draw (-8,0) arc (180:270:2 and 0.4);

\draw [thick=0.6] (-6,-0.4) [rotate=-50] arc (-90:11:1.68 and 0.4);

\draw [dashed] (-8,0) arc (-180:-270:2 and 0.4);

\draw [thick=0.6, dashed] (-4,0) arc (0:90:2 and 0.4);

\draw [thick=0.6, dashed] (-5.95,0.38) [rotate=-45] arc (-270:-349:2.32 and 0.4);

\draw [thick=0.6] (-4,0) arc (0:-90:2 and 0.4);


\draw [black] (-6,0.4) node {\tiny{$\bullet$}};
\draw [black] (-6,-0.4) node {\tiny{$\bullet$}};

\draw (-4,0.15) node [right] {$S_0^+$};
\draw (-4.75,-1.75) node [right] {$S_{-\frac{\pi}{4}}^+$};
\draw (-7.7,1.6) node { $\S^3$};
\draw (-2.8,-0.9) node { $\cup_{\theta\in[-\frac{\pi}{4},0]} S_\theta^+$};

\end{tikzpicture}
\caption{A schematic representation of the fundamental domain for the $\langle\Sigma,T\rangle$-action on the sphere $\S^3$. }
\label{figure2}
\end{center}
\end{figure}

\begin{ejem}
\label{eje4}
For $n=4$, the quotient $\S^3/\langle\Sigma, T\rangle$ is homeomorphic to $\S^3$. The homeomorphism can be described as follows: For $\theta\in[0,2\pi]$, the set $S^+_\theta:=\{(z_1,te^{i\theta})\in\S^3\colon t\geq 0\}$ is a copy of the closed upper hemisphere in $\S^2$. Notice that$$\S^3=\bigcup_{\theta\in[0,2\pi]} S^+_\theta\, ,\qquad \cS:=\bigcap_{\theta\in[0,2\pi]} S^+_\theta=\{(z_1,0)\in\S^3\}\cong\S^1$$is the boundary circle of all the hemispheres $S^+_\theta$, and for a generic interval $[\theta_1,\theta_2]$, the union $\cup_{\theta\in[\theta_1,\theta_2]} S_\theta^+$ is a copy of the closed disc $\D^3$ with boundary sphere $S_{\theta_1}^+\cup S_{\theta_2}^+\cong\S^2$ (Figure \ref{figure2} illustrates the union $\cup_{\theta\in[-\frac{\pi}{4},0]} S_\theta^+$). As $\Sigma(z,w) = (-z, e^{3 \pi i /2} w)$ and $T(z,w) = (\overline{z}, e^{- \pi i /2} \overline{w})$, then a fundamental domain for the action of $\langle \Sigma,T\rangle$ in $\S^3$, is the angular region $\cup_{\theta\in[-\frac{\pi}{4},0]} S_\theta^+$, homeomorphic to a closed ball $\B^3$. We can split this domain as its interior $\bigcup_{\theta\in(-\frac{\pi}{4},0)} S_\theta^+ \smallsetminus \cS $, where no two elements are identified with an elemen of $\langle \Sigma , T \rangle$ and its boundary $S_0^+ \cup S_{-\pi/4}^+ \cong \S^2$, where $T$ acts in $S_{-\pi/4}^+$ as $(z,w) \mapsto (\overline{z},w)$, and $T \circ \Sigma$ acts in $S_0^+$ as $(z,w) \mapsto (-\overline{z},w)$. It is not hard to see that taking the quotient of  $\cup_{\theta\in[-\frac{\pi}{4},0]} S_\theta^+$ under the identifications given by $T$ gives us a space homeomorphic to $\B^3$ and then taking the quotient of this space under the identifications given by $T \circ \Sigma$ gives us a space homeomorphic to $\S^3$.
\end{ejem}

\begin{obs}
The constructions for quadrilaterals in Example \ref{eje4}, can be done analogously for all $n$. It is necessary to analyze the quotient of the $(2k-1)$-dimensional waist sphere $S^{2k-1}:=\{(z_1,\dots,z_k,0,\dots,0)\in\S^{2n-5}\}$ by the action of the dihedral group $\langle\Sigma,T\rangle$. In this way we construct a fundamental domain in $\S^{2n-5}$ for the action of $\Sigma$, which is invariant by the action of $T$. Such fundamental domain makes evident that there is a homeomorphism between $\S^{2n-5}/\langle\Sigma,T\rangle$ and the double quotient$$\big\{\S^{2n-5}/\langle\Sigma\rangle\big\}/\langle T\rangle\cong L(n:2,3,\dots,n-1)/\langle T\rangle,$$where $L(n:2,3,\dots,n-1)=\S^{2n-5}/\langle\Sigma\rangle$ is an Orbifold Lens Space. 
\end{obs}

The reason behind the fact that $\S^{2n-5}/\langle\Sigma,T\rangle$ is homeomophic to the double quotient $\big\{\S^{2n-5}/\langle\Sigma\rangle\big\}/\langle T\rangle$, is that the cyclic group $\langle\Sigma\rangle$ is normal in the dihedral group $\langle\Sigma,T\rangle$, see Proposition \ref{prop:A.1}.

\section{Height function in quadrilaterals}
\label{height}

In this section, we define the height function in the space of simple quadrilaterals, which will be used to describe its topology, and give some of its properties. We also describe the analogous function in the case of triangles and use it to describe its topology. 

\begin{defi}
\label{distances}
For every $Z=(z_1,z_2,z_3,z_4)\in\C^4$,
\begin{enumerate}[(a)]
    \item the longest side function is $\ell\colon\C^4\to\R$, $\ell(Z)=\max_{1\leq j\leq 4}|z_j-z_{j+1}|$, and
    \item the minimum distance function is $r\colon\C^4\to\R$, $r(Z)=\min_{1\leq j\leq 4}\{r_j(Z)\}$, where $r_j(Z):=d(z_j,\cL_j)$ with $\cL_j:=\overline{z_{j+1}z_{j+2}}\cup\overline{z_{j+2}z_{j+3}}$.
\end{enumerate} 
\end{defi}
Notice that the functions $\ell$ and $r$ are continuous and non-negative. Moreover, $r(Z) \leq \ell(Z)$ for all $Z\in\C^4$, and $\ell(Z)=0$ if and only if $Z = (z,z,z,z)$, for some $z \in \C$. Thus, the function $h\colon\C^4\smallsetminus\{(z,z,z,z)\}\to\R$ defined as $h(Z)=\frac{r(Z)}{\ell(Z)}$ is continuous and $0\leq h(Z)\leq 1$. For all $Z\in\C^4\smallsetminus\{(z,z,z,z)\}$ and $a,b\in\C$ with $a\neq 0$, the properties 
    \[  r(aZ+b)=|a|r(Z) \quad \textrm{and} \quad \ell(aZ+b)=|a|\ell(Z)   \]
hold, so that $h$ is invariant under the action of the complex affine group $\cA_{\C}$ on $\C^4$. We conclude that there is a well defined and continuous function $h\colon P(4)\to\R$ in the space of shapes of quadrilaterals.

\begin{defi}
The height function in the space of shapes of simple and positively oriented quadrilaterals is the restriction $h\colon S(4)\to\R$, defined by
    \[  h[z_1,z_2,z_3,z_4]=\frac{r(z_1,z_2,z_3,z_4)}{\ell(z_1,z_2,z_3,z_4)}.    \]
\end{defi}

\begin{prop}
\label{globalmax}
The height function $h\colon S(4)\to\R$ is positive and has a global maximum in the square $\mathfrak{R}_4=[0,1,1+i,i]$, where $h$ equals $1$.
\end{prop}

\begin{proof}
It is straightforward to check from Definition \ref{distances} that $h(\mathfrak{R}_4)=1$ and $h(Z)>0$, for every $Z\in S(4)$. If $Z$ has one side strictly smaller than the others, say $|z_j-z_{j+1}|$ for some $j\in\{1,2,3,4\}$, then $r(Z)\leq r_j(Z)\leq|z_j-z_{j+1}|<\ell(Z)$, and therefore, $h(Z)<1$. If all of the sides of $Z\in S(4)$ have equal length (we may choose a representative whose sides are equal to 1), then it is a rhombus. If $Z \neq \mathfrak{R}_4$, then there exists $j\in\{1,2,3,4\}$, such that the interior angle at the vertex $z_j$ satisfies $0<\theta_j<\frac{\pi}{2}$ which implies that $h(Z)=r(Z)=r_{j+1}(Z)=\sin(\theta_j)<1$.
\end{proof}

\begin{prop}
\label{zeroboundary}
The height function $h\colon S(4)\to\R$ is continuously extended to the closure $\overline{S(4)}$ such that $h(Z) = 0$ if and only if $Z \in \partial S(4)$.
\end{prop}

\begin{proof}
For every $Z=[z_1,z_2,z_3,z_4] \in \overline{S(4)}$, we have that $h(Z) = 0$ if and only if there is an index $j \in \{1,2,3,4\}$, such that $r_j(Z) = d(z_j, \cL_j) = 0$, with $\cL_j=\overline{z_{j+1}z_{j+2}}\cup\overline{z_{j+2}z_{j+3}}$ the opposite polygonal to the vertex $z_j$. Thus $h(Z) = 0$ if and only if $z_j \in \cL_j$, for some $j \in \{1,2,3,4\}$, and this is true if and only if $Z \in \partial S(4)$. We may observe that the property $\ell(Z)>0$ holds for every $Z \in \overline{S(4)}$.
\end{proof}

\begin{obs}
\label{invarianza}
We can observe that $h(Z)=h(\sigma(Z))=h(\tau(Z))$, for all $Z\in P(4)$, so that the level sets of the height function $h\colon S(4)\to\R$ are invariant under the action of the dihedral group $\langle\sigma,\tau\rangle$.
\end{obs}

\begin{figure}[h]
\begin{center}
\begin{tikzpicture}[scale=1.5]
\draw [black!50, dashed] (-1,0)--(2,0);
\draw (-.1,0) node [below] {\scriptsize{$[0,1,0]$}};
\draw (1,0) node [below] {\scriptsize{$[0,1,1]$}};
\draw (2,0) node [right] {\small{$\partial S(3)\cong\R\cup\{\infty\}$}};
\draw (-1.7,2.3) node [right] {\small{$S(3)\cong\H$}};

\draw (0,0) node {\scriptsize{$\bullet$}};
\draw (1,0) node {\scriptsize{$\bullet$}};
\draw (.5,.866) node {\scriptsize{$\bullet$}};

\draw [thick=0.4] (1,0) arc (0:60:1);
\draw [thick=0.4] (.5,.866) arc (120:180:1);
\draw [thick=0.4] (.5,.866)--(.5,2.5);
\draw (.47,.95) node [right] {$\mathfrak{R}_3$};

\draw [black!70] (.292,.7071)--(.7071,.7071);
\draw [black!70] (.7071,.7071) arc (-40:75:.4);
\draw [black!70] (.292,.7071) arc (-140:-255:.4);

\draw [black!70] (.1335,1/2)--(.866,.5);
\draw [black!70] (.866,.5) arc (-40:70:.85);
\draw [black!70] (.1335,.5) arc (-140:-250:.85);

\draw [black!70] (1-.966,.258)--(.966,.258);
\draw [black!70] (.966,.258) arc (-40:64:1.4);
\draw [black!70] (.044,.258) arc (-140:-244:1.4);
\end{tikzpicture}
\caption{Three level curves of the height function $h\colon S(3)\to\R$, from the inside out, these correspond to $h^{-1}(\sin(\frac{\pi}{4}))$, $h^{-1}(\sin(\frac{\pi}{6}))$, and $h^{-1}(\sin(\frac{\pi}{12}))$.}
\label{figure3}
\end{center}
\end{figure}
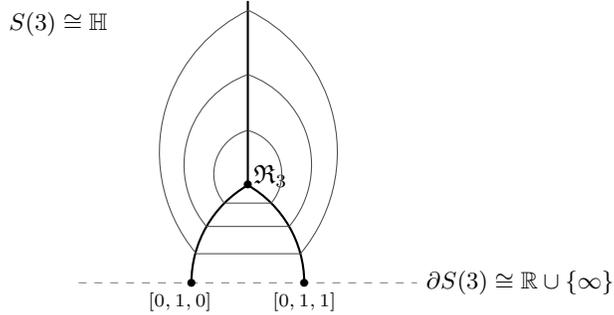

\begin{ejem}
\label{eje5}
The analogous height function for the case of triangles is given by $h\colon S(3)\to\R$, $h(Z)=\frac{r(Z)}{\ell(Z)}$, with $Z=[z_1,z_2,z_3]$, $r(Z)=\min_jd(z_j,\overline{z_{j+1}z_{j+2}})$ and $\ell(Z)=\max_j|z_j-z_{j+1}|$. In this case the function $h$ also has a global maximum at the equilateral triangle $\mathfrak{R}_3=[0,1,(1+i\sqrt{3})/2]$, where it attains the value $\sqrt{3}/2$, and can be continuously extended by the zero function in $\partial S(3)=\R\cup\{\infty\}$. If we consider the transverse curve $\mu : [0,\infty) \rightarrow S(3)$ given by $\mu(t) = [0,1,1/2 + i e^t \sqrt{3}/2]$, then the level curves of $h$ together with the curve $\mu(t)$ provide polar coordinates to the space $S(3)$, giving us the homeomorphism with the open cone $S(3) \cong C(\S^1) \cong \R^2$, see Figure \ref{figure3}.
\end{ejem}

\begin{obs}\label{triangle_space_decomposition}
The three regions bounded by the vertical line and the two arcs of circles of radius $1$ passing through $0, 1 \in \R \cup \{\infty\} \cong \partial S(3)$ in Figure \ref{figure3}, correspond to the triangles with fixed largest side, that is, those regions are the sets 
\[  \mathrm{cl}(\cU_j) = \{ Z \in S(3) \colon \ell(Z) = \ell_j(Z) \}.  \]
The intersection of each one of these sets with a level set is a topological segment $h^{-1}(s) \cap \mathrm{cl}(\cU_j) \cong \D^1$ and thus, $h^{-1}(s)$ is a union of three segments with pairwise intersection in a single point. This trivial fact on the space of triangles generalizes to the space of quadrilaterals in a more complicated combinatorial decomposition.
\end{obs}

 \section{Level sets of the height function}
\label{levelcurves}

In this section we prove that for all $s\in(0,1)$, the level set $h^{-1}(s)$ is homeomorphic to the sphere $\S^3$. As this section is much more technical than the previous ones, before we start, we summarize the main points of the proof in a general strategy:
\begin{enumerate}
    \item Introduce the sets in $S(4)$, where $\ell$ and $r$ are uniquely reached, this is Definition \ref{conos}. In $S(3)$ there are three such sets, completely determined by $\ell$, as seen in Remark \ref{triangle_space_decomposition}.
    
    \item Show that for all $s \in (0,1)$, the level set $h^{-1}(s)$, intersected with the closure of each of the sets in the previous step, is homeomorphic to the closed ball $\D^3$, this is Lemma \ref{nivelconos}. In the case of $S(3)$, it is clear that the corresponding intersections are topological segments, Figure \ref{figure3}.
    
    \item Describe the decomposition of $h^{-1}(s)$ into closed balls, by codifying the combinatorics of the intersections in a graph, this is Lemma \ref{pegados}. The corresponding behaviour in $S(3)$ is the fact that $h^{-1}(s)$, as a union of closed segments, is a topological triangle, and the non-trivial intersection of the segments are described by the dual polygon which is again a triangle.
    
    \item Conclude that $h^{-1}(s) \cong \S^3$ by constructing a  decomposition of $\S^3$ into closed balls which is equivalent to the one given for $h^{-1}(s)$. This is given in Theorem \ref{levelspheres}.
\end{enumerate}

\begin{defi}
\label{min}
Let $Z=(z_1,z_2,z_3,z_4)\in\S(4)$. Then
\begin{enumerate}[(a)]
    \item we say that $\ell(Z)$ is unique if there exists $j\in\{1,2,3,4\}$, such that $\ell(Z)=|z_j-z_{j+1}|$ and for $k\neq j$, $|z_k-z_{k+1}|<|z_j-z_{j+1}|$.
    
    \item we say that $r(Z)$ is unique if there exist $j\in\{1,2,3,4\}$ and $w\in\cL_j$, such that $r(Z)=|z_j-w|$ and $r(Z)<|z_k-w'|$, for every $k\in\{1,2,3,4\}$ and $w'\in\cL_k$ such that $\{z_k,w'\}\neq\{z_j,w\}$.
\end{enumerate} 
\end{defi}

Notice that in a quadrilateral, $r$ cannot be realized as the distance between two opposite vertices, because if two non-adjacent vertices are close, then $r$ is reached between one of them and a side non-adjacent to it. However, $r$ could be reached from a vertex to a two distinct non-adjacent sides.

\begin{prop}
\label{property}
If $Z\in S(4)$ is such that $\ell(Z)=|z_j-z_{j+1}|$ is unique, then $d(z_j,\overline{z_{j+1}z_{j+2}})>r(Z)$ and $d(z_{j+1},\overline{z_{j+3}z_j})>r(Z)$. 
\end{prop} 

\begin{proof}
We can assume that $j=1$ and prove that $d(z_2,\overline{z_4z_1})>r(Z)$, as the other cases are analogous. We work with the representative $Z=[0,1,z_3,z_4]$.

If the interior angle at the vertex $z_1$, denoted by $\theta_1$ here, satisfies $\theta_1\geq\frac{\pi}{2}$, then $d(z_2,\overline{z_4z_1})=|z_2-z_1|=\ell(Z)=1>r(Z)$, by Proposition \ref{globalmax}. 

For $\theta_1<\frac{\pi}{2}$. If $\measuredangle z_1z_4z_2\geq\frac{\pi}{2}$, then $r(Z)\leq d(z_4,\overline{z_1z_2})=Im(z_4)<|z_2-z_4|=d(z_2,\overline{z_4z_1})$. If $\measuredangle z_1z_4z_2<\frac{\pi}{2}$, then there exists a point $w\in\overline{z_4z_1}$, such that $\measuredangle z_1wz_2=\frac{\pi}{2}$. We conclude that $d(z_2,\overline{z_4z_1})=d(z_2,w)=\ell(Z)\sin\theta_1>|z_4-z_1|\sin\theta_1=d(z_4,\overline{z_1 z_2})\geq r(Z)$.
\end{proof}

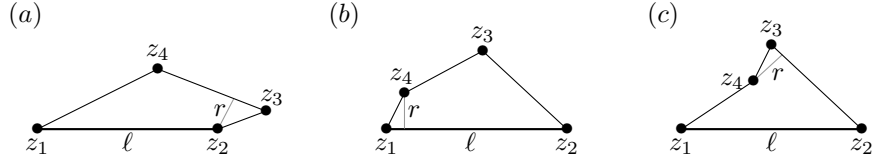
\begin{figure}[h]
\begin{center}
\begin{tikzpicture}[scale=.8]
\draw [thick] (-5,0)--(-2,0);
\draw (-3.5,0.1) node [below] {$\ell$};

\draw [black!40] (-2,0)--(-1.74,0.5);
\draw (-1.7,0.32) node [left] {$r$};

\draw (-5,0)--(-2,0)--(-1.2,0.3)--(-3,1)--cycle;
\draw (-5,0) node {$\bullet$};
\draw (-2,0) node {$\bullet$};
\draw (-1.2,0.3) node {$\bullet$};
\draw (-3,1) node {$\bullet$};
\draw (-5,0) node [below] {$z_1$};
\draw (-2,0) node [below] {$z_2$};
\draw (-1.1,0.25) node [above] {$z_3$};
\draw (-3,1) node [above] {$z_4$};
\draw (-5.2,1.9) node {$(a)$};

\draw [thick] (0.8,0)--(3.8,0);
\draw (2.2,0.1) node [below] {$\ell$};

\draw [black!40] (1.1,0.6)--(1.1,0);
\draw (1,0.3) node [right] {$r$};

\draw (0.8,0)--(3.8,0)--(2.4,1.3)--(1.1,0.6)--cycle;
\draw (0.8,0) node {$\bullet$};
\draw (3.8,0) node {$\bullet$};
\draw (2.4,1.3) node {$\bullet$};
\draw (1.1,0.6) node {$\bullet$};
\draw (0.8,0) node [below] {$z_1$};
\draw (3.8,0) node [below] {$z_2$};
\draw (2.4,1.3) node [above] {$z_3$};
\draw (1.05,0.6) node [above] {$z_4$};
\draw (0.1,1.9) node {$(b)$};

\draw [thick] (5.7,0)--(8.7,0);
\draw (7.2,0.1) node [below] {$\ell$};

\draw [black!40] (6.9,0.8)--(7.38,1.24);
\draw (7.05,0.92) node [right] {$r$};

\draw (5.7,0)--(8.7,0)--(7.2,1.4)--(6.9,0.8)--cycle;
\draw (5.7,0) node {$\bullet$};
\draw (8.7,0) node {$\bullet$};
\draw (7.2,1.4) node {$\bullet$};
\draw (6.9,0.8) node {$\bullet$};
\draw (5.7,0) node [below] {$z_1$};
\draw (8.7,0) node [below] {$z_2$};
\draw (7.2,1.4) node [above] {$z_3$};
\draw (6.9,0.8) node [left] {$z_4$};
\draw (5.4,1.9) node {$(c)$};
\end{tikzpicture}
\caption{Three quadrilaterals with largest side uniquely reached at $\overline{z_1z_2} $. The quadrilaterals in $(a)$, $(b)$ and $(c)$ belong to $\cU_1$, $\cV_1$ and $\cW_1$ respectively.}
\label{figure4}
\end{center}
\end{figure}

\vspace{-0.2cm}
\begin{defi}
\label{conos}
Let $j\in\{1,2,3,4\}$. The set of quadrilaterals $Z\in S(4)$ such that $\ell(Z)=|z_j-z_{j+1}|$ and $r(Z)$ are unique, is divided in the next three subsets (Figure \ref{figure4}):
\begin{enumerate}[(a)]
\item $\cU_j:=\{r(Z)=d(z_{j+1},\overline{z_{j+2}z_{j+3}})\text{ or to }r(Z)=d(z_{j+2},\overline{z_jz_{j+1}})\}.$
\item $\cV_j:=\{r(Z)=d(z_{j+3},\overline{z_jz_{j+1}})\text{ or }r(Z)=d(z_j,\overline{z_{j+2}z_{j+3}})\}$.
\item $\cW_j:=\{r(Z)=d(z_{j+2},\overline{z_{j+3}z_j})\text{ or }r(Z)=d(z_{j+3},\overline{z_{j+1}z_{j+2}})\}.$
\end{enumerate}
\end{defi}

For all $j\in\{1,2,3,4\}$, the sets $\cU_j$, $\cV_j$ and $\cW_j$ are disjoint and every simple and positively oriented quadrilateral with $\ell$ and $r$ unique, belongs to the union $\bigcup_{1\leq j\leq4}(\cU_j\cup\cV_j\cup\cW_j)\subset S(4)$ (Proposition \ref{property}).

\begin{obs}
\label{permutaconos}
The sets $\cU_j,\cV_j$ and $\cW_j$ satisfy $\sigma^k(\cU_j)=\cU_{j+k}$, $\sigma^k(\cV_j)=\cV_{j+k}$, $\sigma^k(\cW_j)=\cW_{j+k}$, $\tau(\cU_1)=\cV_1$, and $\tau(\cW_1)=\cW_1$. In particular, we have that $\sigma^j \circ \tau \circ \sigma^{-j}(\cU_j) = \cV_j$ and $\sigma^j \circ \tau \circ \sigma^{-j}(\cW_j) = \cW_j$.
\end{obs}

\begin{prop}
The set $\bigcup_{j = 1}^{4} \left(\cU_j \cup  \cV_j \cup \cW_j \right)$ is open and dense in $S(4)$.
\end{prop}

\begin{proof}
Observe that every simple quadrilateral can be approximated by a quadrilateral with $\ell$ and $r$ unique, so that $\bigcup_{j = 1}^{4} \left(\cU_j \cup  \cV_j \cup \cW_j \right) \subset S(4)$ is dense.

Since the functions $r_j$ and $\ell_k$ are continuous, then if the inequalities
\[\ell_j(Z)<\ell_k(Z), \qquad  r_j(Z)<r_k(Z)  \]
hold in a quadrilateral $Z$, then the same inequalities remain valid in an open
neighborhood of $Z$. If the minimum satisfies $r(Z)=r_j(Z)=r_{j+1}(Z)$, \textit{i.e.} $r(Z)=|z_j-z_{j+1}|$, then it can be verified that $r$ remains unique in a neighbourhood of $Z$ and the result follows.
\end{proof}

\begin{obs}\label{remark_differentiability_height}
We may observe that the functions $r$ and $\ell$ have unique expressions as $r_k$ and $\ell_l$, for some $k,l \in \{1,2,3,4\}$, when restricted to one of the sets $\{\mathrm{cl}(\cU_j), \mathrm{cl}(\cV_j), \mathrm{cl}(\cW_j) : j \in \{1,2,3,4\} \}$ and this implies in particular that the height function $h$ is differentiable in such restrictions.
\end{obs}

\begin{lem}
\label{nivelconos}
If  $j\in\{1,2,3,4\}$ and $s\in(0,1)$, then the closures $\mathrm{cl}(h^{-1}(s)\cap\cU_j)$, $\mathrm{cl}(h^{-1}(s)\cap\cV_j)$ and $\mathrm{cl}(h^{-1}(s)\cap\cW_j)$ are homeomorphic to the closed disc $\D^3$.
\end{lem}

\begin{proof}
From Remarks \ref{invarianza} and \ref{permutaconos}, it is enough to show the result for $\mathrm{cl}(h^{-1}(s)\cap\cU_1)$ and $\mathrm{cl}(h^{-1}(s)\cap\cW_1)$. We will construct these sets separately by working on the representatives $Z=[0,1,z_3,z_4]$. By Definition \ref{conos}, in both cases the constrains $\ell(Z)=1$, $r(Z)=s$, $s\leq|z_4|\leq 1$ and $s\leq|z_3-1|\leq 1$ must be satisfied.

For each set $\mathrm{cl}(h^{-1}(s)\cap\cU_1)$ and $\mathrm{cl}(h^{-1}(s)\cap\cW_1)$, the proof will be carried out in three steps: 1 - Show that the set $\Omega$ of possible points for $z_4$ is homeomorphic to $\D^2$ by giving a precise description of it. 2 - For a fixed $z_4\in\Omega$, determine the curve $\gamma_{z_4}$ of possible values for $z_3$. 3 - Prove that steps 1 and 2 endow the set with a structure of an $I$-bundle over $\Omega$ (with some collapsed fibers), which will imply the homeomorphism with $\D^3$.

\smallskip
\noindent
\textit{Construction of $\mathrm{cl}(h^{-1}(s)\cap\cU_1)$}. Suppose that $[0,1,z_3,z_4]\in\mathrm{cl}(h^{-1}(s)\cap\cU_1)$. 

\noindent
\textit{Step 1.} By Definition \ref{conos}.$(a)$, $s=d(z_3,\overline{01})$ or $s=d(1,\overline{z_3z_4})$, then either the vertex $z_3$ belongs to the segment $\overline{(is)(1+is)}$ or to the circle $\cC_1=\{|z-1|=s\}$, or the side $\overline{z_3z_4}$ is tangent to $\cC_1$. These conditions and the inequalities $d(0,\overline{z_3z_4})\geq s$ and $d(z_4,\overline{01})\geq s$, imply that $\Omega_{\cU s}=\{z\in\D^2\colon Im(z)\geq s\}$ is the set of possible points for $z_4$, see Figure \ref{Omega_U}.

\noindent
\textit{Step 2.} Fix a point $z_4\in\Omega_{\cU s}$. Let $p_{z_4}\in\cC_1$ be the point satisfying the condition $Re(p_{z_4})\geq 1$ and such that $\overline{p_{z_4} z_4}$ is tangent to $\cC_1$. Denote by $\cA_{z_4}\subset\cC_1$, the arc in the negative sense from $1+is$ to $p_{z_4}$, and $\cR_{z_4}=\{t(p_{z_4}-z_4)+p_{z_4}\colon t\geq0\}$ the ray starting at $p_{z_4}$. The set $\gamma_{z_4}$ of possible values for $z_3$ is$$\{z\in\overline{(is)(1+is)}\cup\cA_{z_4}\cup\cR_{z_4}\colon d(z_4,\overline{1z})\geq s,\ d(z,\overline{z_40})\geq s,~\text{and}~~|z-z_4|\leq 1 \}.$$Observe that for $z_4\in\Omega_{\cU s}$ such that $|z_4-1|=\sqrt{2}s$ (equivalently the circles $\cC_1$ and $\{|z-z_4|=s\}$ are orthogonal),  $d(z_4,\overline{1 z})=s$ if and only if $d(z,\overline{z_40})=s$, and therefore, the constraints $d(z_4,\overline{1 z})\geq s$ and $d(z,\overline{z_40})\geq s$ are equivalent. Then, for the case $|z_4 -1|\geq\sqrt{2}s$ it is satisfied that$$\gamma_{z_4}=\{z\in\overline{(is)(1+is)}\cup\cA_{z_4}\cup\cR_{z_4}\colon d(z,\overline{z_40})\geq s,~\text{and}~~|z-z_4|\leq 1\},$$ (illustrated by the construction on the points $z_4$ and $z_4''$ in Figure \ref{Omega_U}) and for the case $|z_4-1|\leq\sqrt{2}s$, it is satisfied that$$\gamma_{z_4}=\{z\in\overline{(is)(1+is)}\cup\cA_{z_4}\cup\cR_{z_4}\colon d(z_4,\overline{1z})\geq s,~\text{and}~~|z-z_4|\leq 1 \},$$(illustrated by the construction on the point $z_4'$ in Figure \ref{Omega_U}). As the curve $\gamma_{z_4}$ is constructed by an arc and two segments tangent to it, it is a $C^1$-curve whose construction only depends on the point $z_4$ and the parameter $s\in(0,1)$.

\begin{figure}[h]
\begin{center}
\begin{tikzpicture}[scale=1]
\fill [gray!30]
(2.83,1)--(-2.83,1)--(2.83,1) arc (19.5:160.5:3);
\draw [black!50, dashed] (-3.5,0)--(6.5,0);
\draw [black!50, dashed] (0,-0.5)--(0,3.3);
\draw [black!50] (4,0) arc (0:180:1);
\draw [black!50] (4.4142,0) arc (0:180:1.4142);
\draw [black!60] (2.83,1)--(-2.83,1);
\draw (3,0) arc (0:180:3);
\draw (6,0) arc (0:180:3);

\draw [black!60] (0,0)--(2.55,1.1);
\draw (2.4,1.3) node {\footnotesize{$z'_4$}};
\draw [black!15] (3.55,1.1) arc (0:360:1);
\draw [black!15,shift={(3,0)},rotate=54.5] (0,0)--(2,0);
\draw [black!15] (3.37,0.55) node {\tiny{$\bullet$}};
\draw [dotted,shift={(2.55,1.1)},rotate=-9] (0,0)--(3.3,0);
\draw (2.55,1.1) node {\tiny{$\bullet$}};
\draw (3.66,0.92) node {\tiny{$\bullet$}};
\draw [thick=0.8,shift={(3.68,0.92)},rotate=-9] (0,0)--(1.84,0);
\draw (5.5,0.63) node {\tiny{$\bullet$}};
\draw (4.8,0.9) node {\footnotesize{$\gamma_{z'_4}$}};
\draw (3.8,2.1) node {\footnotesize{$\{|z\!-\!z'_4|\!=\!s\}$}};

\draw [black!30,shift={(2.225,1)},rotate=142] (0,0)--(1,0);
\draw (1.8,1.47) node {\footnotesize{$s$}};
\draw [black!60] (0,0)--(1.8,2);
\draw [dotted,shift={(1.8,2)},rotate=-33.35] (0,0)--(3.5,0);
\draw (1.8,2) node {\tiny{$\bullet$}};
\draw (1.65,2.15) node {\footnotesize{$z_4$}};
\draw [thick=0.8] (2.225,1)--(3,1);
\draw [thick=0.8] (3,1) arc (-270:-302:1);
\draw [thick=0.8,shift={(3.55,0.84)},rotate=-33.35] (0,0)--(0.9,0);
\draw (2.225,1) node {\tiny{$\bullet$}};
\draw (3.55,0.84) node {\tiny{$\bullet$}};
\draw (4.28,0.35) node {\tiny{$\bullet$}};
\draw (3.45,0.76) node {\footnotesize{$p_{z_4}$}};
\draw (2.95,0.8) node {\footnotesize{$\gamma_{z_4}$}};

\draw [black!60] (0,0)--(-2,1.3);
\draw (-2,1.3) node {\tiny{$\bullet$}};
\draw (-2.2,1.4) node {\footnotesize{$z''_4$}};
\draw [thick=0.8] (0.28,1)--(0.98,1);
\draw [black!15,shift={(0.28,1)},rotate=-130] (0,0)--(1,0);
\draw (-0.23,0.55) node {\footnotesize{$s$}};
\draw (0.58,0.8) node {\footnotesize{$\gamma_{z''_4}$}};

\draw (0.28,1) node {\tiny{$\bullet$}};
\draw (0.98,1) node {\tiny{$\bullet$}};
\draw (-0.15,1.15) node {\footnotesize{$is$}};
\draw (0,1) node {\tiny{$\bullet$}};
\draw [black!60] (2.83,1) node {\tiny{$\bullet$}};
\draw (2.85,1.16) node {\footnotesize{$q_{_1}$}};
\draw [black!60] (-2.83,1) node {\tiny{$\bullet$}};
\draw (-3,1.18) node {\footnotesize{$q_{_2}$}};
\draw [black!60] (5.83,1) node {\tiny{$\bullet$}};
\draw (6.1,1.2) node {\footnotesize{$q_{_1}\!\!+\!1$}};
\draw (-3,0.03) node [below] {$-1$};
\draw (0,3) node {\tiny{$\bullet$}};
\draw (0.1,3.2) node [left] {$i$};
\draw (-3,0) node {\tiny{$\bullet$}};
\draw (0,0.03) node [below] {$0$};
\draw (0,0) node {\tiny{$\bullet$}};
\draw (3,0.03) node [below] {$1$};
\draw (3,0) node {\tiny{$\bullet$}};
\draw (6,0.03) node [below] {$2$};
\draw (6,0) node {\tiny{$\bullet$}};
\draw (-1.3,2.2) node {$\Omega_{\cU s}$};

\draw (4.4142,0) node {\tiny{$\bullet$}};
\draw (4.4,-0.18) node {\scriptsize{$1\!\!+\!\!s\sqrt{2}$}};
\draw (2,0) node {\tiny{$\bullet$}};
\draw (2,-.18) node {\footnotesize{$1\!\!-\!\!s$}};
\end{tikzpicture}
\caption{The shaded region is the set $\Omega_{\cU s}$. For the points $z_4,z'_4,z''_4\in\Omega_{\cU s}$, the curves $\gamma_{z_4}, \gamma_{z'_4}$ and $\gamma_{z''_4}$ of possible values for the vertex $z_3$ are represented with the darkest lines.}
\label{Omega_U}
\end{center}
\end{figure}
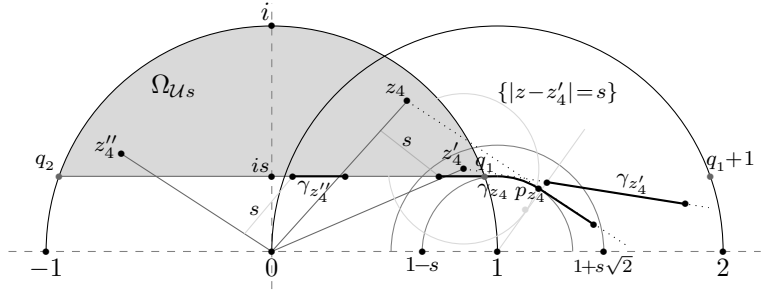

\noindent
\textit{Step 3.} As noted above, for every $z_4\in\Omega_{\cU s}$, $\gamma_{z_4}$ is a $C^1$-curve and so, it can be parametrized by a rescaling of the arc-length as a continuous function $\gamma_{z_4}\colon[0,1] \rightarrow \mathrm{cl}(h^{-1}(s)\cap\cU_1)$. Observe that for the points $q_1=e^{i\arcsin(s)}$ and $q_2=-e^{-i\arcsin(s)}$ (Figure \ref{Omega_U}), we have that $\gamma_{q_1}=\{q_1+1\}$ and $\gamma_{q_2}=\{q_2+1\}$, and these are the only points with collapsed curves. The continuous function \[\Omega_{\cU s} \times [0,1] \rightarrow\mathrm{cl}(h^{-1}(s)\cap\cU_1), \qquad (z_4,t) \rightarrow [0,1,\gamma_{z_4}(t),z_4],\] is a homeomorphism when restricted to $\mathrm{cl}(h^{-1}(s)\cap\cU_1)\smallsetminus\{[0,1,q_1+1,q_1],[0,1,q_2+1,q_2]\}$. This map endows the set $\mathrm{cl}(h^{-1}(s)\cap\cU_1)$ with the structure of an $I$-bundle with two collapsed fibers at the points $q_1$ and $q_2$ in the boundary of $\Omega_{\cU s}$. Then, we conclude that $\mathrm{cl}(h^{-1}(s)\cap\cU_1)$ is homeomorphic to $\D^3$. 

\smallskip
\noindent
\textit{Construction of $\mathrm{cl}(h^{-1}(s)\cap\cW_1)$}. Suppose that $[0,1,z_3,z_4]\in\mathrm{cl}(h^{-1}(s)\cap\cW_1)$.

\noindent
\textit{Step 1.} Since $s=d(z_3,\overline{z_40})$ or $s=d(z_4,\overline{1z_3})$ (Definition \ref{conos}.$(c)$), then the vertices $z_3$ and $z_4$ have analogous roles. As $d(z_4,\overline{01})\geq s$, then for $0\leq Re(z_4)\leq 1$, necessarily $Im(z_4)\geq s$. In order to determine the possible values for $z_4$ with $Re(z_4)\leq 0$, we analyze the constraints of $z_4$ when belonging to the radii $r_{\alpha}=\{te^{i\alpha}\colon 0\leq t\leq1\}$ with argument $\alpha\geq\frac{\pi}{2}$. We denote by $\cC_0=\{|z|=s\}$ the circle of radius $s$ centered at $0$, and by $L_{\alpha}=\{(t-is)e^{i\alpha}\colon t\in\R\}$ the line parallel to $r_{\alpha}$, at distance $s$ which is tangent to $\cC_0$ at $-ise^{i\alpha}$. The tangent to $\cC_0$ at $se^{i\alpha}$ intersects $L_\alpha$ in the point $p_{\alpha}=(s-is)e^{i\alpha}$ (Figure \ref{construccion_chipote}.$(a)$). There are three cases:

\begin{enumerate}[i.]
\item If $|p_{\alpha}-1|\leq 1$, then $[0,1,p_\alpha,se^{i \alpha}]\in  \mathrm{cl}(h^{-1}(s)\cap\cW_1)$. The set of possible values for $z_4$ in $r_\alpha$ is $\{te^{i\alpha}\colon s\leq t\leq 1\}$. This case occurs for arguments $\frac{\pi}{2}\leq\alpha\leq\alpha_0$, with$$\alpha_0=\arccos\left(-\sqrt{\frac{1-s\sqrt{2-s^2}}{2}}\right),$$the angle such that $|p_{\alpha_0}-1|=1$ holds (see Figure \ref{construccion_chipote}.$(a)$).
   
\item If $|p_{\alpha}-1|>1$, then $p_\alpha$ is not a valid value for $z_3$ anymore. Let $w_\alpha \in L_{\alpha}\cap\{|z-1|=1\}$ be the intersection point with $\Im(w_\alpha) > 0$, and $q_\alpha \in r_{\alpha}$ so that the segment $\overline{q_\alpha w_\alpha}$ is tangent to $\cC_0$ (Figure \ref{construccion_chipote}.$(b)$). If $$\lambda_\alpha=\cos(\alpha)+\sqrt{\cos(\alpha)^2+s(2\sin(\alpha)-s)},$$
then $w_\alpha=(\lambda_{\alpha}-is)e^{i \alpha}$ and $q_\alpha=t_{\alpha}e^{i\alpha}\in r_{\alpha}$, with $t_\alpha=\frac{s^2+\lambda_\alpha^2}{2 \lambda_\alpha}$. The quadrilateral $[0,1,q_\alpha,w_\alpha]$ belongs to $\mathrm{cl}(h^{-1}(s)\cap\cW_1)$, and the set of possible values for $z_4$ in $r_{\alpha}$ is $\{te^{i\alpha}\colon t_{\alpha}\leq t\leq1\}$. This case occurs for arguments $\alpha_0<\alpha\leq\pi-\arcsin(s)$. Notice that for $\alpha=\alpha_0$, $w_\alpha=p_\alpha$ and $q_\alpha=se^{i\alpha}$, and in the limit angle $\alpha=\pi-\arcsin(s)$, $q_\alpha=-e^{-i\arcsin(s)}=q_2$ (see Figure \ref{construccion_chipote}.$(b)$).
   
\item If $\alpha>\pi-\arcsin(s)$, then for all $s\leq t\leq1$, the tangent line to $\cC_0$ through $te^{i\alpha}\in r_\alpha$ intersects $L_\alpha$ outside the circle $\{|z-1|=1\}$. We conclude that there are no allowed values of $z_4$ in the radii $r_{\alpha}$.
\end{enumerate}

\begin{figure}[h]
\begin{center}
\begin{tikzpicture}[scale=0.75]
\draw (-11.3,4.3) node [below] {$(a)$};
\draw [black!20,shift={(-9,0.58)},rotate=31.5] (0,0)--(2,0);
\draw [black!20,shift={(-8,0)},rotate=119] (0,0)--(3,0);
\draw [black!20,shift={(-6.95,0.2)},rotate=119] (-0.5,0)--(4,0);

\draw (-5,0) arc (0:180:3);
\draw (-7,0) arc (0:180:1);
\draw (-8,0) arc (-180:-280:3);
\draw [black!50,dashed] (-11.5,0)--(-4.5,0);
\draw [black!50,dashed] (-8,-0.5)--(-8,3.5);
\draw [dotted] (-10.82,1.02)--(-5.18,1.02);

\draw [black!70,shift={(-6.965,0)},rotate=103] (-0.5,0)--(4,0);
\draw [black!70,shift={(-6.9,1.3)},rotate=90] (0,0)--(-0.356,1.52);
\draw [black!70,shift={(-8,0)},rotate=103] (0,0)--(1,0);
\draw [thick=0.8,shift={(-8.2249,0.9743)},rotate=103] (0,0)--(2,0);

\draw [black] (-9,0) node {\tiny{$\bullet$}};
\draw (-9.1,-0.2) node {{\footnotesize $-s$}};

\draw (-7.35,1.45) node [right] {{\footnotesize $p_{\alpha}$}};
\draw [black] (-7.24,1.22) node {\tiny{$\bullet$}};

\draw (-7.95,1.25) node [above] {{\footnotesize{\color{gray} $p_{\alpha_0}$}}};
\draw [black] (-7.63,1.42) node {\tiny{$\bullet$}};

\draw (-7.27,.65) node [right] {{\footnotesize{\color{gray} $-ise^{i\alpha_0}$}}};
\draw [black] (-7.14,0.52) node {\tiny{$\bullet$}};

\draw (-7.15,.25) node [right] {{\footnotesize $-ise^{i\alpha}$}};
\draw [black] (-7.0251,.2225) node {\tiny{$\bullet$}};

\draw (-8.82,3.1) node {{\footnotesize $e^{i\alpha}$}};
\draw [black] (-8.6675,2.9247) node {\tiny{$\bullet$}};

\draw (-9.6,2.8) node {{\footnotesize{\color{gray} $e^{i\alpha_0}$}}};
\draw [black] (-9.45,2.625) node {\tiny{$\bullet$}};

\draw [black] (-8.5,.85) node {\tiny{$\bullet$}};

\draw [black] (-8.2225,.9749) node {\tiny{$\bullet$}};

\draw (-8,0.05) node [below] {$0$};
\draw [black] (-8,0) node {\tiny{$\bullet$}};
\draw (-5,0.05) node [below] {$1$};
\draw [black] (-5,0) node {\tiny{$\bullet$}};
\draw (-11.1,0.05) node [below] {$-1$};
\draw [black] (-11,0) node {\tiny{$\bullet$}};
\draw (-8.1,3.2) node {$i$};
\draw [black] (-8,3) node {\tiny{$\bullet$}};

\draw (-10.7,1.2) node [left] {$q_{_2}$};
\draw [black!60] (-10.82,1.02) node {\tiny{$\bullet$}};
\draw (-5.3,1.2) node [right] {$q_{_1}$};
\draw [black!60] (-5.18,1.02) node {\tiny{$\bullet$}};
\draw (-8.8,2.4) node {{\footnotesize $r_{\alpha}$}};
\draw (-9.39,2) node {{\footnotesize {\color{gray} $r_{\alpha_0}$}}};
\draw (-7.5,3.75) node {{\footnotesize $L_{\alpha}$}};
\draw (-9.25,3.7) node {{\footnotesize {\color{gray} $L_{\alpha_0}$}}};

\draw (-3.3,4.3) node [below] {$(b)$};

\draw [black!20,shift={(-0.5,0.58)},rotate=31.5] (0,0)--(2,0);
\draw [black!20,shift={(0.5,0)},rotate=119] (0,0)--(3,0);
\draw [black!20,shift={(1.55,0.2)},rotate=119] (-0.4,0)--(4,0);

\draw (3.5,0) arc (0:180:3);
\draw (1.5,0) arc (0:180:1);
\draw (0.5,0) arc (-180:-280:3);
\draw [black!50,dashed] (-3,0)--(4,0);
\draw [black!50,dashed] (0.5,-0.5)--(0.5,3.5);
\draw [dotted] (-2.345,1.02)--(3.3,1.02);

\draw [black!70,shift={(0.5,0)},rotate=140] (0,0)--(1.2,0);
\draw [thick=0.8,shift={(-0.4192,0.7713)},rotate=140] (0,0)--(1.8,0);
\draw [black!70,shift={(1.81,0.22)},rotate=50] (0,0)--(0,4.2);
\draw [black!70,shift={(-0.7,0.675)},rotate=18.05] (0,0)--(2,0);

\draw (-0.6,-0.2) node {{\footnotesize $-s$}};
\draw [black] (-0.5,0) node {\tiny{$\bullet$}};

\draw [black] (0.5,3) node {\tiny{$\bullet$}};
\draw (-0.35,0.7) node [above] {{\footnotesize $q_{\alpha}$}};
\draw [black] (-0.41,0.768) node {\tiny{$\bullet$}};

\draw (0.72,1.05) node [right] {{\footnotesize $w_{\alpha}$}};
\draw [black] (0.72,1.125) node {\tiny{$\bullet$}};
\draw [black] (0.18,0.935) node {\tiny{$\bullet$}};

\draw (1.05,.82) node [right] {{\footnotesize $-ise^{i\alpha}$}};
\draw [black] (1.12,0.78) node {\tiny{$\bullet$}};

\draw (-2.08,2) node {{\footnotesize $e^{i\alpha}$}};
\draw [black] (-1.795,1.915) node {\tiny{$\bullet$}};

\draw (1.32,1.2) node [above] {{\footnotesize{\color{gray} $p_{\alpha_0}$}}};
\draw [black] (0.87,1.42) node {\tiny{$\bullet$}};

\draw (1.23,.5) node [right] {{\footnotesize{\color{gray}  $-ise^{i\alpha_0}$}}};
\draw [black] (1.36,0.52) node {\tiny{$\bullet$}};

\draw (3.5,0.05) node [below] {$1$};
\draw [black] (3.5,0) node {\tiny{$\bullet$}};
\draw (0.4,3.2) node {$i$};
\draw (0.5,0.05) node [below] {$0$};
\draw [black] (0.5,0) node {\tiny{$\bullet$}};
\draw (-2.6,0.05) node [below] {$-1$};
\draw [black] (-2.5,0) node {\tiny{$\bullet$}};

\draw (-2.2,1.2) node [left] {$q_{_2}$};
\draw [black!60] (-2.32,1.02) node {\tiny{$\bullet$}};
\draw (3.2,1.2) node [right] {$q_{_1}$};
\draw [black!60] (3.32,1.02) node {\tiny{$\bullet$}};
\draw [black] (0.02,.85) node {\tiny{$\bullet$}};
\draw (-0.75,2.95) node {{\footnotesize{\color{gray}  $e^{i\alpha_0}$}}};
\draw [black] (-0.95,2.625) node {\tiny{$\bullet$}};
\draw (-1.53,1.48) node {{\footnotesize $r_{\alpha}$}};
\draw (-1.7,3.1) node {{\footnotesize $L_{\alpha}$}};
\draw (-0.89,2) node {{\footnotesize{\color{gray}  $r_{\alpha_0}$}}};
\draw (-0.75,3.7) node {{\footnotesize{\color{gray}  $L_{\alpha_0}$}}};

\end{tikzpicture}
\caption{The darkest lines represent the set of possible values for the vertex $z_4$ in the radius $r_{\alpha}$ with $\frac{\pi}{2}\leq\alpha\leq\pi-\arcsin(s)$. The cases $(a)$ and $(b)$ correspond to $\frac{\pi}{2}\leq\alpha\leq\alpha_0$ and $\alpha_0\leq\alpha\leq\pi-\arcsin(s)$ respectively. The construction for the angle $\alpha_0$ is represented in light gray in both cases.}
\label{construccion_chipote}
\end{center}
\end{figure}
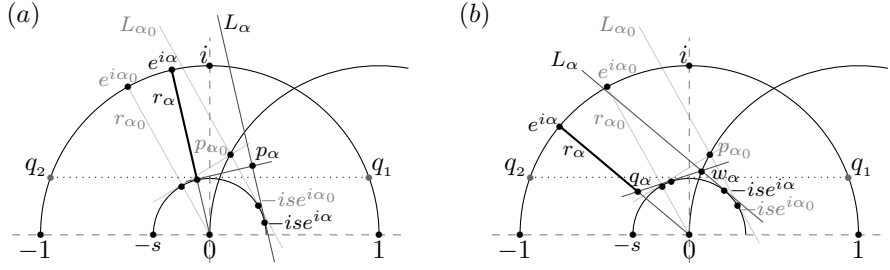

We may define $t_\alpha=\frac{s}{\sin{\alpha}}$ for $\arcsin{s}\leq\alpha\leq\frac{\pi}{2}$, and $t_\alpha=s$ for $\frac{\pi}{2}\leq\alpha\leq\alpha_0$, so that the set of possible values for $z_4$ is (Figure \ref{Omega_W}):$$\Omega_{\cW s}=\{te^{i\alpha}\in r_\alpha\colon t_\alpha\leq t\leq 1,\,\arcsin{s}\leq \alpha\leq\pi-\arcsin(s)\}.$$

\noindent
\textit{Step 2.} Fix a point $z_4\in\Omega_{\cW s}$. Let $\cC_{z_4}=\{|z-z_4|=s\}$ and $\{q_{z_4} \}= L_{\arg(z_4)}\cap\cC_{z_4}$. If $p_{z_4}\in\cC_{z_4}$ is the point such that the line containing $\overline{1p_{z_4}}$ is tangent to $\cC_{z_4}$ and does not intersect the segment $\overline{z_40}$, then $\cA_{z_4}\subset\cC_{z_4}$ denotes the arc in the positive sense from $q_{z_4}$ to $p_{z_4}$ and $\cR_{z_4}=\{t(p_{z_4}-1)+p_{z_4}\colon t\geq0\}$ denotes the ray tangent to $\cC_{z_4}$ starting at $p_{z_4}$. The curve $\gamma_{z_4}$ of possible values for $z_3$ is contained in the union $\overline{(-ise^{i\arg(z_4)})(q_{z_4})}\cup\cA_{z_4}\cup\cR_{z_4}$. 
The boundaries of of $\gamma_{z_4}$ are determined different in two cases: 

\begin{enumerate}[i.]
\item Suppose that $Im(z_4)\geq s$. If $|z_4 -1|\geq\sqrt{2}s$, we define
    \[  \gamma_{z_4}=\{z\in\overline{(-ise^{i\arg(z_4)})(q_{z_4})}\cup\cA_{z_4}\cup\cR_{z_4}\colon d(z,\overline{01})\geq s,~\text{and}~~|z-1|\leq 1\},   \]
and if $|z_4-1|\leq\sqrt{2}s$, we define
    \[\gamma_{z_4}=\{z\in\overline{(-ise^{i\arg(z_4)})(q_{z_4})}\cup\cA_{z_4}\cup\cR_{z_4}\colon d(1,\overline{z z_4})\geq s,~\text{and}~~|z-1|\leq 1 \}.\]
We observe that the conditions $d(1,\overline{z z_4})\geq s$ and $d(z,\overline{01})\geq s$ coincide when $|z_4 - 1| = \sqrt{2}s$.

\item Suppose that $Im(z_4)\leq s$. If $\omega_{z_4}\in L_{\arg(z_4)}$ is the point such that the line $\{ t \omega_{z_4} + (1-t) z_4 : t \in \R\}$ has positive slope and is tangent $\cC_0$, then
    \[  \gamma_{z_4}=\{z\in\overline{(-ise^{i\arg(z_4)})(q_{z_4})}\cup\cA_{z_4}\cup\cR_{z_4}\colon Im(z)\geq Im(\omega_{z_4}),~\text{and}~~|z-1|\leq 1\}.  \]
In the limit case, when $z_4 \in \Gamma$, then with our previous notation $z_4 = q_{\arg(z_4)}$ and $\omega_{z_4} = \omega_{\arg(z_4)}$, see Figure \ref{construccion_chipote}.
\end{enumerate}    
As before, the curve $\gamma_{z_4}$ is constructed by an arc and two segments tangent to it, so it is a $C^1$-curve whose construction only depends on the point $z_4$ and the parameter $s$.

\begin{figure}[h]
\begin{center}
	\includegraphics[scale=0.24]{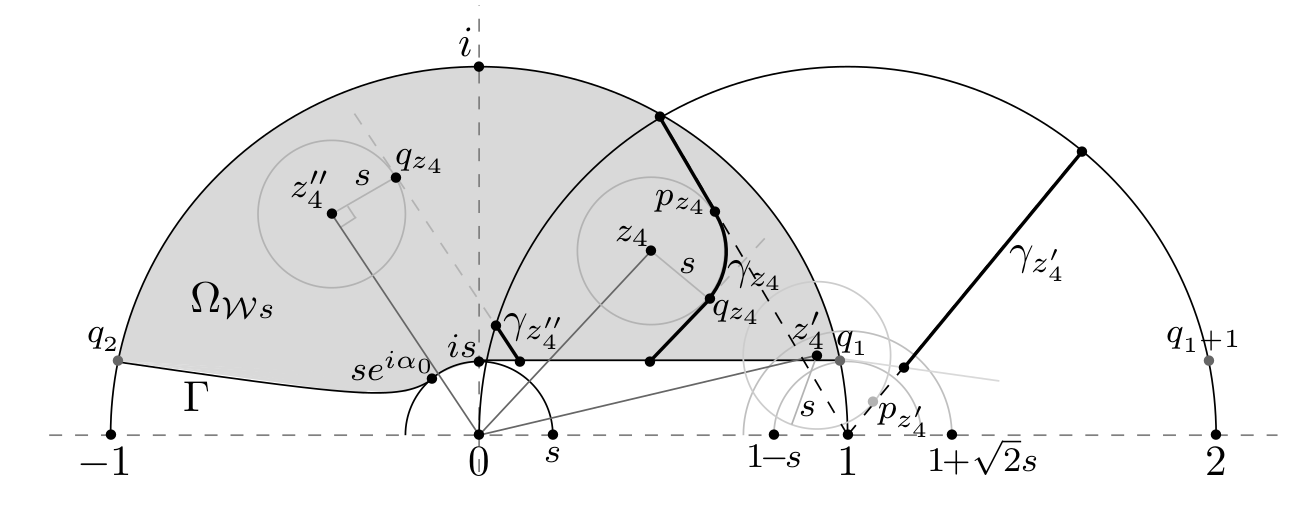}
\caption{The shaded region is the set $\Omega_{\cW s}$. For the points $z_4,z'_4,z''_4\in\Omega_{\cW s}$, the curves $\gamma_{z_4}, \gamma_{z'_4}$ and $\gamma_{z''_4}$ of possible values for the vertex $z_3$ are represented with the darkest lines. The curve $\Gamma=\{q_\alpha\in r_\alpha\colon\alpha_0\leq\alpha\le\pi-\arcsin(s)\}$ is the boundary section between the points $q_2$ and $se^{i\alpha_0}$.}
\label{Omega_W}
\end{center}
\end{figure}

\noindent
\textit{Step 3.} For every $z \in \Omega_{\cW s}$, the curve $\gamma_{z}$ can be parametrized by a rescaling of the arc-length as a continuous function $\gamma_{z}\colon [0,1]\rightarrow \mathrm{cl}(h^{-1}(s)\cap\cW_1)$. We have a projection
    \[\Omega_{\cW s}\times[0,1]\rightarrow \mathrm{cl}(h^{-1}(s)\cap\cW_1),\qquad
    (z,t)\rightarrow [0,1,\gamma_z(t),z],\]
where the collapsed fibers correspond to the points in $\Gamma\cup\{q_1\}$ (Figure \ref{Omega_W}). As the set $\Gamma \cup \{q_1\}$ is contained in $\partial\Omega_{\cW s}$, a homeomorphism $\Omega_{\cW s}\times [0,1] \cong \D^3$ induces a homeomorphism from $\mathrm{cl}(h^{-1}(s)\cap\cW_1)$ to $\D^3$.
\end{proof}

\begin{obs}
\label{recapitulacion_fibrado1}
Three consequences of the proof of Lemma \ref{nivelconos}.
\begin{enumerate}[$I$.]
    \item For $1\leq j\leq 4$, the sets $\mathrm{cl}(h^{-1}(s)\cap\cU_j), \mathrm{cl}(h^{-1}(s)\cap\cV_j), \mathrm{cl}(h^{-1}(s)\cap\cW_j)$ are $I$-bundles with basis homeomorphic to the disc $\D^2$ and with collapsed fibers at some closed regions on their boundaries.

    \item The boundaries of the balls $\mathrm{cl}(h^{-1}(s)\cap\cU_j), \mathrm{cl}(h^{-1}(s)\cap\cV_j), \mathrm{cl}(h^{-1}(s)\cap\cW_j)$ correspond to the loss of uniqueness of the values $\ell$ or $r$ (see Definitions \ref{min} and \ref{conos}). Then for $A\in\{\cU_j, \cW_j, \cV_j : 1 \leq j \leq 4\}$, it is satisfies that\[\partial\mathrm{cl}(h^{-1}(s) \cap A)=\bigcup_{B \neq A} \mathrm{cl}(h^{-1}(s)\cap B)\cap \mathrm{cl}(h^{-1}(s) \cap A),\]where the union is over the set $\{\cU_j, \cW_j, \cV_j : 1 \leq j \leq 4\}$. The decompositions of the boundaries of the balls given by these intersections are permuted by the action of the dihedral group $\langle\sigma,\tau\rangle$.
    
    \item The points $q_1=e^{i\arcsin(s)}$ and $q_2=-e^{-i\arcsin(s)}$ (see Figures \ref{Omega_U} and \ref{Omega_W}) are special because the quadrilaterals $Q_1=[0,1,q_1+1,q_1]$ and $Q_2=[0,1,q_2+1,q_2]$ satisfy\[ \{Q_1, Q_2\} = h^{-1}(s)\cap \bigcap_{j= 1}^4 \mathrm{cl}(\cU_j) \cap \mathrm{cl}(\cV_j) \cap \mathrm{cl}(\cW_j).\] Moreover, the action of the group $\langle\sigma,\tau\rangle$ in this points reduces to $\sigma(Q_1)=\tau(Q_1)=Q_2$,  $\sigma(Q_2)=\tau(Q_2)=Q_1$ and the fibers of the $I$-bundle structures on each sets collapse on these quadrilaterals.
\end{enumerate}
\end{obs}

In order to give a homeomorphism between $h^{-1}(s)$ and $\S^3$, we describe the intersection pattern of the twelve balls of Lemma \ref{nivelconos} through their boundaries. These intersections are codified in a combinatorial structure that we encounter naturally in $\S^3$, in the remainder of the section we describe this combinatorial pattern.

\begin{defi}\label{bigonization-def}
A bigon is a topological closed disc with two marked points on the boundary. The marked points are the vertices of the bigon and the two arcs in the boundary bounded by the vertices are the sides of the bigon.  A decomposition of a surface into bigons is a covering of topological bigons such that the intersection of two such bigons is on the vertices or in one of its sides.
\end{defi}

\begin{figure}[h]
\begin{multicols}{2}
\begin{center}
\begin{tikzpicture}[scale=0.6]
\draw [thick=0.4] (-3,-3)--(3,-3)--(3,3)--(-3,3)--cycle;
\draw (-3,-3) node {$\bullet$};
\draw (3,-3) node {$\bullet$};
\draw (3,3) node {$\bullet$};
\draw (-3,3) node {$\bullet$};

\draw (-3,-3) node [below] {$\cV_3$};
\draw (3,-3) node [below] {$\cW_4$};
\draw (3,3) node [above] {$\cU_4$};
\draw (-3,3) node [above] {$\cW_3$};

\draw [thick=0.4] (-1,-1)--(1,-1)--(1,1)--(-1,1)--cycle;
\draw (-1,-1) node {$\bullet$};
\draw (1,-1) node {$\bullet$};
\draw (1,1) node {$\bullet$};
\draw (-1,1) node {$\bullet$};
\draw (-0.95,-1.1) node [left] {$\cV_1$};
\draw (1,-1.1) node [right] {$\cW_1$};
\draw (1,1.1) node [right] {$\cU_2$};
\draw (-0.95,1.1) node [left] {$\cW_2$};

\draw [thick=0.4] (1,1)--(0,2)--(-1,1);
\draw (0,2) node {$\bullet$};
\draw (0,2) node [above] {$\cV_2$};
\draw [thick=0.4] (-1,1)--(-2,0)--(-1,-1);
\draw (-2,0) node {$\bullet$};
\draw (-2,0) node [left] {$\cU_3$};
\draw [thick=0.4] (-1,-1)--(0,-2)--(1,-1);
\draw (0,-2) node {$\bullet$};
\draw (0,-2) node [below] {$\cU_1$};
\draw [thick=0.4] (1,-1)--(2,0)--(1,1);
\draw (2,0) node {$\bullet$};
\draw (2,0) node [right] {$\cV_4$};

\draw [thick=0.4] (3,3)--(0,2)--(-3,3);
\draw [thick=0.4] (-3,3)--(-2,0)--(-3,-3);
\draw [thick=0.4] (-3,-3)--(0,-2)--(3,-3);
\draw [thick=0.4] (3,-3)--(2,0)--(3,3);
\end{tikzpicture}

\vspace{2cm}
\caption{Graph $\cG$ representing the bigon decomposition of the boundary spheres of the balls $\mathrm{cl}(h^{-1}(s)\cap\cU_j), \mathrm{cl}(h^{-1}(s)\cap\cV_j),$ and $\mathrm{cl}(h^{-1}(s)\cap\cW_j)$.}
\label{Diagrama_pegados}
\end{center}
\end{multicols}
\end{figure}
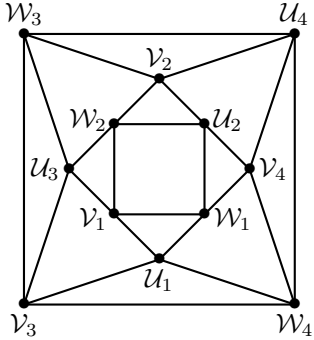

\begin{lem}
\label{pegados}
The spheres $\partial\mathrm{cl}(h^{-1}(s)\cap\cU_j)$, $\partial \mathrm{cl}(h^{-1}(s)\cap\cV_j)$ and $\partial \mathrm{cl}(h^{-1}(s)\cap\cW_j)$ admit a decomposition into four bigons determined by the intersection with their neighboring spheres. More precisely, the graph $\mathcal{G}$ in Figure \ref{Diagrama_pegados} represents such a decomposition as follows:
\begin{enumerate}
    \item each vertex of $\mathcal{G}$ represents a ball in \[\big\{\mathrm{cl}(h^{-1}(s)\cap\cU_j), \mathrm{cl}(h^{-1}(s)\cap\cV_j), \mathrm{cl}(h^{-1}(s)\cap\cW_j) : 1 \leq j \leq 4\big\},\]
    
    \item each edge of $\cG$ represents a bigon which is the intersection of the boundaries of the balls corresponding to its vertices.    
    
    \item each edge of $\mathcal{G}$ is contained in two cycles of the graph of length three and four respectively, those cycles represent the intersection of the balls and determine the sides of the bigon.
    \end{enumerate}
The quadrilaterals $Q_1$ and $Q_2$ are the vertices of all the bigons and the bigon decomposition is preserved by the action of the group $\langle \sigma , \tau \rangle $, see Remark \ref{recapitulacion_fibrado1}.$III$.
\end{lem}

\begin{proof}
Let $A \in \{\cU_j, \cW_j, \cV_j\colon 1 \leq j \leq 4\}$. The result follows if we prove that the intersections $\mathrm{cl}(h^{-1}(s) \cap B) \cap \mathrm{cl}(h^{-1}(s) \cap A)$ with $B\in\{\cU_j,\cW_j,\cV_j\colon 1 \leq j \leq 4\}$ define a bigon decomposition of the boundary sphere $\partial \mathrm{cl}(h^{-1}(s) \cap A)$. As the bigon decomposition of each set is defined by pairwise intersections, the Dihedral group $\langle \sigma, \tau \rangle$ will act preserving such bigon decompositions. Then, for Remarks \ref{invarianza} and \ref{permutaconos}, it is enough to prove the result for the balls $\mathrm{cl}(h^{-1}(s) \cap \cU_1)$ and $\mathrm{cl}(h^{-1}(s) \cap \cW_1)$.

\smallskip
\noindent
\textit{Bigon decomposition of $\partial \mathrm{cl}(h^{-1}(s)\cap\cU_1)$.} The base space of the $I$-bundle structure in $\mathrm{cl}(h^{-1}(s)\cap\cU_1)$ is the topological disc $\Omega_{\cU_s} = \{ z \in \C \colon |z| \leq 1, \ \Im(z) \geq s\}$ (Remark \ref{recapitulacion_fibrado1}), and the fibers are $C^1$-curves $\gamma_{z}$ parametrized by a rescaling of their arc-length (Figure \ref{Omega_U}). The induced projection
\[\Omega_{\cU_s} \times [0,1] \rightarrow\mathrm{cl}(h^{-1}(s)\cap\cU_1), \qquad (z,t) \rightarrow [0,1,\gamma_z(t),z],\]collapses the fibers at the points $q_1=e^{i\arcsin(s)},q_2=-e^{-i\arcsin(s)} \in \partial \Omega_{\cU_s}$. Let $\partial\Omega_{\cU_s}=\cL_s\cup\cA_s$ with $\cL_s=\{ z \in \Omega_{\cU_s}\colon \Im(z) = s \}$ and $\cA_s = \{ z \in \Omega_{\cU_s} \colon |z| = 1\}$, these are the segment and the arc sections of the boundary of the shaded region in Figure \ref{Omega_U}. Then the projection of the boundary $\partial(\Omega_{\cU_s} \times [0,1])=\Omega_{\cU_s}\times\{0,1\}\cup\partial\Omega_{\cU_s}\times[0,1]$ determines a bigon decomposition in the sphere $\partial\mathrm{cl}(h^{-1}(s)\cap\cU_1)$ with the four bigons:
\begin{multicols}{2}
\begin{itemize}
    \item $\cB_1$ the projection of $\Omega_{\cU_s} \times \{0\}$.
    \item $\cB_2$ the projection of $\cL_s\times[0,1]$
    \item $\cB_3$ the projection of $\Omega_{\cU_s} \times \{1\}$
    \item $\cB_4$ the projection of $\cA_s\times[0,1]$
\end{itemize}
\end{multicols}
\noindent
and with two vertices at $Q_1, Q_2$ (Remark \ref{recapitulacion_fibrado1}.$III$). From constructions in the proof of Lemma \ref{nivelconos}, a quadrilateral $Z\in\partial\mathrm{cl}(h^{-1}(s)\cap\cU_1)$ satisfy that (see Definition \ref{conos} and Figure \ref{Omega_U}):
\begin{itemize}
    \item if $Z\in\cB_1$, then $r(Z)=d(z_3,\overline{z_4z_1})$ and therefore, $Z\in\partial\mathrm{cl}(h^{-1}(s)\cap\cW_1)$.
    \item if $Z\in\cB_2$, then $r(Z)=d(z_4,\overline{z_1z_2})$ or $r(Z)=d(z_1,\overline{z_3z_4})$ and therefore, $Z\in\partial\mathrm{cl}(h^{-1}(s)\cap\cV_1)$.
    \item if $Z\in\cB_3$, then $\ell(Z)=|z_3-z_4|$ and therefore, $Z\in\partial\mathrm{cl}(h^{-1}(s)\cap\cV_3)$.
    \item if $Z\in\cB_4$, then $\ell(Z)=|z_1-z_4|$ and therefore, $Z\in\partial\mathrm{cl}(h^{-1}(s)\cap\cW_4)$.
\end{itemize}
\noindent
The sides of the bigons are (Figure \ref{edges_bigons}.$a)$):

\begin{itemize}
    \item $\cB_1\cap\cB_2 = h^{-1}(s)\cap\mathrm{cl}(\cU_1)\cap\mathrm{cl}(\cV_1)\cap\mathrm{cl}(\cW_1)$, the projection of $\cL_s \times \{0\}$.
    \item $\cB_2\cap\cB_3=h^{-1}(s)\cap\mathrm{cl}(\cU_1)\cap\mathrm{cl}(\cV_1)\cap\mathrm{cl}(\cU_3)\cap\mathrm{cl}(\cV_3)$, the projection of $\cL_s\times\{1\}$.
    \item $\cB_3\cap\cB_4= h^{-1}(s) \cap \mathrm{cl}(\cU_1) \cap \mathrm{cl}(\cV_3) \cap \mathrm{cl}(\cW_4)$, the projection of $\cA_s\times\{1\}$.
    \item $\cB_4\cap\cB_1= h^{-1}(s)\cap \mathrm{cl}(\cU_1)\cap \mathrm{cl}(\cV_4) \cap \mathrm{cl}(\cW_1) \cap \mathrm{cl}(\cW_4)$, the projection of $\cA_s\times\{0\}$.
\end{itemize}

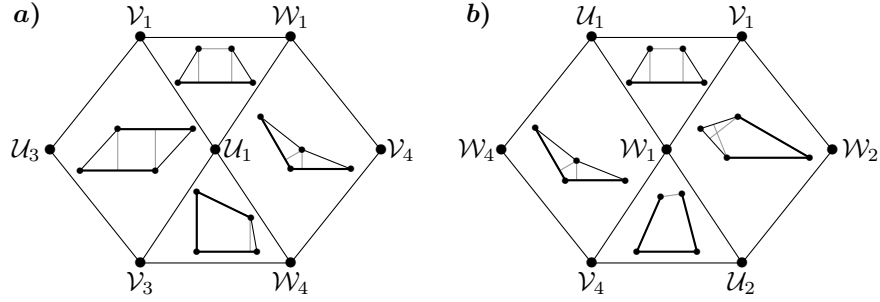
\begin{figure}
\begin{center}
\begin{tikzpicture}[scale=1]
\draw (-5.5,1.8) node {$\boldsymbol{a)}$};
\draw (-0.8,0)--(-2,1.5)--(-4,1.5)--(-5.2,0)--(-4,-1.5)--(-2,-1.5)--cycle;
\draw (-4,1.5)--(-2,-1.5);
\draw (-2,1.5)--(-4,-1.5);

\draw (-0.8,0) node {$\bullet$};
\draw (-0.85,0) node [right] {$\cV_4$};
\draw (-2,1.5) node {$\bullet$};
\draw (-2,1.5) node [above] {$\cW_1$};
\draw (-4,1.5) node {$\bullet$};
\draw (-4,1.5) node [above] {$\cV_1$};
\draw (-5.2,0) node {$\bullet$};
\draw (-5.2,0) node [left] {$\cU_3$};
\draw (-4,-1.5) node {$\bullet$};
\draw (-4,-1.5) node [below] {$\cV_3$};
\draw (-2,-1.5) node {$\bullet$};
\draw (-2,-1.5) node [below] {$\cW_4$};
\draw (-3,0) node {$\bullet$};
\draw (-3,0) node [right] {$\cU_1$};

\draw (-4.8,-0.27)--(-3.8,-0.27)--(-3.3,0.28)--(-4.3,0.28)--cycle;
\draw [black!40] (-3.8,-0.27)--(-3.8,0.28);
\draw [black!40] (-4.3,0.28)--(-4.3,-0.27);
\draw [thick=0.4] (-3.3,0.28)--(-4.3,0.28);
\draw [thick=0.4] (-4.8,-0.27)--(-3.8,-0.27);
\draw (-4.8,-0.27) node {\tiny{$\bullet$}};
\draw (-3.8,-0.27) node {\tiny{$\bullet$}};
\draw (-3.3,0.28) node {\tiny{$\bullet$}};
\draw (-4.3,0.28) node {\tiny{$\bullet$}};

\draw (-3.225,1.35)--(-3.5,0.9)--(-2.5,0.9)--(-2.785,1.35);
\draw [black!40] (-3.225,1.35)--(-3.225,0.9);
\draw [black!40] (-2.785,1.35)--(-2.785,0.9);
\draw [black!40] (-2.785,1.35)--(-3.225,1.35);
\draw [thick=0.4] (-3.5,0.9)--(-2.5,0.9);
\draw (-3.225,1.35) node {\tiny{$\bullet$}};
\draw (-3.5,0.9) node {\tiny{$\bullet$}};
\draw (-2.5,0.9) node {\tiny{$\bullet$}};
\draw (-2.785,1.35) node {\tiny{$\bullet$}};

\draw [shift={(-2,-0.25)}] (0.8,0)--(0,0)--(-0.4,0.69282)--(0.15, 0.259808)--cycle;
\draw [black!40,shift={(-2,-0.25)}] (0.15,0.259808)--(0.15,0);
\draw [black!40,shift={(-1.85,-0.0001)},rotate=-60] (0,0)--(0,-0.256);
\draw [thick=0.4,shift={(-2,-0.25)}] (0.82,0)--(0,0);
\draw [thick=0.4,shift={(-2,-0.25)}] (0,0)--(-0.4,0.69282);
\draw (-1.2,-0.25) node {\tiny{$\bullet$}};
\draw (-2,-0.25) node {\tiny{$\bullet$}};
\draw (-2.4,0.44282) node {\tiny{$\bullet$}};
\draw (-1.85,0.009808) node {\tiny{$\bullet$}};

\draw (-3.25,-1.35)--(-2.45,-1.35)--(-2.53,-0.9)--(-3.25,-0.55)--cycle;
\draw [black!40] (-2.54,-0.9)--(-2.54,-1.35);
\draw [thick=0.4] (-3.25,-1.35)--(-2.45,-1.35);
\draw [thick=0.4] (-3.25,-1.35)--(-3.25,-0.55);
\draw [thick=0.4] (-2.53,-0.9)--(-3.25,-0.55);
\draw (-3.25,-1.35) node {\tiny{$\bullet$}};
\draw (-2.45,-1.35) node {\tiny{$\bullet$}};
\draw (-2.53,-0.9) node {\tiny{$\bullet$}};
\draw (-3.25,-0.55) node {\tiny{$\bullet$}};

\draw (0.5,1.8) node {$\boldsymbol{b)}$};
\draw (0.8,0)--(2,1.5)--(4,1.5)--(5.2,0)--(4,-1.5)--(2,-1.5)--cycle;
\draw (4,1.5)--(2,-1.5);
\draw (2,1.5)--(4,-1.5);

\draw (0.8,0) node {$\bullet$};
\draw (0.85,0) node [left] {$\cW_4$};
\draw (2,1.5) node {$\bullet$};
\draw (2,1.5) node [above] {$\cU_1$};
\draw (4,1.5) node {$\bullet$};
\draw (4,1.5) node [above] {$\cV_1$};
\draw (5.2,0) node {$\bullet$};
\draw (5.2,0) node [right] {$\cW_2$};
\draw (4,-1.5) node {$\bullet$};
\draw (4,-1.5) node [below] {$\cU_2$};
\draw (2,-1.5) node {$\bullet$};
\draw (2,-1.5) node [below] {$\cV_4$};
\draw (3,0) node {$\bullet$};
\draw (3,0) node [left] {$\cW_1$};

\draw [shift={(6,0)}] (-3.225,1.35)--(-3.5,0.9)--(-2.5,0.9)--(-2.785,1.35);
\draw [black!40,shift={(6,0)}] (-3.225,1.35)--(-3.225,0.9);
\draw [black!40,shift={(6,0)}] (-2.785,1.35)--(-2.785,0.9);
\draw [black!40,shift={(6,0)}] (-2.785,1.35)--(-3.225,1.35);
\draw [thick=0.4,shift={(6,0)}] (-3.5,0.9)--(-2.5,0.9);
\draw [shift={(6,0)}] (-3.225,1.35) node {\tiny{$\bullet$}};
\draw [shift={(6,0)}] (-3.5,0.9) node {\tiny{$\bullet$}};
\draw [shift={(6,0)}] (-2.5,0.9) node {\tiny{$\bullet$}};
\draw [shift={(6,0)}] (-2.785,1.35) node {\tiny{$\bullet$}};

\draw [shift={(1.65,-0.4)}] (0.8,0)--(0,0)--(-0.4,0.69282)--(0.15, 0.259808)--cycle;
\draw [black!40,shift={(1.65,-0.4)}] (0.15,0.259808)--(0.15,0);
\draw [black!40,shift={(1.8,-0.141202)},rotate=-60] (0,0)--(0,-0.256);
\draw [thick=0.4,shift={(1.65,-0.4)}] (0.82,0)--(0,0);
\draw [thick=0.4,shift={(1.65,-0.4)}] (0,0)--(-0.4,0.69282);
\draw (2.45,-0.4) node {\tiny{$\bullet$}};
\draw (1.65,-0.4) node {\tiny{$\bullet$}};
\draw (1.25,0.29282) node {\tiny{$\bullet$}};
\draw (1.8,-0.141202) node {\tiny{$\bullet$}};

\draw [black!40] (2.92,-0.62)--(3.2,-0.58);
\draw [thick] (3.4,-1.35)--(2.6,-1.35);
\draw [thick] (2.6,-1.35)--(2.92,-0.62);
\draw [thick] (3.4,-1.35)--(3.2,-0.58);
\draw (2.92,-0.62) node {\tiny{$\bullet$}};
\draw (2.6,-1.35) node {\tiny{$\bullet$}};
\draw (3.4,-1.35) node {\tiny{$\bullet$}};
\draw (3.2,-0.58) node {\tiny{$\bullet$}};

\draw [black!40,shift={(4.9-0.952628,0.45)},rotate=-142] (0,0)--(0.48,0);
\draw [black!40,shift={(3.8,-0.1)},rotate=112] (0,0)--(.48,0);
\draw [thick=0.4] (3.8,-0.1)--(4.9,-0.1);
\draw [thick=0.4] (4.9,-0.1)--(4.9-.952628,0.45);
\draw (3.8,-0.1)--(4.9,-0.1)--(4.9-0.952628,0.45)--(4.9-1.44889,0.288229)--cycle;
\draw (3.8,-0.1) node {\tiny{$\bullet$}};
\draw (4.9,-0.1) node {\tiny{$\bullet$}};
\draw (4.9-0.952628,0.45) node {\tiny{$\bullet$}};
\draw (4.9-1.44889,0.288229) node {\tiny{$\bullet$}};
\end{tikzpicture}
\caption{Quadrilaterals belonging to the sides of the bigon decompositions on spheres $\partial \mathrm{cl}(h^{-1}(s)\cap \cU_1)$ in $a)$ and $\partial \mathrm{cl}(h^{-1}(s)\cap \cW_1)$ in $b)$. The sides of the quadrilateral with maximum length $\ell$ are the darkest ones and the length $r$ is represented with light gray segments.}
\label{edges_bigons}
\end{center}
\end{figure}

\noindent
\textit{Bigon decomposition of $\partial \mathrm{cl}(h^{-1}(s)\cap\cW_1)$}. Similarly to te previous case, the base of the $I$-bundle structure in $ \mathrm{cl}(h^{-1}(s)\cap\cW_1)$ is the disc $\Omega_{\cW_s}$ in Figure \ref{Omega_W}, and the fibers are $C^1$-curves $\gamma_z$ parametrized by a rescaling of their arc-length. The projection$$\Omega_{\cW_s}\times[0,1]\to\mathrm{cl}(h^{-1}(s)\cap\cW_1),\qquad (z,t)\to[0,1,\gamma_z(t),z]$$collapses the fibers of points in the set $\Gamma\cup\{q_1\}$. Let $\partial\Omega_{\cW_s}=\cA_s\cup\cL_s$ with $\cA_s = \{ z \in \Omega_{\cU_s} \colon |z| = 1\}$ and$$\cL_s=\Gamma\cup\{z \in \Omega_{\cW_s}\colon|z|=s\}\cup\{z\in \Omega_{\cW_s}\colon\Im(z)=s,\Re(z)\geq 0\},$$be a division of the boundary $\partial\Omega_{\cW_s}$ in two curves joining $q_1$ with $q_2$, see Figure \ref{Omega_W}. The four bigons in the sphere $\partial\mathrm{cl}(h^{-1}(s)\cap\cW_1)$ are:

\begin{itemize}
    \item $\cB_1$ the projection of $\{z \in \Omega_{\cW_s} : \Im(z) \geq s\} \times \{0\}$
    \item $\cB_2$ the projection of $\{z \in \Omega_{\cW_s} : \Im(z) \leq s\} \times \{0\} \cup \cL_s\times[0,1]$ 
\end{itemize}

\begin{multicols}{2}
\begin{itemize}
    \item $\cB_3$ the projection of $\Omega_{\cW_s}\times\{1\}$
    \item $\cB_4$ the projection of $\cA_s\times[0,1]$
\end{itemize}
\end{multicols}
\noindent
whose vertices are $Q_1, Q_2$. These bigons satisfy
\begin{multicols}{2}
\begin{itemize}
    \item $\cB_1=h^{-1}(s)\cap\mathrm{cl}(\cW_1)\cap\mathrm{cl}(\cU_1)$
    \item $\cB_2=h^{-1}(s)\cap\mathrm{cl}(\cW_1)\cap\mathrm{cl}(\cV_1)$
    \item $\cB_3=h^{-1}(s)\cap\mathrm{cl}(\cW_1)\cap\mathrm{cl}(\cU_2)$
    \item $\cB_4=h^{-1}(s)\cap\mathrm{cl}(\cW_1)\cap\mathrm{cl}(\cV_4)$
\end{itemize}
\end{multicols}
\noindent
and the sides of the bigons are (Figure \ref{edges_bigons}.$b)$):
\begin{itemize}
    \item $\cB_1\cap\cB_2= h^{-1}(s)\cap\mathrm{cl}(\cU_1)\cap\mathrm{cl}(\cV_1)\cap\mathrm{cl}(\cW_1)$, the projection of $\{z \in \Omega_{\cW_s} : \Im(z) = s \} \times \{0\}$.
    \item $\cB_2\cap\cB_3=h^{-1}(s)\cap\mathrm{cl}(\cW_1)\cap\mathrm{cl}(\cV_1)\cap\mathrm{cl}(\cU_2)\cap\mathrm{cl}(\cV_2)$, the projection of $\cL_s\times\{1\}$.
    \item $\cB_3\cap\cB_4= h^{-1}(s) \cap \mathrm{cl}(\cW_1) \cap \mathrm{cl}(\cU_2) \cap \mathrm{cl}(\cV_4)$, the projection of $\cA_s\times\{1\}$.
    \item $\cB_4\cap\cB_1= h^{-1}(s)\cap \mathrm{cl}(\cW_1)\cap \mathrm{cl}(\cV_4) \cap \mathrm{cl}(\cU_1) \cap \mathrm{cl}(\cW_4)$ is the projection of $\cA_s\times\{0\}$.
\end{itemize}
\end{proof}

\begin{obs}
The sides of the bigons in the decomposition of Lemma \ref{pegados} correspond to either repetition of three extreme values $\ell$ or $r$, obtained by intersecting three balls, or to repeat two times $\ell$ and two times $r$, obtained by intersecting four balls. This is captured in the diagram of Figure \ref{edges_bigons}.
\end{obs}

\begin{thm}
\label{levelspheres}
If $s\in(0,1)$, then the level set $h^{-1}(s)\subset S(4)$ is homeomorphic to the sphere $\S^3$.
\end{thm}

\begin{figure}[h]
\begin{center}
\begin{tikzpicture}
\draw (-3.3,2.6) node {$a)$};

\draw [black!40] (-3,-1.5)--(0,-1.5)--(0,1.5)--(-3,1.5)--cycle;
\draw [black!40] (1,-0.5)--(1,2.5)--(-2,2.5);
\draw [black!40,dashed] (-2,2.5)--(-2,-0.5)--(1,-0.5);
\draw [black!40,dashed] (-3,-1.5)--(-2,-0.5);
\draw [black!40] (0,-1.5)--(1,-0.5);
\draw [black!40] (0,1.5)--(1,2.5);
\draw [black!40] (-3,1.5)--(-2,2.5);

\draw [black!40] (-3,-1.5) node {\tiny{$\bullet$}};
\draw [black!40] (0,-1.5) node {\tiny{$\bullet$}};
\draw [black!40] (0,1.5) node {\tiny{$\bullet$}};
\draw [black!40] (-3,1.5) node {\tiny{$\bullet$}};

\draw [black!40] (1,-0.5) node {\tiny{$\bullet$}};
\draw [black!40] (1,2.5) node {\tiny{$\bullet$}};
\draw [black!40] (-2,2.5) node {\tiny{$\bullet$}};
\draw [black!40] (-2,-0.5) node {\tiny{$\bullet$}};

\draw [thick] (-1.5,-1.5)--(0,0)--(-1.5,1.5)--(-3,0)--cycle;
\draw [thick] (0,0)--(0.5,-1)--(1,1)--(0.5,2)--cycle;
\draw [thick] (-1.5,1.5)--(0.5,2)--(-0.5,2.5)--(-2.5,2)--cycle;
\draw (-1.5,-1.5) node {\small{$\bullet$}};
\draw (0,0) node {\small{$\bullet$}};
\draw (-1.5,1.5) node {\small{$\bullet$}};
\draw (-3,0) node {\small{$\bullet$}};
\draw (0.5,-1) node {\small{$\bullet$}};
\draw (1,1) node {\small{$\bullet$}};
\draw (0.5,2) node {\small{$\bullet$}};
\draw (-1.5,1.5) node {\small{$\bullet$}};
\draw (-0.5,2.5) node {\small{$\bullet$}};
\draw (-2.5,2) node {\small{$\bullet$}};

\draw [thick,dashed] (-1.5,-1.5)--(0.5,-1)--(-0.65,-0.5)--(-2.5,-1)--cycle;
\draw [thick,dashed] (-2.5,-1)--(-3,0)--(-2.5,2)--(-2,0.85)--cycle;
\draw [thick,dashed] (-0.65,-0.5)--(1,1)--(-0.5,2.5)--(-2,0.85)--cycle;
\draw (-1.5,-1.5) node {\small{$\bullet$}};
\draw (0.5,-1) node {\small{$\bullet$}};
\draw (-0.65,-0.5) node {\small{$\bullet$}};
\draw (-2,0.85) node {\small{$\bullet$}};
\draw (-2.5,-1) node {\small{$\bullet$}};

\draw (2.2,2.6) node {$b)$};

\fill [gray!50] (5.3,1.5)--(5.64,2)--(6.5,0.5)--(5.95,0.5);
\fill [gray!15] (4.62,0.5)--(5.3,1.5)--(5.95,0.5);

\fill [gray!50] (5.3,-1.5)--(5.504,-2.1)--(6.5,0.5)--(5.95,0.5);
\fill [gray!15] (4.62,0.5)--(5.3,-1.5)--(5.95,0.5);

\fill [gray!50] (5.3,-1.5)--(5.504,-2.1)--(3.628,-0.3)--(4,0);
\fill [gray!15] (4.62,0.5)--(4,0)--(5.3,-1.5);

\fill [gray!50] (5.3,1.5)--(4,0)--(3.5,-0.4)--(5.65,2.05);
\fill [gray!15] (4.62,0.5)--(5.3,1.5)--(4,0);

\draw [black!40,shift={(4,0)}] (-1.3,-1.5)--(1.3,-1.5)--(1.3,1.5)--(-1.3,1.5)--cycle;
\draw [black!40,shift={(4,0)}] (2.6,-0.5)--(2.6,2.5)--(0,2.5);
\draw [black!40,dashed,shift={(4,0)}] (0,-0.5)--(2.6,-0.5);
\draw [black!40,dashed,shift={(4,0)}] (0,-0.5)--(0,2.5);
\draw [black!40,shift={(4,0)}] (-1.3,1.5)--(0,2.5);
\draw [black!40,shift={(4,0)}] (1.3,-1.5)--(2.6,-0.5);
\draw [black!40,shift={(4,0)}] (1.3,1.5)--(2.6,2.5);
\draw [black!40,dashed,shift={(4,0)}] (-1.3,-1.5)--(0,-0.5);

\draw [thick] (5.3,1.5)--(5.65,2.05);
\draw [thick,dashed] (4.62,0.5)--(5.3,1.5);

\draw [thick] (5.3,-1.5)--(5.504,-2.1);
\draw [thick,dashed] (4.62,0.5)--(5.3,-1.5);

\draw [thick] (4,0)--(3.628,-0.3);
\draw [thick,dashed] (4.62,0.5)--(4,0);

\draw [thick] (5.95,0.5)--(6.5,0.5);
\draw [thick,dashed] (4.62,0.5)--(5.28,0.5);
\draw [thick,dashed] (5.32,0.5)--(5.95,0.5);

\draw [shift={(4,0)}] (0.62,0.5) node {\small{$\bullet$}};
\draw [shift={(4,0)}] (1.3,1.5) node {\small{$\bullet$}};
\draw [shift={(4,0)}] (1.95,0.5) node {\small{$\bullet$}};
\draw [shift={(4,0)}] (1.3,-1.5) node {\small{$\bullet$}};
\draw [shift={(4,0)}] (0,0) node {\small{$\bullet$}};
\end{tikzpicture}
\caption{Cuboctahedron graph $\cG_0$, constructed from the middle points of the regular cube. The right image is one of the cones, constructed using the dual structure of the graph which gives a bigon on the sphere $\S^3$ as a one-point compactification of $\R^3$.}
\label{Diagrama_pegados_tridimensional}
\end{center}
\end{figure}
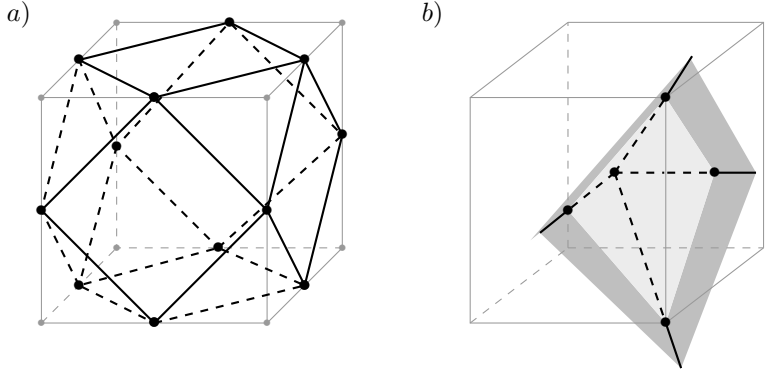

\begin{proof}
Let $\mathcal{C}_0\subset\mathbb{R}^3$ be the cube centered at $(0,0,0)$ and with vertices at $\{(\pm 1, \pm 1, \pm 1)\}$. The graph $\mathcal{G}_0$ obtained by joining the middle points of adjacent edges at each square face of $\cC_0$ is a regular cuboctahedron with vertices in the set
\[V=\{(\pm 1, \pm 1, 0),(\pm 1 , 0 ,\pm 1), (0,\pm 1, \pm 1)\},\]see Figure \ref{Diagrama_pegados_tridimensional}.$a)$. For every $v \in V$, $P_v = \{ w \in \mathbb{R}^3 : \langle v, w \rangle = 0\}$ denotes the orthogonal plane to $v$. These planes determine twelve closed discs $D_j$ in the sphere $\S^2\subset\R^3$, such that, for all $1\leq j\leq12$, the infinite cone region
\[C_j = \{ tw \colon w \in D_j, \ t \geq 0 \}\subset\R^3,\]
contains one and only one edge of the cube $\mathcal{C}_0$, see Figure \ref{Diagrama_pegados_tridimensional}.$b)$. Then, the one point compactification  $\R^3\cup\{\infty\}$ have a decomposition into twelve closed balls $C_j^*:=C_j \cup \{\infty\} \cong \D^3$ that are the one point compactifications of the cones $C_j$. This decomposition satisfies:

\begin{enumerate}
    \item There is exactly one ball $C_j^*$ for each vertex of the cuboctahedron $\mathcal{G}_0$.
    
    \item The intersection $C_j^*\cap C_k^*$ is a bigon if and only if the corresponding vertices of the graph  $\mathcal{G}_0$ share an edge. This bigon is contained in the compactification of a plane $P_v$ and its vertices are $\{0,\infty\}$.
    
    \item Let $E=\{(\pm 1 , \pm 1 , \pm 1)\} \cup \{(\pm 1, 0,0)\} \cup \{(0,\pm 1,0)\}\cup \{(0,0,\pm 1)\}\subset\cC_0$. The rays $\{t e\colon t\geq 0\}$ with $e\in E$, are the sides of the bigons $\{C_j^*\cap C_k^*\}$. Moreover, in one hand, the vertices of the cube $\{(\pm 1 ,\pm 1 ,\pm 1)\}$, correspond to the triangular faces of $\mathcal{G}_0$ and thus parameterize triple intersections $C^*_j \cap C^*_k \cap C^*_l$, and in the other hand, the centers of the faces in the cube $\{(\pm 1, 0,0)\} \cup \{(0,\pm 1,0)\}\cup \{(0,0,\pm 1)\}$, correspond to square faces of $\mathcal{G}_0$ and parameterize quadruple intersections $C^*_i \cap C^*_j \cap C^*_k \cap C^*_l$.

    \item The intersection of all the balls is
    \[  \bigcap_{j = 1}^{12} C_j^* = \{0,\infty\}.    \]
\end{enumerate}

\noindent
Summarizing: For all $1\leq j\leq 12$, the boundary of the closed ball $C_j^*\subset\R^3\cup\{\infty\}$ is decomposed into four bigons determined by the intersection with its neighboring balls, the intersection pattern is given by the cuboctahedron graph $\mathcal{G}_0$. Since $\cG_0$ is isomorphic to the graph $\cG$ in Figure \ref{Diagrama_pegados}, then the decomposition of $\R^3\cup\{\infty\}$ into closed balls is combinatorically isomorphic to the decomposition of the level sets $h^{-1}(s)$ given in Lemma \ref{pegados}, for all $0<s<1$. As the union of two discs through a disc at their boundaries is again a disc (Theorem \ref{thm_gluing_balls}), we conclude that $h^{-1}(s)$ is homeomorphic to the sphere $\S^3\cong\R^3\cup\{\infty\}$.
\end{proof}

\begin{obs}\label{polyhedron_structure}
It is possible to refine the decomposition of $h^{-1}(s)$, into the cones $\{\textrm{cl}(\cU_j), \textrm{cl}(\cV_j), \textrm{cl}(\cW_j)\}$, to a triangulation. One way to make this refinement explicit is by adding the faces of the cube in Figure \ref{Diagrama_pegados_tridimensional} and the planes defined by the diagonals so that the resulting pieces are tetrahedrons. For this reason, we take the liberty to call such balls, with the combinatorial structure given by Lemma \ref{pegados}, polyhedrons (even though the faces are bigons, not triangles).
\end{obs}

\begin{obs}\label{almost_conjugation_actions}
Conjugating the action of $\Sigma$ and $T$ on $\mathbb{C}^2$ (see Section \ref{features}) with the involution $(z_1,z_2)\mapsto(z_1,\overline{z_2}),$ we obtain $\hat{\Sigma}(z_1,z_2)=(-z_1,iz_2)$ and $\hat{T}(z_1,z_2)=(\overline{z_1},i\overline{z_2})$. Let $\varphi\colon\S^3\to\R^3$ denote the stereographic projection $(a+ib,c+id)\mapsto\frac{1}{1-b} (a,c,d)$. The compositions $\varphi\circ\hat{\Sigma}\circ\varphi^{-1}, \varphi\circ\hat{T}\circ\varphi^{-1}\colon\R^3\to\R^3$ are given by\[(x,y,z)\mapsto\frac{1}{x^2+y^2+z^2}(-x,-z,y)\quad\text{and}\quad (x,y,z)\mapsto\frac{1}{x^2+y^2+z^2}(x,z,y),\]respectively. The Dihedral group $\langle\varphi\circ\hat{\Sigma}\circ\varphi^{-1}, \varphi\circ\hat{T}\circ\varphi^{-1} \rangle$ acts in $\R^3$, permutes the planes $\{P_v\colon v \in V\}$, and preserve the decomposition $\mathbb{R}^3=\bigcup_{i = 1}^{12} C_i$ given in the proof of Theorem \ref{levelspheres}. Then, the identification $h^{-1}(s) \cong \mathbb{S}^3$ given by Theorem \ref{levelspheres} is equivariant under the corresponding actions of the Dihedral groups in the cones.
\end{obs}

\section{Cone structure on $S(4)$ and $\cM(4)$}
\label{special}
In this section we prove that $S(4)$ is homeomorphic to $\R^4$ by showing that the level sets of the height function $h\colon S(4)\to\R$ provide in $S(4)$ a structure of open cone over the sphere $\S^3$, we also describe the moduli space $\cM(4)$. It is worth mentioning that in \cite{AHT-JO} it is already proved that $S(4)\cong\R^4$, however, the result provided in this section is stronger because the cone structure in $S(4)$ naturally descends to an open cone structure of $\cM(4)=S(4)/\langle \sigma,\tau\rangle$. 

\begin{lem}
\label{eucliconos}
If $1\leq j\leq 4$, then the sets $\textrm{cl}(\cU_j), \textrm{cl}(\cV_j), \textrm{cl}(\cW_j)$ are homeomorphic to the open cone over the closed ball $C(\mathbb{D}^3)$. Moreover, under such homeomorphism, the height function $h$ corresponds to the function on $C(\D^3) \rightarrow (0,1]$ induced from the projection
    \[  \D^3 \times (0,1] \rightarrow (0,1].    \]
\end{lem}

\begin{proof} In order to prove the result for a set $\textrm{cl}(\cU) \in \{\textrm{cl}(\cU_j), \textrm{cl}(\cV_j), \textrm{cl}(\cW_j)\}$, we construct a family of homeomorphisms $f_s\colon h^{-1}(s) \cap \textrm{cl}(\cU) \rightarrow h^{-1}(1/2) \cap \textrm{cl}(\cU)$ that vary continuously on the parameter $s$ and a curve $\mu\colon(0,1) \rightarrow \mathrm{cl}(\cU)$ such that $h(\mu(t))$ is strictly decreasing (in fact, after a reparametrization of $\mu$, we can make $h(\mu(t)) = t$). As the curve $\mu$ is transversal to the level sets of the height function $h$, we can conclude that $\mathrm{cl}(\cU) \smallsetminus \{\mathfrak{R}_4\}$ has a structure of a fiber bundle with base $(0,1)$ and fiber $h^{-1}(1/2) \cap \textrm{cl}(\cU)$, homeomorphic to $\mathbb{D}^3$ by Lemma \ref{nivelconos}; this fiber bundle is trivial because it has a contractible base space, see \cite[Theorem 9.9]{HUS}). We can write down the trivialization as
    \[  h^{-1}\left(\frac{1}{2}\right)\cap \mathrm{cl}(\cU) \times (0,1) \rightarrow \mathrm{cl}(\cU) \smallsetminus \{\mathfrak{R}_4\}, \qquad (x,t) \mapsto f^{-1}_{h(\mu(t))}(x) = f^{-1}_{t}(x),\]
with continuous inverse $Z \mapsto (f_{h(Z)}(Z), h(Z))$. As the fibers $h^{-1}(s) \cap \textrm{cl}(\cU)$ collapse to $\{\mathfrak{R}_4\}$, when $s \rightarrow 1$, the product structure on $\mathrm{cl}(\cU) \smallsetminus \{\mathfrak{R}_4\}$ induces the structure of an open cone on $\mathrm{cl}(\cU)\subset S(4)$, where the height function $h$ corresponds to the projection onto the second factor.

Recall from Lemma \ref{nivelconos} that for all $s\in(0,1)$, the level set $h^{-1}(s) \cap \textrm{cl}(\cU)$ have the structure of a fiber bundle over a topological disc with some collapsed fibers at the boundary. Using this structure, we can construct the function $f_s : h^{-1}(s) \cap \textrm{cl}(\cU) \rightarrow h^{-1}(1/2) \cap \textrm{cl}(\cU)$ by giving first the homeomorphism between the base spaces and then extending such homeomorphism to the fibers as linear maps on the arc-length parameter defining them. If the homeomorphisms on the base spaces depend continuously on the parameter $s$, then the extensions on the fibers depend continuously on the parameter $s$ as well. 

If $\cU$ is either $\cU_1$ or $\cW_1$, then the base space of the fiber bundle structure of $\mathrm{cl}(h^{-1}(s) \cap \cU)$, is a topological disc contained in $\D^2$, which we can denote temporarily here by $\Omega_s$. We decompose the domain $\Omega_s$ into vertical intervals as follows: observe that for every $z \in \Omega_s$, the real part satisfies $-\cos(\arcsin(s)) \leq \Re(z) \leq \cos(\arcsin(s))$, and every $x$ in the interval $[-\cos(\arcsin(s)),\cos(\arcsin(s))]$, defines the vertical interval $I_x = \{y \in \R : x + i y \in \Omega_s \}$. Using this decomposition of $\Omega_s$, we can give the homeomorphism $f_s\colon\Omega_s\to \Omega_{1/2}$ as $x+iy\mapsto x'(x,s)+iy'(x,y,s)$, where \[x'(x,s)=\left(\frac{\cos(\arcsin(\frac{1}{2}))}{\cos(\arcsin(s))}\right)x,\]
and $y \mapsto y'(x,y,s)$ is the linear map between the intervals $I_x\subset\Omega_s$ and $I_{x'}\subset\Omega_{1/2}$. The expression for the imaginary part can be made explicit for $\cU_1$ as
\[ y'(x,y,s)=\left(\frac{\sqrt{1-(x')^2~}-s}{\sqrt{1-x^2~}-s}\right)(x-s)+s.\]
The corresponding function for $\cW_1$ has a much more involved expression, due to the complicated behaviour of the curve $\Gamma$ in the lower part of the boundary (determined by the points $q_\alpha$ in the proof of Lemma \ref{nivelconos}, see
Figure \ref{Omega_W}). This homeomorphism depends continuously on $s$ and thus, as discussed above, can be extended to the homeomorphism $f_s : h^{-1}(s) \cap \textrm{cl}(\cU) \rightarrow h^{-1}(1/2) \cap \textrm{cl}(\cU)$ which still depends continuously on the parameter $s$. The curves $\mu_1\colon[0,1)\to\mathrm{cl}(\cU_1)$ and $\mu_2\colon[0,1)\to\mathrm{cl}(\cW_1)$ given by
\[ \mu_1(t)=[0,1,1-\frac{t}{2}+i(1-t),i(1-\frac{t}{2})]~~\text{and}~~~\mu_2(t)=[0,1,(1-\frac{t}{2})(1+i),\frac{t}{2}+i(1-\frac{t}{2})],\] satisfy that $\mu_j(0)=\mathfrak{R}_4$, $h(\mu_j(t))=1-t$, and they are transversal to the level sets of the height function $h$. By the previous discussion, these curves give us the missing part for the construction of the homeomorphisms $f\colon \mathrm{cl}(\cU_1) \to C(\mathbb{D}^3)$ and $g\colon\mathrm{cl}(\cW_1)\to C(\D^3)$.

Finally, from Remark \ref{permutaconos}, it follows that for all $j\in\{1,2,3,4\}$, the homeomorphisms sought can be constructed of the form $f\circ\sigma^{-j}\colon\mathrm{cl}(\cU_j)\to C(\mathbb{D}^3)$, $f\circ\sigma^{-j}\circ\tau\colon \mathrm{cl}(\cV_j)\to C(\mathbb{D}^3),$ and $g\circ\sigma^{-j}\colon \mathrm{cl}(\cW_j)\to C(\mathbb{D}^3)$.
\end{proof}

The strategy in the proof of Lemma \ref{eucliconos}, can be used to prove the cone structure of $S(4)$ over $\S^3$. For this, we should be able to construct a family of functions $f_s\colon h^{-1}(s) \rightarrow h^{-1}\left(\frac{1}{2}\right)$ which depend continuously on the $s$ parameter. Such type of functions were constructed in Lemma \ref{eucliconos} for each piece $\mathrm{cl}(\cU_j),\mathrm{cl}(\cV_j),\mathrm{cl}(\cW_j)$ in the decomposition of the level set $h^{-1}(s)$, however they do not coincide in their intersections, and remains unsolved to construct the family $f_s\colon h^{-1}(s) \rightarrow h^{-1}\left(\frac{1}{2}\right)$. We overcome this difficulty by using the combinatorics of the well-understood decomposition of $S(4)$ into the cones $\mathrm{cl}(\cU)$. 

\begin{thm}\label{theo_cone_structure_simple_quadr}
There is a homeomorphism $S(4) \cong C(\S^3)\cong\R^4$, such that the height function $h$ corresponds to the function $C(\S^3) \rightarrow (0,1]$, induced from the projection $\S^3 \times (0,1] \rightarrow (0,1]$.
\end{thm}

\begin{proof}
We will proceed in a similar way as in the proof Theorem \ref{levelspheres}, by showing two combinatorially equivalent decompositions.

By taking the decomposition $\S^3=\cup_{j=1}^{12} C_j^*$ given in the proof of Theorem \ref{levelspheres} over each level sphere, we define decompositions of $\R^4 \cong C(\S^3)$ and $\S^4 \cong \Sigma(\S^3)$. In both cases, the combinatorics of the decomposition is codified by the graph $\mathcal{G}_0$, Figure \ref{Diagrama_pegados_tridimensional}.

The one point compactification $S(4)^*=S(4) \cup\{\infty\}$ satisfies that\[S(4)^*= \bigcup_{i=1}^4 \mathrm{cl}(\cU_i)^* \cup \mathrm{cl}(\cV_i)^* \cup \mathrm{cl}(\cW_i)^*,\]
where for every $\cU \in \{\cU_i, \cV_i, \cW_i\}$, $\mathrm{cl}(\cU)^* = \mathrm{cl}(\cU) \cup \{\infty\}$ is the one point compactification of $\mathrm{cl}(\cU)$. From Lemma \ref{eucliconos}, it follows that $\mathrm{cl}(\cU)^*$ is homeomorphic to the suspension $\Sigma(h^{-1}(1/2) \cap \mathrm{cl}(\cU)) \cong \Sigma(\D^3)\cong\D^4$. Since the graph $\cG$ (Figure \ref{Diagrama_pegados}) codifies the combinatorics of the decomposition of $h^{-1}(s)$ into the sets $\{h^{-1}(s) \cap \mathrm{cl}(\cU_i), h^{-1}(s) \cap \mathrm{cl}(\cV_i), h^{-1}(s) \cap \mathrm{cl}(\cW_i)\}$ (Lemma \ref{pegados}), then $\cG$ also codifies the decomposition of $S(4)^*$ determined by taking the suspension structure of each set $\mathrm{cl}(\cU)^*$.

We conclude that these decompositions of $S(4)^*$ and $\S^4$ into twelve copies of the closed disc $\D^4$ are combinatorically isomorphic, so that, $S(4)^*\cong\S^4$ and $S(4)$ is homeomorphic to $\R^4$. Moreover, each transversal section in $S(4)\cong C(\S^3)$ corresponds to a level set of the height function $h\colon S(4)\to\R$ and by Lemma \ref{eucliconos}, $h$ corresponds to the function on the cone $C(\S^3)$ corresponding to the projection $\S^3 \times (0,1] \rightarrow (0,1]$ under such identifications, thus the result follows.
\end{proof}

In the rest of this section we use the cone structure in $S(4)$ just proved, to describe the topology of the moduli space $\cM(4)$. By a combination of Theorem \ref{local} and Example \ref{eje4}, we can see that the quotients $h^{-1}(s)/\langle\sigma,\tau\rangle$ are homeomorphic to the sphere $\S^3$, for $s$ close enough to $1$. We prove this fact for every $s \in (0,1)$, because such construction will pay off when we study the boundary $h^{-1}(0) = \partial S(4)$.

\begin{lem}\label{Lema:level_sphere_moduli}
For every $s \in (0,1)$, $h^{-1}(s)/\langle\sigma,\tau\rangle$ is homeomorphic to $\S^3$.
\end{lem}

\begin{proof}
A fundamental domain in the level set $h^{-1}(s)$  for the action of the subgroup $\langle \sigma \rangle<\langle\sigma,\tau\rangle$ is the set
\[ \big[\mathrm{cl}(\cU_1) \cup \mathrm{cl}(\cV_1) \cup \mathrm{cl}(\cW_1) \big] \cap h^{-1}(s).  \]
The only reflection that leaves setwise fixed this subset is $\tau$, the action of which permutes $\mathrm{cl}(h^{-1}(s) \cap \cU_1)$ with $\mathrm{cl}(h^{-1}(s) \cap \cV_1)$ and leaves setwise fixed $\mathrm{cl}(h^{-1}(s) \cap \cW_1)$. Thus, the set
\[ B_s =\big\{[0,1,z_3,z_4] \in \mathrm{cl}(h^{-1}(s) \cap \cW_1) : \Im(z_3) \leq \Im(z_4)\big\}, \]
is a closed ball (Lemma \ref{nivelconos}), and the union $B_s\cup\mathrm{cl}(h^{-1}(s) \cap \cU_1)$ is a fundamental domain for the action of the dihedral group $\langle \sigma, \tau \rangle$.

If we take the identifications on the sphere $\partial \mathrm{cl}(h^{-1}(s) \cap \cU_1)$ determined by the actions of elements of the dihedral group $\langle\sigma,\tau\rangle$, we obtain a closed ball with boundary sphere splitted into two discs obtained from the intersections $\mathrm{cl}(h^{-1}(s) \cap \cU_1) \cap \mathrm{cl}(\cW_1)$ and $\mathrm{cl}(h^{-1}(s) \cap \cU_1) \cap \mathrm{cl}(\cW_4)$. If we do the analogous identifications on $\partial B_s$, we obtain a closed ball with boundary sphere splitted into two discs obtained from the intersections $\mathrm{cl}(h^{-1}(s) \cap \cW_1) \cap \mathrm{cl}(\cU_1)$ and $\mathrm{cl}(h^{-1}(s) \cap \cW_1) \cap \mathrm{cl}(\cV_4)$. As $\sigma^3 \tau (\textrm{cl}(\cW_1) \cap \textrm{cl}(\cV_4)) = \textrm{cl}(\cW_4) \cap \textrm{cl}(\cU_1)$, we conclude that after performing the identifications by elements of the dihedral group $\langle \sigma , \tau \rangle$ in the fundamental domain, we obtain two closed balls glued together through their boundary spheres. The space obtained in this way is homeomorphic to $\S^3$, finishing the proof.

As the  statements on the identifications on $\partial \mathrm{cl}(h^{-1}(s) \cap \cU_1)$ and $\partial B_s$ by elements of $\langle \sigma , \tau \rangle$ are not obvious, we make a detailed account of them separately.

\noindent
\textit{Identifications on $\partial \mathrm{cl}(h^{-1}(s) \cap \cU_1)$}. Recall from Lemma \ref{pegados} that $\mathrm{cl}(h^{-1}(s) \cap \cU_1)$ is homeomorphic to a closed ball with boundary sphere divided the bigons, 
\begin{enumerate}
    \item $\mathrm{cl}(h^{-1}(s) \cap \cU_1) \cap \mathrm{cl}(\cV_1)$, consisting of trapezoids with height $s$. 
    
    \item $\mathrm{cl}(h^{-1}(s) \cap \cU_1) \cap \mathrm{cl}(\cW_4)$, whose elements have two adjacent sides of maximal length.

    \item $\mathrm{cl}(h^{-1}(s) \cap \cU_1) \cap \mathrm{cl}(\cV_3)$, whose elements have two oposite sides of maximal length

    \item $\mathrm{cl}(h^{-1}(s) \cap \cU_1) \cap \mathrm{cl}(\cW_1)$, whose elements satisfy that $z_3$ lies in the bisector of the angle $\measuredangle z_4 z_1 z_2$ and it is at distance $s$ of both sides $\overline{z_1 z_2}$ and $\overline{z_1 z_4}$.
\end{enumerate}

As the geometric properties of the quadrilaterals at each intersection are mutually exclusive, we see that there are no identifications among different bigons under the action of the dihedral group $\langle\sigma,\tau\rangle$. We have two maps leaving fixed setwise two faces:
\[ \tau(\textrm{cl}(\cU_1) \cap \textrm{cl}(\cV_1)) = \textrm{cl}(\cV_1) \cap \textrm{cl}(\cU_1) \quad \textrm{and} \quad \sigma^2 \tau (\textrm{cl}(\cU_1) \cap \textrm{cl}(\cV_3)) = \textrm{cl}(\cV_3) \cap \textrm{cl}(\cU_1). \]
The map $\tau$ folds the face $\mathrm{cl}(h^{-1}(s) \cap \cU_1) \cap \textrm{cl}(\cV_1)$ in half along the curve of symmetric trapezoids
\[ \big\{ [0,1,1/2 + t+is, 1/2 - t+is]\colon s \leq 2 t \leq 1 \big\}. \]
Analogously, the map $\sigma^2 \tau$ is a reflection in $\mathrm{cl}(h^{-1}(s) \cap \cU_1) \cap \textrm{cl}(\cV_3)$ that swaps the angles $\measuredangle z_3z_40$ with $\measuredangle z_401$ and $\measuredangle z_4z_31$ with $\measuredangle 01z_3$, so it folds this face along the one-dimensional curve 
    \[  \left\{[0,1,z_3,z_4] : \arccos\left(\frac{1-s}{2}\right) \leq \measuredangle z_3z_40= \measuredangle z_401 \leq \pi/2 \right\},  \]
corresponding to quadrilaterals that are truncated isosceles triangles. The two curves, along which the maps $\tau$ and $\sigma^2 \tau$ fold the corresponding faces, are connected in the rectangle $[0,1,1+is,is]$, so that, the two identifications can be combined to be seen as folding in half the topological disc 
    \[  \Big[\mathrm{cl}(h^{-1}(s) \cap \cU_1) \cap \textrm{cl}(\cV_1) \Big] \bigcup \Big[\mathrm{cl}(h^{-1}(s) \cap \cU_1) \cap \textrm{cl}(\cV_3) \Big]. \] 
After doing such identifications, we are left with a space homeomorphic to a closed $3$-ball, whose boundary sphere is decomposed as the union of two topological discs, coming from the intersections $\mathrm{cl}(h^{-1}(s) \cap \cU_1) \cap \textrm{cl}(\cW_1)$ and $\mathrm{cl}(h^{-1}(s) \cap \cU_1) \cap \textrm{cl}(\cW_4)$.

\noindent
\textit{Identifications on $B_s$}. Because of the condition $\Im(z_3) \leq \Im(z_4)$, there are no longer elements with $\Im(z_4)<s$ in the closed ball $B_s\subset\mathrm{cl}(h^{-1}(s)\cap\cW_1)$ (see Figure \ref{Omega_W}). Thus, just as in the case of $\mathrm{cl}(h^{-1}(s) \cap \cU_1)$, an element $[0,1,z_3,z_4] \in B_s$ satisfies that
\[ z_4 \in \Omega_{\cU_s} = \{ z \in \D : \Im(z) \geq s\}. \]

\begin{figure}[h]
\begin{center}
\begin{tikzpicture}[scale=1]
\fill [gray!30]
(2.83,1)--(-2.83,1)--(2.83,1) arc (19.5:160.5:3);
\draw [black!50, dashed] (-3.5,0)--(6.5,0);
\draw [black!50, dashed] (0,-0.5)--(0,3.3);
\draw [black!50] (4,0) arc (0:180:1);
\draw [black!60] (2.83,1)--(-2.83,1);
\draw (3,0) arc (0:180:3);
\draw (6,0) arc (0:180:3);

\draw [black!30,shift={(2.225,1)},rotate=142] (0,0)--(1,0);
\draw (1.8,1.47) node {\footnotesize{$s$}};
\draw [black!60] (0,0)--(1.8,2); 
\draw [thick=0.8, shift={(2.225,1)}] (0,0)--(0.9,1); 
\draw [dotted, shift={(2.225,1)}] (0,0)--(1.8,2); 
\draw [black!50, dashed,shift={(1.8,2)}] (-0.5,0)--(1.8,0); 
\draw [dotted,shift={(1.8,2)},rotate=-33.35] (0,0)--(3.5,0);
\draw (1.8,2) node {\tiny{$\bullet$}};
\draw (1.65,2.15) node {\footnotesize{$z_4$}};
\draw [thick=0.8] (2.225,1)--(3,1);
\draw [thick=0.8] (3,1) arc (-270:-302:1);
\draw [thick=0.8,shift={(3.55,0.84)},rotate=-33.35] (0,0)--(0.9,0);
\draw (2.225,1) node {\tiny{$\bullet$}};
\draw (3.55,0.84) node {\tiny{$\bullet$}};
\draw (4.28,0.35) node {\tiny{$\bullet$}};
\draw (3.12,0.8) node {\footnotesize{$\gamma_{z_4}$}};

\draw [black!60] (0,0)--(-2,1.3);
\draw [thick=0.8, shift={(0.28,1)}] (0,0)--(-2/21,1.3/21);
\draw [dotted, shift={(0.28,1)}] (0,0)--(-2/21,1.3/21);
\draw [black!50, dashed,shift={(-2,1.3)}] (-0.5,0)--(2.4,0); 
\draw (0.28,1) node {\tiny{$\bullet$}};
\draw (-2,1.3) node {\tiny{$\bullet$}};
\draw (-1.9,1.5) node {\footnotesize{$z'_4$}};
\draw [thick=0.8] (0.28,1)--(0.98,1);
\draw [black!15,shift={(0.28,1)},rotate=-130] (0,0)--(1,0);
\draw (-0.23,0.55) node {\footnotesize{$s$}};
\draw (0.58,1.2) node {\footnotesize{$\gamma_{z'_4}$}};

\draw [black!60] (2.83,1) node {\tiny{$\bullet$}};
\draw [black!60] (-2.83,1) node {\tiny{$\bullet$}};
\draw (-3,0.03) node [below] {$-1$};
\draw (0,3) node {\tiny{$\bullet$}};
\draw (0.1,3.2) node [left] {$i$};
\draw (-3,0) node {\tiny{$\bullet$}};
\draw (0,0.03) node [below] {$0$};
\draw (0,0) node {\tiny{$\bullet$}};
\draw (3,0.03) node [below] {$1$};
\draw (3,0) node {\tiny{$\bullet$}};
\draw (6,0.03) node [below] {$2$};
\draw (6,0) node {\tiny{$\bullet$}};
\draw (-1.3,2.2) node {$\Omega_{\mathcal{U} s}$};

\end{tikzpicture}
\caption{Curves of possible values for the $z_3$ coordinate of a quadrialteral, with a fixed $z_4$ coordinate. These curves are cut from above with either the condition $\Im(z_4) \geq \Im(z_3)$, or $\ell_3=|z_3 - 1| \leq 1$.}
\label{Figure:I_bundle_moduli}
\end{center}
\end{figure}
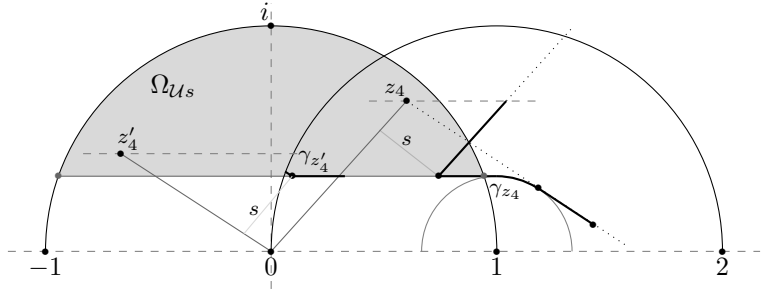

\begin{figure}[h]
\begin{center}
	\includegraphics[scale=0.24]{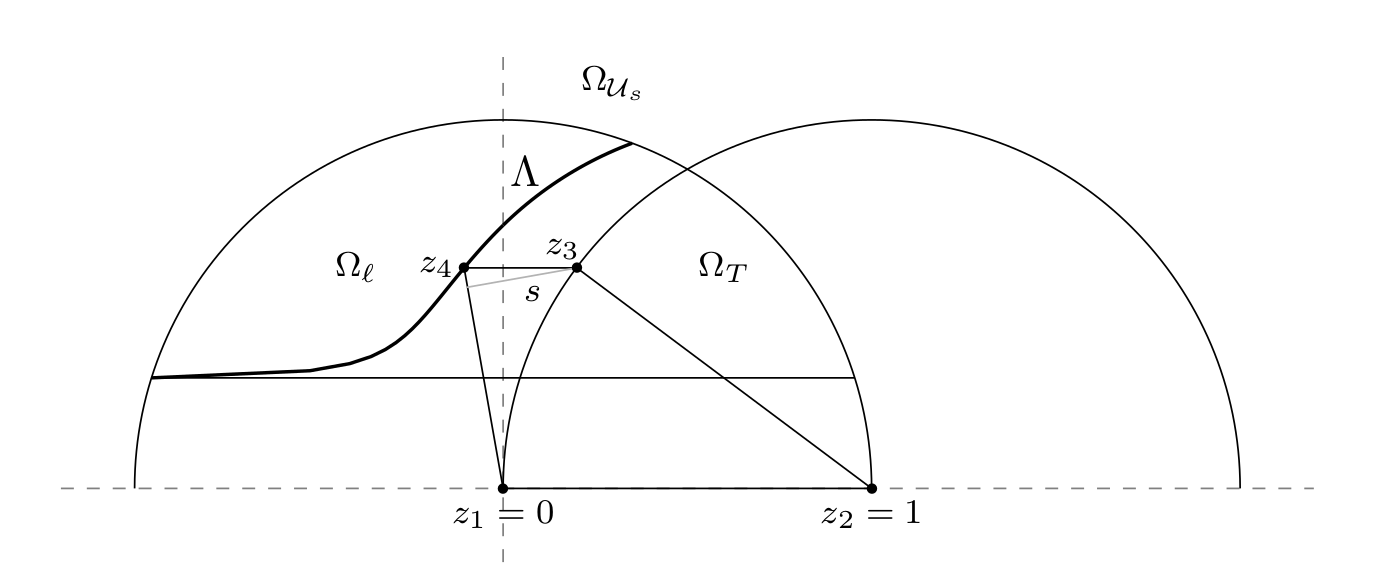}
\caption{A vertex $z_4 \in \Lambda$, satisfies that $[0,1,z_3,z_4] \in \textrm{cl}(h^{-1} \cap \cW_1)$ is a trapezoid and $|1 - z_3| = 1$. For $z_4\in\Omega_T$, we cut the fibers of the $I$-bundle at trapezoids, and for $z_4\in\Omega_{\ell}$, the fibers of the $I$-bundle are not cut by the new condition, so they are still cut by the condition $\ell_3=1= \ell_2$.}
\label{Figure:parallel_curve}
\end{center}
\end{figure}

There is a curve $\Lambda\subset \Omega_{\cU_s}$, which is determined by the property: $z_4 \in \Lambda$, if and only if, there is a $z_3 \in \{z\in\C\colon|z-1|=1\}$ such that\[[0,1,z_3,z_4]\in \mathrm{cl}(h^{-1}(s)\cap\cW_1), \quad \text{and}\quad\Im(z_3) =\Im(z_4), \]
see Figure \ref{Figure:parallel_curve}. The curve $\Lambda$ splits the set $\Omega_{\cU_s}$ into two connected components $\Omega_{\cU_s} = \Omega_\ell \cup \Omega_T$, so that, for $z_4  \in \Omega_\ell$, the fiber of the $I$-bundle of $\mathrm{cl}(h^{-1}(s) \cap \cW_1)$ is completely contained in $B_s$ and for $z_4 \in \Omega_T$, the fiber is cut at the point $z_3$, such that $[0,1,z_3,z_4]$ is a trapezoid, see Figure \ref{Figure:I_bundle_moduli} and Figure \ref{Figure:parallel_curve}. We have thus, a splitting of the boundary of $B_s$ into the following boundary components, determined by the domination of a constraint
\begin{enumerate}
    \item The intersection $\{[0,1,r+s+it,r+it] : 0 \leq r,t \leq 1\} \cap \mathrm{cl}(h^{-1}(s) \cap \cW_1)$, consisting of trapezoids with parallel sides of lengths $\ell = 1$ and $s$.

    \item $[0,1,z_3,z_4] \in  \mathrm{cl}(h^{-1}(s) \cap \cW_1) \cap \mathrm{cl}(\cV_4)$, such that $\Im(z_4) \geq \Im(z_3)$ (or equivalently $z_4 \in \Omega_T$) whose elements have two adjacent sides of maximal length $\ell=1$.  

    \item $[0,1,z_3,z_4] \in \mathrm{cl}(h^{-1}(s) \cap \cW_1) \cap \mathrm{cl}(\cU_2)$, such that $z_4 \in \Omega_\ell$, whose elements have two adjacent sides of maximal length. 

    \item $\mathrm{cl}(h^{-1}(s) \cap \cW_1) \cap\mathrm{cl}(\cU_1)$. 
\end{enumerate}

We can observe that for $Z = [0,1,z_3,z_4]$ belonging to the face described in the point 2, it is satisfied that $\tau(Z) = [0,1,z_3',z_4'] \in \mathrm{cl}(h^{-1}(s) \cap \cW_1) \cap \mathrm{cl}(\cU_2)$ with $z_4' \in \Omega_T$. On the other side, for $Z = [0,1,z_3,z_4]$ in the face described in point 3, it is satisfied that $\tau(Z)= [0,1,z_3',z_4']\in  \mathrm{cl}(h^{-1}(s) \cap \cW_1) \cap \mathrm{cl}(\cV_4)$ with $\Im(z_4') \leq \Im(z_3')$. Then these two faces of $B_s$ complement each other (under the action of $\tau$), to obtain a face equivalent to $\mathrm{cl}(h^{-1}(s) \cap \cW_1) \cap \mathrm{cl}(\cV_4)$ in $\mathrm{cl}(h^{-1}(s) \cap \cW_1)$.

It is immediate from the previous descriptions that the reflection $\tau$, acts as a folding map on the face of $B_s$ described in point 1, defined by trapezoids. After doing this folding in $B_s$, we obtain a space homeomorphic to a closed $3$-ball, whose boundary sphere is decomposed as the union of two topological discs, one coming from the intersection $\textrm{cl}(h^{-1}(s)\cap \cW_1) \cap \mathrm{cl}(\cU_1)$ and the other one, as a combination of faces described in points 2 and 3, which is a face equivalent to $\textrm{cl}(h^{-1}(s)\cap \cW_1) \cap \mathrm{cl}(\cV_4)$.

Finally, as we have that $\sigma^3 \tau(\cW_1 \cap \cV_4) =\cU_1\cap \cW_4$, then we conclude that the quotient $h^{-1}(s) / \langle \sigma, \tau\rangle$ is obtained by gluing together two copies of a closed $3$-ball by a homeomorphism between their boundaries
\[  \big( \mathrm{cl}(h^{-1}(s) \cap \cU_1) / \{\tau, \sigma^2 \tau\} \big) \cup \big(B_s / \{\tau, \sigma \} \big)\cong (\D^3 \cup \D^3) / \sim.  \]
This space is homeomorphic to the sphere $\S^3$.
\end{proof}

\begin{obs}\label{Remark:I_bundle_moduli}
The fundamental domain of the action of the dihedral group $\langle\sigma,\tau\rangle$ in $h^{-1}(s)$, constructed in Lemma \ref{Lema:level_sphere_moduli}, is characterized by two properties, on their elements as follows: $Z = [0,1,z_3,z_4]$ belongs to the fundamental domain if and only if
\begin{itemize}
    \item $\Im(z_4) \geq s$ and $|z_4| \leq 1$, and 

    \item $z_3$ belongs to a piecewise differentiable curve $\gamma_{z_4}$, which deppends continuously on the parameter $z_4$.
\end{itemize}
This is shown by two generic elements in Figure \ref{Figure:I_bundle_moduli}, and endows naturally with the structure of an $I$-bundle on the fundamental domain, over the topological disc $\Omega_{\cU_s}=\{z\in\D\colon\Im(z)\geq s\}$.
\end{obs}

\begin{thm}
There is a homeomorphism $\cM(4) \cong C(\S^3) \cong \R^4$.
\end{thm}
 
\begin{proof}
For every $0<s<1$, the fundamental domain $B_s\cup \mathrm{cl}(h^{-1}(s)\cap\cU_1)$ for the action of $\langle \sigma, \tau \rangle$ on the level set $h^{-1}(s)$, constructed in the proof of Lemma \ref{Lema:level_sphere_moduli}, is an $I$-bundle over the topological disc $\Omega_{\cU_s}=\{z\in\D\colon\Im(z)\geq s\}$, see Remark \ref{Remark:I_bundle_moduli} and Figure \ref{Figure:I_bundle_moduli}. By using the strategy used in Lemma \ref{eucliconos}, we can conclude that these fundamental domains can be uniformized with a continuous function as follows,
\[ \big(\D^2 \times [0,1] \big) \times (0,1) \rightarrow S(4), \qquad (w,t,s) \mapsto [0,1,\gamma_{\varphi_s(w)}(t),\varphi_s(w)]. \]
From Proposition \ref{globalmax}, it follows that this function can be extended by the constant with value $\mathfrak{R}_4\in S(4)$ in $\big(\D^3 \times [0,1] \big) \times \{1\}$.

Finally, it is possible to define an action of $\langle \sigma, \tau \rangle$ into $\S^3$ so that the uniformizing function $C(\S^3) \rightarrow S(4)$ is equivariant. By Lemma \ref{Lema:level_sphere_moduli}, we have $h^{-1}(s)/\langle\sigma,\tau\rangle\cong\S^3$ and so, the result follows.
\end{proof}

\section{Asymptotic description of $\overline{S(4)}$ and $\overline{\cM(4)}$}\label{section_asymptotic_desc}

As seen in the previous section, we have the cone structures $S(4) \cong C(\S^3) \cong \R^4$, but despite such well behaved structure, the closure $\overline{S(4)}$ is not even a manifold. In this section we describe precisely the region in $\partial S(4)$  which captures this bad asymptotic behaviour and prove that outside of such set, we have a very well behaved structure. As this section is harder to follow than the previous ones, we break down the ideas into four subsections.

Recall that for a subset $A \subset S(4) \subset \C\P^2$, we are considering two types of closures, as a subset of $\C\P^2$ denoted as $\overline{A}$, and as a subset of $S(4)$, denoted as $\textrm{cl}(A) = \overline{A} \cap S(4)$.

\subsection{Bad asymptotic behaviour at the boundary}\label{subsection_bad_asymptotic}
Here we describe the different shapes a quadrilateral can have in the boundary and show a bad asymptotic behaviour as the level sets $h^{-1}(s)$ approaches to $h^{-1}(0) = \partial S(4)$.

\begin{obs}
\label{cuadrilaterons_frontera}
A quadrilateral $Z=[z_1,z_2,z_3,z_4]\in\partial S(4)$ satisfies that there exists $j\in\{1,2,3,4\}$ such that $z_j\in\cL_j=\overline{z_{j+1}z_{j+2}}\cup\overline{z_{j+2}z_{j+3}}$ (Proposition \ref{zeroboundary}). There are four types of quadrilaterals with this property,
\begin{enumerate}[a)]
    \item Flags: The vertex $z_j$ belongs to $\cL_j\smallsetminus\{z_{j+1},z_{j+2},z_{j+3}\}$.
    \item Triangles: It is satisfied that $z_j=z_{j+1}$ or $z_j=z_{j+3}$.
    \item Wedges: It is satisfied that $z_j=z_{j+2}$. We denote the set of wedges as $\mathcal{K}$.
    \item 4-segments. $\cL_j$ degenerates to a segment and therefore, the four vertices are collinear.
\end{enumerate}
\end{obs}

Notice that in general, flags and wedges are limits of sequences of necessarily non-convex quadrilaterals.

\begin{defi}
\label{def:4-segments}
Let $Z\in\partial S(4)$ be a 4-segment.
\begin{enumerate}[(a)]
    \item We say that $Z$ is convex if it can be obtained as the limit of convex quadrilaterals.
    \item We say that $Z$ is non-convex if it cannot be obtained as limit of convex quadrilaterals. We denote as $\cN\subset\partial S(4)$ the closure of the space of non-convex $4$-segments.
\end{enumerate}
\end{defi}

We can observe that in a non-convex $4$-segment, the vertices must be in a zig-zag configuration, for example, the quadrilateral $\left[0,\frac{3}{4},\frac{1}{4},1\right]$.

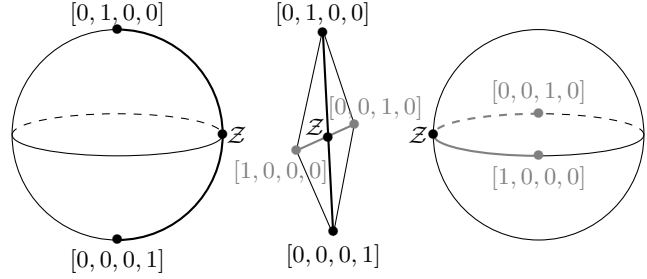
\begin{figure}[h]
\begin{center}
\begin{tikzpicture}[scale=0.7]
\draw (-4,2) arc (90:270:2);
\draw [thick] (-4,-2) arc (270:450:2);
\draw (-6,0) arc (180:360:2 and 0.4);
\draw [dashed] (-6,0) arc (-180:-360:2 and 0.4);

\draw (-4,2) node {\small{$\bullet$}};
\draw (-4,1.9) node [above] {\small{$[0,1,0,0]$}};
\draw (-2,0) node {\small{$\bullet$}};
\draw (-2.1,0) node [right] {\small{$\cZ$}};
\draw (-4,-2) node {\small{$\bullet$}};
\draw (-4,-2) node [below] {\small{$[0,0,0,1]$}};

\draw (-0.6,-0.3)--(-0.1,2)--(0.5,0.2)--(0.1,-1.9)--cycle;

\draw [black!50] (-0.6,-0.3) node {\small{$\bullet$}};
\draw [thick,black!50] (-0.9,-0.3) node [below] {\small{$[1,0,0,0]$}};

\draw (-0.1,1.95) node {\small{$\bullet$}};
\draw (-0.1,1.9) node [above] {\small{$[0,1,0,0]$}};

\draw [black!50] (0.5,0.2) node {\small{$\bullet$}};
\draw [thick,black!50] (0.9,0.1) node [above] {\small{$[0,0,1,0]$}};

\draw (0.1,-1.85) node {\small{$\bullet$}};
\draw (0.1,-1.9) node [below] {\small{$[0,0,0,1]$}};

\draw [thick,black!50] (-0.6,-0.3)--(0.5,0.2);
\draw [thick] (-0.1,2)--(0.1,-1.9);

\draw (0,-0.05) node {\small{$\bullet$}};
\draw (-0.25,0.15) node {\small{$\cZ$}};

\draw (4,0) circle (2);
\draw [thick,black!50] (2,0) arc (180:270:2 and 0.4);
\draw (4,-0.4) arc (270:360:2 and 0.4);
\draw [dashed,thick,black!50] (4,0.4) arc (90:180:2 and 0.4);
\draw [dashed] (6,0) arc (0:90:2 and 0.4);

\draw [thick,black!50] (4,-0.4) node {\small{$\bullet$}};
\draw [thick,black!50] (4,-0.4) node [below] {\small{$[1,0,0,0]$}};
\draw [thick,black!50] (4,0.4) node {\small{$\bullet$}};
\draw [thick,black!50] (4,0.4) node [above] {\small{$[0,0,1,0]$}};
\draw (2.08,0) node [left] {\small{$\cZ$}};

\draw [thick] (2,0) node {\small{$\bullet$}};
\end{tikzpicture}
\caption{$\cD_0$ as a wedge of two spheres, corresponding to the set of wedges $\mathcal{K}$, union with a square in the middle, corresponding to the set of non-convex 4-segments $\cN$. The union is through the diagonal segments of $\cN$ and the two spheres intersect at the center of the square that correspond to the zig-zag $\cZ=[0,1,0,1]$.}
\label{wedge_square}
\end{center}
\end{figure}

\begin{prop}\label{prop_topology_wedges_segments}
The set of wedges $\cK$ is homeomorphic to a wedge of two spheres $\S^2 \vee \S^2$. The space of non-convex $4$-segments $\cN$ is homeomorphic to a topological square and each sphere of $\mathcal{K}$ intersects $\mathcal{N}$ into one of the diagonals of the square. The intersection of the two diagonals of $\cN$, which is also the intersection of the two spheres of $\cK$, is the ``zig-zag'' $\cZ=[0,1,0,1]$, Figure \ref{wedge_square}.
\end{prop}

\begin{proof}
The set of wedges can be written as the union
\[  \cK = \{[0,1,0,w] \in \partial S(4)\colon w \in \mathbb{CP}^1\} \cup \{[0,1,z,1] \in \partial S(4) : z \in \mathbb{CP}^1 \}, \]
where both sets in the union are homeomorphic to $\S^2$ and intersect at $\{\cZ\}$, giving us the homeomorphism $\mathcal{K}\cong \S^2 \vee \S^2$. The set of non-convex $4$-segments can be triangulated as $\cN = \bigcup_{j = 0}^3 \sigma^j T$, with the simplex $T=\{ [0,1,t,x]\colon 0 \leq t\leq x\leq 1 \}.$ As there are four triangles in such description, we get a topological square, with diagonals
$$\{ [0,1,0,t] : 0 \leq t\leq 1 \} \cup \{ [0,t,0,1] : 0 \leq t\leq 1 \}, \text{~ and}$$ 
$$\{ [0,1,t,1] : 0 \leq t\leq 1 \} \cup \{ [0,1-t,-t,1-t] : 0 \leq t\leq 1 \},$$
which correspond to the intersection with each sphere of $\cK$, Figure \ref{wedge_square}. \end{proof}

We introduce a decomposition of the boundary into subsets analogous to the ones given in Definition \ref{conos} for a level set $h^{-1}(s)$, with $0 < s < 1$.

\begin{defi}\label{def:boundary_cones}
Let $j\in\{1,2,3,4\}$, then the sets $\widehat{\cU}_j, \widehat{\cV}_j,\widehat{\cW}_j \subset \partial S(4)$ consist of the points which are limit of quadrilaterals in $\mathrm{cl}(\cU_j)$,  $\mathrm{cl}(\cV_j)$, and $\mathrm{cl}(\cW_j)$, respectively.
\end{defi}

\begin{obs}
\label{permutaconos_frontera}
From Remark \ref{cuadrilaterons_frontera} it follows that every quadrilateral in $\partial S(4)$ belongs to one of the sets $\widehat{\cU}_j$, $\widehat{\cV}_j,$ or $ \widehat{\cW}_j$, with $j\in\{1,2,3,4\}$. This gives us the identity
\[h^{-1}(0)\cap \overline{S(4)}=\partial S(4)=\bigcup_{j = 1}^4 \left( \widehat{\cU}_j \cup  \widehat{\cV}_j \cup \widehat{\cW}_j \right).\]
The dihedral group $\langle\sigma,\tau\rangle$ permutes the sets $\widehat{\cU}_j$, $\widehat{\cV}_j,$ and $ \widehat{\cW}_j$ just like permutes the sets $\cU_j$, $\cV_j,$ and $\cW_j$ in the level sets $h^{-1}(s)\cap S(4)$, with $0<s<1$, see Remark \ref{permutaconos}.
\end{obs}

In order to understand the sets $\widehat{\cU}_j$, $\widehat{\cV}_j,$ and $ \widehat{\cW}_j$, for $j\in\{1,2,3,4\}$, we proceed to describe $\widehat{\cU}_1$ and $\widehat{\cW}_1$.

\begin{con}
\label{const}
We construct $\widehat{\cU}_1$ and $\widehat{\cW}_1$, the other sets can be obtained by using the action of the dihedral group $\langle\sigma,\tau\rangle$.

\noindent
\textit{$a)$ Construction of $\widehat{\cU}_1$}. From Lemma \ref{nivelconos} we know that a quadrilateral $Z(s) =[0,1,z_3(s),z_4(s)]\in\textrm{cl}(h^{-1}(s) \cap \cU_1)$, satisfies that $z_4(s)$ lies in the set $\Omega_{\cU s}$, and $z_3(s)$ lies in the $C^1$-curve $\gamma_{z_4(s)}$ constructed as a concatenation of an arc of a circle and two segments tangent to it, as shown in Figure \ref{Omega_U}. These conditions become as $s \rightarrow 0$, to the fact that a quadrilateral $Z=[0,1,z_3,z_4]\in\widehat{\cU}_1$ must satisfy that $z_4$ lies in the closed set $\Omega_{\cU_0}=\{z \in\mathbb{D}^2\colon\Im(z) \geq 0\}$, and the set of possible values for $z_3$ is the limit $\lim_{s\to 0}\gamma_{z_4(s)}$. For $z_4\in \Omega_{\cU_0}\smallsetminus\{1\}$, the limit $\lim_{s\to 0}\gamma_{z_4(s)}$ is a curve $\gamma_{z_4}$ obtained by cutting the union $\overline{01}\cup\{ (1-z_4)t+1\colon t\geq 0\}$ with the constraint $|z_3-z_4| \leq 1$, Figure \ref{ball_boundary}.$a)$. In the remaining case the quadrilaterals are wedges with $z_4 = z_2=1$, in this case the limit $\lim_{s\to 0}\gamma_{z_4(s)}$ behaves somewhat pathologically, because it is the two-dimensional set
\[\overline{01}\cup\big\{z+1\in\C\colon |z|\leq 1, Re(z)\geq 0, \Im(z)\leq 0\big\}.\]

\noindent
\textit{$b)$ Construction of $\widehat{\cW}_1$}. A quadrilateral $Z(s) = [0,1,z_3(s),z_4(s)] \in \textrm{cl}(h^{-1}(s) \cap \cW_1)$ is determined by $z_4(s) \in \Omega_{\cW s}$ and $z_3(s)$ belonging to a $C^1$-curve $\gamma_{z_4(s)}$ which is again the union of an arc and two segments tangent to it, see Figure \ref{Omega_W}. As $s \rightarrow 0$, these conditions translate to the fact that a quadrilateral $Z = [0,1,z_3,z_4] \in \widehat{\cW}_1$, must satisfy that $z_4$ lies in the set $\Omega_{\cW_0} = \{z \in \mathbb{D}^2: \Im(z) \geq 0\}$ (observe $\Omega_{\cW_0}=\Omega_{\cU_0}$ and it is better behaved than $\Omega_{\cW s}$), and $z_3$ lies in the limit $\lim_{s\to 0}\gamma_{z_4(s)}$. For every $z_4\in\Omega_{\cW_0} \smallsetminus \{1\}$, the set of possible values for $z_3$ so that $Z =[0,1,z_3,z_4] \in \widehat{\cW}_1$, is the curve $\gamma_{z_4}$ contained in the union $\overline{0z_4}\cup\{z_4+t(z_4-1)\colon t\geq 0\}$ and cut with the constraint $|z_3-1|\leq 1$. The case $z_4=z_2=1$, which corresponds to wedges, is again pathological because the set of possible values for $z_3$ is
\[\overline{01}\cup\{z+1\in\C\colon |z|\leq1, Re(z)\geq 0,Im(z)\geq 0\}.\]
\end{con}

\begin{lem}\label{segment_wedge_characterization}
The set of non-convex $4$-segments can be written as the union
\[  \mathcal{N} = \bigcup_{j = 1}^4 \left(\widehat{\cW}_j \cap \widehat{\cU}_j \cap \widehat{\cV}_j\right),  \]
and the set of wedges $\mathcal{K}$ is the closure of 
\[   \bigcup_{j = 1}^4 \left[ \left(\widehat{\cW}_j \cap \widehat{\cU}_j \right) \cup \left(\widehat{\cW}_j \cap \widehat{\cV}_j \right) \right] \smallsetminus \left(\widehat{\cW}_j \cap \widehat{\cU}_j \cap \widehat{\cV}_j\right).    \]
\end{lem}

\begin{proof}
We proceed by using the same strategy as in previous lemmas, \textit{i.e.} we describe the sets of $4$-segments and wedges contained in $\widehat{\cU}_1$, $\widehat{\cV}_1$ and $\widehat{\cW}_1$, and then using the dihedral group $\langle\sigma,\tau\rangle$, we obtain the description of the whole spaces $\cN$ and $\cK$.

Notice from Construction \ref{const}$.a)$ that a quadrilateral $Z \in \widehat{\cU}_1$ is a $4$-segment if and only if $Z = [0,1,s,t]$, with $-1\leq t \leq 1$ and $0 \leq s \leq 2$. If $Z$ is non-convex, then the inequality $s \leq t$ is satisfied. Also from Construction \ref{const}$.b)$, we can see that the same conditions must hold for the non-convex $4$-segments in $\widehat{\cW}_1$, giving us $\widehat{\cU}_1\cap\cN=\widehat{\cW}_1 \cap\cN$. Now, the set $\widehat{\cU}_1 \cap \cN=\{[0,1,s,t]\colon 0 \leq s \leq t \leq 1 \}$ is a $2$-simplex, which is invariant under the homeomorphism $\tau ([0,1,s,t])=[0,1,1-t,1-s]$, and as we have $\tau(\widehat{\cU}_1 \cap\cN)=\widehat{\cV}_1\cap\cN$, then $\widehat{\cU}_1\cap\cN=\widehat{\cV}_1\cap\cN=\widehat{\cW}_1\cap\cN=\widehat{\cW}_1 \cap \widehat{\cU}_1\cap\widehat{\cV}_1$. If we act with the homeomorphism $\sigma$, we get the equivalent conditions for all $j\in\{1,2,3,4\}$
\[\widehat{\cU}_j\cap\cN= \widehat{\cV}_j\cap\cN=\widehat{\cW}_j\cap\cN=\widehat{\cW}_j\cap\widehat{\cU}_j \cap \widehat{\cV}_j,\]
so in particular $\cN= \bigcup_{j = 1}^4 \left(\widehat{\cW}_j \cap \widehat{\cU}_j \cap \widehat{\cV}_j\right)$. As the intersection $\widehat{\cW}_j \cap \widehat{\cU}_j \cap \widehat{\cV}_j$ is a $2$-simplex, then the union of such intersections defines a triangulation of $\cN$.

From Construction \ref{const}$.a)$ we can see that the wedges in $\widehat{\cU}_1$ are of the form $[0,1,0,z_4]$ and the wedges in $\widehat{\cV}_1$ are of the form $[0,1,z_3,1]$, which are obtained by applying the refection $\tau$ to the wedges in $\widehat{\cU}_1$. In both cases, the wedges belong to the set $\widehat{\cW}_1$ as well. Moreover, from the Construction \ref{const}$.b)$ and using the previous case of non-convex 4-segments, we can see that the intersection $\mathcal{K}\cap \widehat{\cW}_1$ is the closure of
    \[  \left[ \left(\widehat{\cW}_1 \cap \widehat{\cU}_1 \right) \cup \left(\widehat{\cW}_1 \cap \widehat{\cV}_1 \right) \right] \setminus \left(\widehat{\cW}_1 \cap \widehat{\cU}_1 \cap \widehat{\cV}_1\right). \]
If we act with the homeomorphism $\sigma$, we obtain the description of $\mathcal{K} \cap \widehat{\cW}_j$, for every $j \in \{1,2,3,4\}$, in terms of double and triple intersections as in the case $j = 1$, and the result follows.
\end{proof}

The previous results motivates the following definition.

\begin{defi}\label{separating_spheres}
We define the set
\[\mathcal{D} := \bigcup_{j = 1}^4 \Big[ \big(\overline{\cW_j} \cap \overline{\cU_j} \big) \cup \big(\overline{\cW_j} \cap \overline{\cV_j} \big) \Big]\subset\overline{S(4)},\]
and for every $s \in [0,1]$, $\cD_s := \cD \cap h^{-1}(s)$.
\end{defi}

Observe that $\cD_0=\bigcup_{j = 1}^4 \big[\big(\widehat{\cW}_j \cap \widehat{\cU}_j \big) \cup \big(\widehat{\cW}_j \cap \widehat{\cV}_j \big) \big]= \cK \cup \cN$ and $\cD_1 =\{\mathfrak{R}_4\}$. We proceed to describe the topology of the sets $\cD_s$.

\begin{prop}\label{prop_limit_spheres}
For every $0 < s < 1$, it is satisfied that $\cD_s \cong \S^2$, and
    \[  \lim_{s \rightarrow 1} \cD_s = \cD_1 = \{\mathfrak{R}_4\}, \qquad  \lim_{s \rightarrow 0} \cD_s = \cD_0\cong (\S^2 \vee \S^2) \cup \D^2.  \]
Thus, $\cD$ is homeomorphic to an open cone $C(\S^2) \cong \R^3$ compactified with $\cD_0$.
\end{prop}

\begin{proof}
Quadrilaterals in $\cD_s$ have a vertex $z_j$ such that $d(z_j,\overline{z_{j+1}z_{j+2}})=s=d(z_j,\overline{z_{j+2}z_{j+3}})$, and therefore, $z_j$ belongs to the bisector of the angle $\measuredangle z_{j+1}z_{j+2}z_{j+3}$. This fact is used to prove that $\lim_{s \rightarrow 0} \mathcal{D}_s=\mathcal{D}_0=\cK\cup\cN$ as follows:

\begin{itemize}
    \item $\mathcal{D}_0 \subset \lim_{s \rightarrow 0} \mathcal{D}_s.$ A non-convex $4$-segment $Z= [0,1,x,y] \in \widehat{\cU}_1 \cap \widehat{\cW}_1$ with $0\leq x \leq y \leq 1$, can be approximated by the quadrilaterals$$Z_{\theta}=\left[0,1,\frac{x}{\cos(\theta/2)} e^{i \theta/2}, \frac{y}{\cos(\theta)} e^{i \theta} \right]\in\Big(\textrm{cl}(\cU_1) \cap \textrm{cl}(\cW_1)\Big),$$satisfying $h(Z_{\theta}) = x \tan(\theta/2) \rightarrow 0$ as $\theta \rightarrow 0$. Thus, $Z \in \lim_{s \rightarrow 0} \mathcal{D}_s$. Analogously, a wedge $Z = [0,1,0,r e^{i \beta}] \in \widehat{\cU}_1 \cap \widehat{\cW}_1$ belongs to $\lim_{s \rightarrow 0} \mathcal{D}_s$ because it can be approximated by the sequence of quadrilaterals $Z_t=[0,1,t e^{i \beta/2},r e^{i \beta}]\in\textrm{cl}(\cU_1) \cap \textrm{cl}(\cW_1)$, satisfying $h(Z_t)= t \sin(\beta/2) \rightarrow 0$, as $t \rightarrow 0$. The statement for every other element in $\cD_0$ can be proved by using the action of the dihedral group $\langle\sigma,\tau\rangle$.
    
    \item $\lim_{s \rightarrow 0} \mathcal{D}_s\subset \mathcal{D}_0.$ Observe that if $Z = [z_1,z_2,z_3,z_4] \in \partial S(4)$ with $z_j\in\cL_j$ is a flag, a triangle or a convex $4$-segment (\textit{i.e.} $Z$ is not in $\cD_0$), see Remark \ref{cuadrilaterons_frontera} and Definition \ref{def:4-segments}, and $\{Z_n = [z_1(n), z_2(n), z_3(n), z_4(n)] \}_{n \in \mathbb{N}}$ is a sequence of simple quadrilaterals which converge to $Z$, then for a sufficiently large subscript $n$, the quadrilateral $Z_n$ have the property that the vertex $z_j(n)$ is far from the bisector of the internal angle $\measuredangle z_{j+1}(n)z_{j+2}(n)z_{j+3}(n)$. Thus, $Z \notin \lim_{s \rightarrow 0} \cD_s.$ 
\end{itemize}

\begin{figure}[h]
\begin{center}
\begin{tikzpicture}
\draw [black!20] (-3,-1.5)--(0,-1.5)--(0,1.5)--(-3,1.5)--cycle;
\draw [black!20] (1,-0.5)--(1,2.5)--(-2,2.5);
\draw [black!20,dashed] (-2,2.5)--(-2,-0.5)--(1,-0.5);
\draw [black!20,dashed] (-3,-1.5)--(-2,-0.5);
\draw [black!20] (0,-1.5)--(1,-0.5);
\draw [black!20] (0,1.5)--(1,2.5);
\draw [black!20] (-3,1.5)--(-2,2.5);

\draw [black!40] (-3,-1.5) node {\tiny{$\bullet$}};
\draw [black!40] (0,-1.5) node {\tiny{$\bullet$}};
\draw [black!40] (0,1.5) node {\tiny{$\bullet$}};
\draw [black!40] (-3,1.5) node {\tiny{$\bullet$}};
\draw [black!40] (1,-0.5) node {\tiny{$\bullet$}};
\draw [black!40] (1,2.5) node {\tiny{$\bullet$}};
\draw [black!40] (-2,2.5) node {\tiny{$\bullet$}};
\draw [black!40] (-2,-0.5) node {\tiny{$\bullet$}};

\draw [black!40] (-1.5,-1.5)--(0,0)--(-1.5,1.5)--(-3,0)--cycle;
\draw [black!40] (0,0)--(0.5,-1);
\draw [black!40] (1,1)--(0.5,2);
\draw [black!40] (-2.5,2)--(-1.5,1.5);
\draw [black!40] (-0.5,2.5)--(0.5,2);

\draw (-0.5,2.5) node {\tiny{$\bullet$}};
\draw (-0.5,2.5) node [above] {\small{$\cU_2$}};

\draw (-2.5,2) node {\tiny{$\bullet$}};
\draw (-2.5,2) node [left] {\small{$\cW_2$}};

\draw (0.5,2) node {\tiny{$\bullet$}};
\draw (0.5,2) node [right] {\small{$\cW_1$}};

\draw (-1.5,1.5) node {\tiny{$\bullet$}};
\draw (-1.5,1.5) node [above] {\small{$\cV_1$}};

\draw (1,1) node {\tiny{$\bullet$}};
\draw (1,1) node [right] {\small{$\cV_4$}};

\draw (-2,0.85) node {\tiny{$\bullet$}};
\draw (-2,0.85) node [right] {\small{$\cV_2$}};

\draw (-3,0) node {\tiny{$\bullet$}};
\draw (-3,0) node [left] {\small{$\cU_3$}};

\draw (0,0) node {\tiny{$\bullet$}};
\draw (0,0) node [right] {\small{$\cU_1$}};

\draw (-0.65,-0.5) node {\tiny{$\bullet$}};
\draw (-0.65,-0.4) node [above] {\small{$\cU_4$}};

\draw (0.5,-1) node {\tiny{$\bullet$}};
\draw (0.5,-1) node [right] {\small{$\cW_4$}};

\draw (-2.5,-1) node {\tiny{$\bullet$}};
\draw (-2.5,-1) node [below] {\small{$\cW_3$}};

\draw (-1.5,-1.5) node {\tiny{$\bullet$}};
\draw (-1.5,-1.5) node [below] {\small{$\cV_3$}};

\draw [black!40,dashed] (-1.5,-1.5)--(0.5,-1);
\draw [black!40,dashed] (-2.5,-1)--(-0.65,-0.5);

\draw [black!40,dashed] (-2.5,2)--(-3,0);
\draw [black!40,dashed] (-2.5,-1)--(-2,0.85);

\draw [black!40,dashed] (-0.65,-0.5)--(1,1)--(-0.5,2.5)--(-2,0.85)--cycle;

\draw [thick] (-2,2.5)--(0,1.5)--(1,-0.5);
\draw [thick,dashed] (-2,2.5)--(-3,-1.5)--(1,-0.5);
\draw (-2,2.5) node {$\bullet$};
\draw (0,1.5) node {$\bullet$};
\draw (1,-0.5) node {$\bullet$};
\draw (-3,-1.5) node {$\bullet$};

\end{tikzpicture}
\caption{For $0<s<1$, the set $\mathcal{D}_s\subset h^{-1}(s)$ is homeomorphic to the suspension over the darkest closed curve that separates the cube. Here, we are removing the edges of the cuboctahedron corresponding to intersections $\overline{\cW_j} \cap \overline{\cU_j}$ and $\overline{\cW_j} \cap \overline{\cV_j}$ in $h^{-1}(s)$.} 
\label{Ds_separating_curve}
\end{center}
\end{figure}
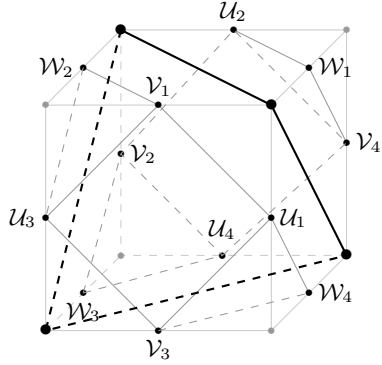

From the decomposition of $h^{-1}(s)$ into the sets $\{\textrm{cl}(\cU_j),\textrm{cl}(\cV_j),\textrm{cl}(\cW_j)\}$ given in Theorem \ref{levelspheres} and represented in Figure \ref{Diagrama_pegados_tridimensional}, we can see that the set $\mathcal{D}_s$ is the suspension over a closed curve in the regular cube shown in Figure \ref{Ds_separating_curve}, and so, $\cD_s$ is homeomorphic to $\mathbb{S}^2$. Thus, the homemomorphism $\cD \cap S(4) \cong C(\S^2)$ follows from the cone structure on $S(4)$ given by Theorem \ref{theo_cone_structure_simple_quadr}.
\end{proof}

Proposition \ref{prop_limit_spheres} justifies the assertion that the quadrilaterals in $\cD_0$ generate bad asymptotic behaviour.

\subsection{Good behaviour at the boundary}\label{subsection_good_behaviour}

As observed in the previous section, there is bad asymptotic behaviour at the set $\cD_0=\cN\cup\mathcal{K}$, as $\cD_s \rightarrow \cD_0$ and $\cD_s \ncong \cD_0$. In this section, we prove that the set $\cD_0$ captures all the bad asymptotic behaviour by proving that the topology of $h^{-1}(s) \smallsetminus \cD_s$ doesn't change as $s \rightarrow 0$.











To prove this, we procede as in the proof of Lemma \ref{nivelconos}, describing first the topology of the sets $\widehat{\cU}_j,\widehat{\cV}_j,\widehat{\cW}_j\subset\partial S(4)$, by identifying a structure of an $I$-bundle on them. We note that this structure exists only outside of the set of wedges and non-convex $4$-segments $\cD_0$.

\begin{lem}\label{lema:boundary_complement_of_wedges}
The sets $\widehat{\cU}_j\smallsetminus\cD_0$, $\widehat{\cV}_j\smallsetminus\cD_0$, and $\widehat{\cW}_j\smallsetminus\cD_0$, are homeomorphic to a closed 3-ball with a closed 2-disc removed from their boundary.
\end{lem}

\begin{proof}
We describe the $I$-bundle structure in the sets $\widehat{\cU}_j\smallsetminus\cD_0$, $\widehat{\cV}_j\smallsetminus\cD_0$, $\widehat{\cW}_j\smallsetminus\cD_0$ as the limit of the $I$-bundle structures in $\textrm{cl}(\cU_j)$, $\textrm{cl}(\cV_j)$, $\textrm{cl}(\cW_j),$ given in Lemma \ref{nivelconos}. From Remark \ref{permutaconos_frontera}, it follows that we only need to describe the $I$-bundle structure for the sets $\widehat{\cU}_1\smallsetminus\cD_0$ and $\widehat{\cW}_1\smallsetminus\cD_0$.

\noindent
\textit{Construction of $\widehat{\cU}_1$}. Recall from Construction \ref{const}, that a quadrilateral $Z=[0,1,z_3,z_4]\in\widehat{\cU}_1$ satisfies that $z_4 \in \Omega_{\cU_0}=\{z \in\mathbb{D}^2\colon\Im(z) \geq 0\}$ and when $z_4 \neq 1$, the set of possible values for $z_3$ is the curve $\gamma_{z_4}$ obtained by cutting the union $\overline{01}\cup\{ (1-z_4)t+1\colon t\geq 0\}$ with the constraint $|z_3-z_4| \leq 1$, Figure \ref{ball_boundary}.$a)$ (the case $z_4=1$ is pathological but consists of wedges, so it belongs to $\cD_0$). Since for $z_4\in \Omega_{\cU_0}\smallsetminus\{1\}$ the curves $\gamma_{z_4}$ are piecewise linear, then they can be parametrized by rescaling of the piecewise arc-length. So that, the function
\[(\Omega_{\cU_0} \smallsetminus \{1\}) \times [0,1] \rightarrow \widehat{\cU}_1, \quad (z,t) \mapsto [0,1,\gamma_{z}(t), z],\]
is a projection map onto its image which we denote by $X_{\widehat{\cU}} \subset \widehat{\cU}_1$. This projection defines a structure of an $I$-bundle, with collapsed fibers when both conditions, $|z_4|=1$ and $Re(z_4) \leq 0$ hold. Thus, $X_{\widehat{\cU}}$ is homeomorphic to a tetrahedron with an edge removed, corresponding to the fiber at $z_4=1$, see Figure \ref{ball_boundary}.$b)$.

\begin{figure}[h]
\begin{center}
\begin{tikzpicture}[scale=0.55]

\fill [gray!10] (0,0) arc (0:180:3) -- (0,0); 
\draw (-6,4.2) node {\large{$a)$}};
\draw (-4,1.5) node {\large{$\Omega_{\cU_0}$}};

\draw (0,0) arc (0:90:3);
\draw [thick,dashed] (-3,3) arc (90:180:3);
\draw (-6,0)--(0,0);

\draw (-3,3) node {{\tiny $\bullet$}};
\draw (-3,0) node {{\tiny $\bullet$}};
\draw (-6,0) node {{\tiny $\bullet$}};

\draw (-3,0) node [below] {$0$};
\draw (0,0) node [below] {$1$};
\draw (-6.1,0) node [below] {$-1$};
\draw (-3,3) node [above] {$i$};

\draw (-1.7,1.275) node {{\tiny $\bullet$}};
\draw (-1.7,1.2) node [above] {$z_4$};
\draw [thick] (-3,0)--(0,0);
\draw [dotted] (-1.7,1.275)--(0,0);
\draw [thick] (0,0)--(0.7,-0.525);
\draw (-1,0.1) node [below] {$\gamma_{z_4}$};

\draw (1.5,4.2) node {\large{$b)$}};
\draw (7.2,4) node {\large{$X_{\widehat{\cU}}$}};

\fill [gray!10] (5,0)--(10.85,0)--(7,2.5); 
\fill [gray!10] (8,0)--(10.75,2)--(10.75,0);

\draw (7.5,1) node {$\Omega_{\cU_0}$};
\draw [thick,dashed] (5,0)--(7,2.5);
\draw (7,2.5)--(10.9,4.1);
\draw (5,0)--(11,0);
\draw (7,2.5)--(7.7,2.0367);
\draw (7.95,1.872)--(9.25,1.01835);
\draw [black!40] (9.4,.9198)--(10.8,0);

\draw (5,0)--(10.9,4.1);
\draw (11,0)--(11,4);
\draw (10.75,0)--(10.75,4);
\draw [fill=white] (11,0) arc (0:360:0.12);
\draw [fill=white] (11,2) arc (0:360:0.12);
\draw [fill=white] (11,4) arc (0:360:0.12);

\draw (8,0) node [below] {$[0,1,0,0]$};
\draw (10.9,2) node [right] {$[0,1,1,1]$};
\draw (8,0)--(10.75,2);

\draw (5,0) node {\small{$\bullet$}};
\draw (7,2.5) node {\small{$\bullet$}};
\draw (8,0) node {\small{$\bullet$}};
\draw (5,0) node [left] {$[0,1,0,-1]$};
\draw (7,2.7) node [left] {$[0,1,0,i]$};
\draw (12,0) node [below] {$\cZ=[0,1,0,1]$};

\draw (3.6,1.4) node {collapsed fibers};
\draw (11.05,0.15) node [right] {removed point};

\end{tikzpicture}
\caption{$a)$ The half of disc $\Omega_{\cU_0}$ of possible values for $z_4$ in $\widehat{\cU}_1$ and an example of a fiber $\gamma_{z_4}$. For points $z_4$ in the dashed arc, the fibers $\gamma_{z_4}$ collapse to a point. $b)$ The topological tetrahedron with and edge removed, corresponding to the ``fiber'' at the zig-zag $\cZ$. The shaded regions represent the intersection $\cD_0\cap X_{\widehat{\cU}}$.}
\label{ball_boundary}
\end{center}
\end{figure}
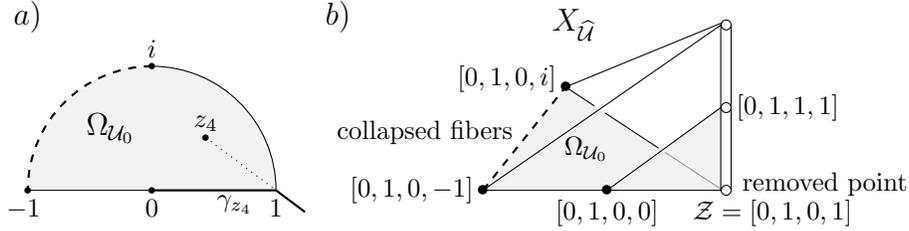

We may observe that the remaining wedges on the set $\widehat{\cU}_1$ are quadrilaterals of the form $[0,1,0,z_4]$, in the $I$-bundle structure they correspond to the bottom part of the fibers, \textit{i.e.} when $z_3=\gamma_{z_4}(0)$. Then removing such wedges is equivalent to remove the bottom face of the topological tetrahedron $X_{\widehat{\cU}}$. 

The quadrilaterals in the front face of the tetrahedron $X_{\widehat{\cU}}$ are $4$-segments of the form $[0,1,s,t]$, with $-1\leq t<1$, $0<s\leq2$ and $|s-t|\leq1$. In this face, the closure of the set of non-convex 4-segments is the triangle $\{[0,1,s,t]\colon 0\leq s\leq t\leq 1\}$, and it is adjacent to the bottom face of wedges previously mentioned.

We conclude that the set $\widehat{\cU}_1\smallsetminus\cD_0$ is equal to the tetrahedron $X_{\widehat{\cU}}$ with two adjacent triangles removed from its boundary, in particular, it is homeomorphic to a closed 3-ball with a disc removed from its boundary, Figure \ref{ball_boundary}.$b)$.

\medskip
\noindent
\textit{Construction of $\widehat{\cW}_1$}. As in the previous case, from Construction \ref{const}, a quadrilateral $Z = [0,1,z_3,z_4] \in \widehat{\cW}_1$ satisfies that $z_4 \in \Omega_{\cW_0} = \{z \in \mathbb{D}^2: \Im(z) \geq 0\}$ and when $z_4 \neq 1$, the set of possible values for $z_3$ is a curve $\gamma_{z_4}$ contained in the union $\overline{0z_4}\cup\{z_4+t(z_4-1)\colon t\geq 0\}$ and cutted with the constraint $|z_3-1|\leq 1$ (here, the case $z_4=1$ is again pathological and corresponds to wedges of the form $[0,1,z_3,1]$). Notice that the curves $\gamma_{z_4}$ are degenerated as point curves when $\Re(z_4) \leq 0$, corresponding to the wedges $[0,1,0,z_4]$. So we may take the region $\Omega'_{\cW_0} = \{z \in \Omega_{\cW_0} : \Re(z) \geq 0, \ z \neq 1\}$ and parameterize the curves $\gamma_{z_4}$ by a rescaling of the piecewise arc-length in order to have a projection map onto its image
    \[\Omega'_{\cW_0} \times [0,1] \rightarrow \widehat{\cW}_1, \quad (z,t) \mapsto [0,1,\gamma_{z}(t), z].\]
This projection defines a structure of an $I$-bundle, with collapsed fibers at the segment $\Re(z_4)=0$. As before, if we denote by $X_{\widehat{\cW}}$ the image of the projection map, then it is homeomorphic to a tetrahedron with an edge removed (the fiber that would correspond to $z_4=1$). We may observe that the remaining wedges on the set $\widehat{\cW}_1$, lying in the $I$-bundle structure are the quadrilaterals in the bottom part of the fibers, \textit{i.e.} when $z_3= \gamma_{z_4}(0)$ and having the form $[0,1,0,z_4]$, and thus, it is the bottom face of the tetrahedron $X_{\widehat{\cW}}$. The $4$-segments on the $I$-bundle structure are of the form $[0,1,s,t]$ with $0\leq s\leq t$, they are all non-convex and correspond to the front face of the tetrahedron $X_{\widehat{\cW}}$. In particular, $\widehat{\cW}_1 \smallsetminus\cD_0$ is homeomorphic to a tetrahedron with two adjacent faces removed from its boundary, and it is thus homeomorphic to a closed ball with a disc removed from its boundary.
\end{proof}

\begin{figure}[h]
\begin{multicols}{2}
\begin{center}
\begin{tikzpicture}[scale=0.6]
\draw [thick] (-3,-3)--(3,-3);
\draw [thick] (3,3)--(-3,3);
\draw [black!30] (-3,3)--(-3,-3);
\draw [black!30] (3,-3)--(3,3);

\draw (-3,-3) node [below] {$\widehat{\cV}_3$};
\draw (3,-3) node [below] {$\widehat{\cW}_4$};
\draw (3,3) node [above] {$\widehat{\cU}_4$};
\draw (-3,3) node [above] {$\widehat{\cW}_3$};

\draw [thick] (-1,1)--(-1,-1);
\draw [thick] (1,-1)--(1,1);
\draw [black!30] (1,1)--(-1,1);
\draw [black!30] (-1,-1)--(1,-1);
\draw (-0.95,-1.1) node [left] {$\widehat{\cV}_1$};
\draw (1,-1.1) node [right] {$\widehat{\cW}_1$};
\draw (1,1.1) node [right] {$\widehat{\cU}_2$};
\draw (-0.95,1.3) node [left] {$\widehat{\cW}_2$};

\draw [thick] (1,1)--(0,2);
\draw [black!30] (0,2)--(-1,1);
\draw (0,2) node [above] {$\widehat{\cV}_2$};

\draw [thick=0.4] (-1,1)--(-2,0)--(-1,-1);
\draw (-2,0) node [left] {$\widehat{\cU}_3$};

\draw [thick=0.4] (-1,-1)--(0,-2);
\draw [black!30] (1,-1)--(0,-2);
\draw (0,-2) node [below] {$\widehat{\cU}_1$};

\draw [thick=0.4] (1,-1)--(2,0)--(1,1);
\draw (2,0) node [right] {$\widehat{\cV}_4$};

\draw [thick] (3,3)--(0,2)--(-3,3);
\draw [thick] (-2,0)--(-3,-3);
\draw [black!30] (-3,3)--(-2,0);
\draw [thick] (-3,-3)--(0,-2)--(3,-3);
\draw [thick] (2,0)--(3,3);
\draw [black!30] (3,-3)--(2,0);

\draw (0,-2) node {$\bullet$};
\draw (-2,0) node {$\bullet$};
\draw (0,2) node {$\bullet$};
\draw (-1,-1) node {$\bullet$};
\draw (1,-1) node {$\bullet$};
\draw (1,1) node {$\bullet$};
\draw (-1,1) node {$\bullet$};
\draw (2,0) node {$\bullet$};
\draw (-3,-3) node {$\bullet$};
\draw (3,-3) node {$\bullet$};
\draw (3,3) node {$\bullet$};
\draw (-3,3) node {$\bullet$};
\end{tikzpicture}
\vspace{2cm}
\caption{Graph describing the combinatorics of the intersections of the sets $\widehat{\cU}_j\smallsetminus\cD_0$, $\widehat{\cV}_j\smallsetminus\cD_0$, $\widehat{\cW}_j\smallsetminus\cD_0$; and the bigon decompositions they induce.}
\label{pegado_frontera}
\end{center}
\end{multicols}
\end{figure}
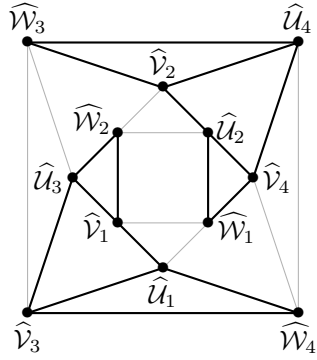

We want to understand the combinatorics of the intersections of the balls in Lemma \ref{lema:boundary_complement_of_wedges}. To do this, we first observe that it is possible to give to these balls the structure of a polyhedron, in the sense of Remark \ref{polyhedron_structure}, and codify the combinatorics of the intersections between neighboring balls into a graph, Figure \ref{pegado_frontera}. We describe this structure in the following Lemma.

\begin{lem}\label{lemma_pegados_frontera}
The boundaries of the balls $\widehat{\cU}_j\smallsetminus\cD_0$, $\widehat{\cV}_j\smallsetminus\cD_0$, and $\widehat{\cW}_j\smallsetminus\cD_0$ are splitted into bigons, corresponding to the intersection of pairs of such balls, but with some sides of the bigons removed. This combinatorial structure is captured in the graph of Figure \ref{pegado_frontera} as follows:
\begin{itemize}
    \item each vertex of the graph represents a ball,

    \item each black edge of the graph represents a bigon, such bigon is the intersection of the boundaries of the balls corresponding to the vertices of the graph,

    \item each grey edge of the graph represents a disc which is removed to the boundaries of both balls,

    \item each triangular or quadrialteral cycle of the graph whose edges are all black, represents the intersection of the corresponding balls and determine a side of the bigon in each ball,

    \item each triangular or quadrialteral cycle of the graph which contains a grey edge, represents a removed side of the bigon.
\end{itemize}
\end{lem}

\begin{ejem}
To illustrate the combinatorial description in Lemma \ref{lemma_pegados_frontera} we give two examples: Figure \ref{pegado_frontera} tells us that the boundary of $\widehat{\cU}_1\smallsetminus \cD_0$ splits into three bigons, corresponding with the intersections with $\widehat{\cV}_1$, $\widehat{\cW}_4$ and $\widehat{\cV}_3$, so that the first two bigons only have one side in $\widehat{\cU}_1\smallsetminus \cD_0$ and the third one have the two sides. Figure \ref{pegado_frontera} also tells us that the boundary of $\widehat{\cW}_1\smallsetminus \cD_0$ splits into two bigons, corresponding with the intersections with $\widehat{\cU}_2$, $\widehat{\cV}_4$ and both bigons only have one side in $\widehat{\cW}_1\smallsetminus\cD_0$. 
\end{ejem}

\begin{proof}
The proof is analogous to the proof of Lemma \ref{pegados}, in which we used the $I$-bundle structure to determine the bigonization of the boundary. It is enough to prove the result for the sets $\widehat{\cU}_1\smallsetminus\cD_0$ and $\widehat{\cW}_1\smallsetminus\cD_0$ because the action of the dihedral group implies the result for the remaining balls. In this case, we use the $I$-bundle structure provided by Lemma \ref{lema:boundary_complement_of_wedges}. The components of the boundary for $\widehat{\cU}_1 \smallsetminus \cD_0$ are:
\begin{enumerate}
    \item $\big(\widehat{\cU}_1\cap\widehat{\cV}_1\big) \smallsetminus \cD_0$ obtained as the $I$-bundle restricted to the segment $z_4\in[-1,1]$, cut with the condition $z_4 \geq z_3$, consisting of convex $4$-segments. There is a side removed in this bigon, corresponding with the condition $z_3=z_4$.
    
    \item $\big(\widehat{\cU}_1\cap \widehat{\cW}_4\big)\smallsetminus \cD_0$ obtained as the $I$-bundle restricted to the arc $|z_4|=1$, with an side removed, corresponding with the condition $z_1=z_3$.

    \item $\big(\widehat{\cU}_1\cap \widehat{\cV}_3\big)\smallsetminus \cD_0$ obtained as the top face of the $I$-bundle structure, that is when $|z_4 - z_3| = 1$. This bigon has the two vertices removed, corresponding to the conditions $z_4 = 1$ and $z_4 = -1$.
\end{enumerate}
The components of the boundary for $\widehat{\cW}_1 \smallsetminus \cD_0$ are:
\begin{enumerate}
    \item $\big(\widehat{\cW}_1\cap \widehat{\cU}_2\big) \smallsetminus \cD_0$ obtained as the top face of the $I$-bundle structure, that is when $|z_3 - z_2| = 1$. This bigon has one side removed, corresponding to the condition $z_4 \in [0,1]$, consisting of non-convex $4$-segments.

    \item $\big(\widehat{\cW}_1\cap \widehat{\cV}_4\big) \smallsetminus \cD_0$ obtained as the $I$-bundle restricted to the arc $|z_4| = 1$, with a side removed, corresponding with the condition $z_3 = z_1$.
\end{enumerate}
\end{proof}

Using this combinatorics we can conclude the following

\begin{cor}\label{boundary_balls}
For every $0 \leq s < 1$, the set $h^{-1}(s) \smallsetminus \cD_s$ is homeomorphic to $\R^3 \sqcup \R^3$.
\end{cor}

\begin{proof}
As observed in the proof of Proposition \ref{prop_limit_spheres}, for every $0 < s < 1$, we recall the decomposition of $h^{-1}(s)$ into the sets $\{\textrm{cl}(\cU_j),\textrm{cl}(\cV_j),\textrm{cl}(\cW_j)\}$ given in Lemma \ref{pegados} and combinatorically equivalent to the one given in Theorem \ref{levelspheres} and represented in Figure \ref{Diagrama_pegados_tridimensional}. From that description, we can see that the set $\mathcal{D}_s$ is the suspension over a closed curve in the regular cube shown in Figure \ref{Ds_separating_curve}, and so, $\cD_s$ is homeomorphic to $\mathbb{S}^2$. Moreover, such curve decomposes the cube into two topological discs, so that $h^{-1}(s) \smallsetminus \cD_s$ is a suspension over two disjoint discs, which is homeomorphic to $\R^3 \sqcup \R^3$. The case $s = 0$ follows from the fact that the description of the decomposition of $h^{-1}(0) \smallsetminus \cD_0$ into the sets $\{\widehat{\cU_j} \smallsetminus \cD_0, \widehat{\cV_j} \smallsetminus \cD_0, \widehat{\cW_j} \smallsetminus \cD_0 \}$ given by Lemma \ref{lemma_pegados_frontera} is combinatorically equivalent to the one given for every other level set $h^{-1}(s)$.
\end{proof}

In particular, the boundary $\partial S(4)$ can be seen as $\partial S(4) \smallsetminus \cD_0 \cong \R^3 \sqcup \R^3$, compactified with a wedge of spheres and a disc $\cD_0 \cong (\S^2 \vee \S^2) \cup \D^2$, and therefore, $\partial S(4)$ is not homeomorphic to the sphere $\S^3$. 

\subsection{Good asymptotic behaviour at the boundary}\label{subsection_good_asymptotic}

Corollary \ref{boundary_balls} suggests that we can extend the cone structure $S(4) \cong C(\S^3)$ to the boundary outside $\cD$. In this section we prove this, and the most difficult part is to give to the sets $\{\overline{\cU_j} \smallsetminus \cD, \overline{\cV_j} \smallsetminus \cD,\overline{\cW_j} \smallsetminus \cD\}$, the structure of a cylinder. We start by looking at an example of another type of bad convergence and explain why this complicates the argument.

\begin{ejem}\label{bad_convergence_ex}
For every $4$-segment $Z = [0,1,x,y] \in \widehat{\cU}_1$ with $0 < x,y < 1$ and for every $0 \leq r \leq y$, there exist a sequence of simple quadrilaterals $Z_k = [0,1,z_k, w_k] \in \textrm{cl}(\cU_1)$ such that $h(Z_k) = s_k$, $Z_k \rightarrow Z$ as $k \rightarrow \infty$, and the fibers of the $I$-bundle in $\textrm{cl}(\cU_1) \cap h^{-1}(s_k)$ at $Z_k$, denoted as $\gamma_{z_4(k)}(t)$, satisfy that $\gamma_{z_4(k)}(t)$ converge to the interval $[r, y+1]$. \textbf{Thus, different sequences converging to $Z$, lie on fibers which converge to different curves at the boundary}. The sequence can be constructed as follows: for every $0 < s < 1$, the complex number $y + i s = \eta e^{i\theta}$ has an element in the bisector of the angle $\theta$ satisfying $\nu e^{i\theta/2} = r + i s'$ for some $s' \in (0,s)$, then the sequence can be constructed with the quadrilaterals $Z_{s'} = [0,1,x + is',y + i s] \in \textrm{cl}(\cU_1) \cap h^{-1}(s')$, as $Z_{s'} \rightarrow Z$ when $s \rightarrow 0$ and the fiber of the $I$-bundle structure at the point $z_4=y + is$ starts with the segment $\{ t + i s' : r \leq t \leq y \}$.
\end{ejem}

\begin{obs}\label{remark_general_bundle_cylinder_extension}
Recall from the proof of Lemma \ref{eucliconos}, where we constructed homeomorphisms $f_s : \textrm{cl}(\cU_1) \cap h^{-1}(s) \rightarrow \D^3$,
induced by the $I$-bundle structure, which vary continuously on the parameter $s \in (0,1)$ and gives us the cylinder structure 
    \[  \textrm{cl}(\cU_1) \smallsetminus \{\mathfrak{R}_4\} \rightarrow \D^3 \times (0,1), \quad Z \mapsto (f_{h(Z)}(Z), h(Z)). \]
If we restrict these homeomorphisms, we get homeomorphisms
$f_s : (\textrm{cl}(\cU_1) \smallsetminus \cD) \cap h^{-1}(s) \rightarrow \H^2 \times (0,1] \cong \H^3$, then by Lemma \ref{pegados} and Lemma \ref{lemma_pegados_frontera}, these homeomorphisms
are combinatorically equivalent to the homeomorphism at the boundary $f_0 : \widehat{\cU}_1 \smallsetminus \cD_0 \rightarrow \H^2 \times (0,1]$. For the particular case of $4$-segments we can write $f_0([0,1,x,y]) = \left(y', \frac{x}{y+1}\right)$ (for some $y' \in \H^2$) and Example \ref{bad_convergence_ex} tells us that for every $r \in [0,y]$, there is a sequence $\{Z_k\}$, such that $Z_k \rightarrow [0,1,x,y]$ and $\lim_{k \rightarrow \infty} f_{s_k}(Z_k) = \left(y', \frac{x-r}{y+1-r}\right)$, so the extended family of functions $\{f_s\}_{s \in [0,1)}$ is no longer continuous when $s \rightarrow 0$. In fact, there is no reasonable way to define $f_0$ using the $I$-bundle structure to obtain such continuity. The problem is of course the lack of uniqueness of convergence of fibers of the $I$-bundle structures when converging to $4$-segments. But, if we restrict the $I$-bundle structure on each level set $h^{-1}(s) \cap \overline{\cU_1}$, so that we have uniqueness on the convergence of fibers, then we have a cylinder structure extended to the boundary. Analogous behaviour occur for the $I$-bundle structures on the rest of the sets $\{h^{-1}(s) \cap \overline{\cU_j}, h^{-1}(s) \cap \overline{\cV_j}, h^{-1}(s) \cap \overline{\cW_j} : 1 \leq j \leq 4\}$.
\end{obs}

We overcome the difficulty illustrated by Remark \ref{remark_general_bundle_cylinder_extension} and Example \ref{bad_convergence_ex} by proving regularity on different regions and gluing the pieces with Theorem \ref{thm_gluing_upper_spaces}

\begin{lem}\label{conos_homeo_h4}
The sets $\overline{\cU_j}\smallsetminus\cD$, $\overline{\cV_j}\smallsetminus\cD$, and $\overline{\cW_j}\smallsetminus\cD$, are homeomorphic to $\H^4$.
\end{lem}

\begin{proof}
As it was done in multiple times already, thanks to the action of the dihedral group $\langle\sigma,\tau\rangle$, it is enough to prove the result for $\overline{\cU_1}\smallsetminus\cD$ and $\overline{\cW_1}\smallsetminus\cD$. 

\noindent
Unlike previous proofs, we start here with the set $\overline{\cW_1} \smallsetminus \cD$, because in this case, it is easier. Recall from the construction on Lemma \ref{eucliconos}, that we can uniformize the $I$-bundle structure of the cone $\textrm{cl}(\cW_1)$ on each level set $h^{-1}(s)$, and construct the function
    \[   \Omega \times [0,1] \rightarrow \textrm{cl}(\cW_1) \cap h^{-1}(s), \qquad  (w,t) \mapsto [0,1,\gamma_{\varphi_s(w)}(t),  \varphi_s(w)] \]
(the curves $\gamma_z$ are the fibers of the $I$-bundle on $\textrm{cl}(\cW_1) \cap h^{-1}(s)$ and $\varphi_s$ is an homeomorphism between $\Omega = \{z \in \D^2 : \Im(z) \geq 0\}$ and the base space of the bundle, cf. Figure \ref{Omega_W}). Taking into account the collapsed fibers at the boundary, the previous function induces the homeomorphism 
    \[    f_s : \textrm{cl}(\cW_1) \cap h^{-1}(s) \rightarrow \frac{\left(\Omega \times [0,1]\right) }{\sim} \cong \D^3,  \]
which varies continuously on the parameter $s \in (0,1)$. We can restrict this homeomorphism to $\Omega'= \{ z \in \Omega : \Re(z) \leq 0, \ \Im(z) > 0 \} \cong \H^2$, and observe that
    \[ \big(\Omega' \times (0,1] \big) \times (0,1) \rightarrow \overline{\cW_1} \smallsetminus \cD, \qquad  (w,t,s) \mapsto [0,1,\gamma_{\varphi_s(w)}(t),  \varphi_s(w)],
    \]
defines an homeomorphism onto its image, giving us the homeomorphism $X = \{[0,1,z,w] \in \overline{\cW_1} \smallsetminus \cD : w \in \Omega'\} \cong \H^4$ (this is immediate from Lemma \ref{eucliconos}, since $X \cap \partial S(4) = \emptyset$, see the construction of $\widehat{\cW}_1$ in proof of Lemma \ref{lema:boundary_complement_of_wedges}). If we now restrict the homeomorphisms $f_s$ to the subset $\Omega'' = \{ z \in \Omega : \Re(z) \geq 0,  \Im(z) > 0\} \cong \H^2$, then $f_s$ restricts to a homeomorphism onto $\frac{\left(\Omega'' \times (0,1]\right) }{\sim} \cong \H^3$ and this is combinatorically equivalent to the homeomorphism 
    \[ f_0 : \widehat{\cW}_1 \smallsetminus \cD_0 \rightarrow \frac{\left(\Omega'' \times (0,1]\right) }{\sim} \cong \H^3,  \]
given in Lemma \ref{lema:boundary_complement_of_wedges}, using the corresponding $I$-bundle structure. Moreover, we can see that for every $z_4 \in \Omega''$, the quadrilaterals $[0,1,z_3,z_4] \in \widehat{\cW}_1$ have the property of unique convergence of fibers as discussed in Remark \ref{remark_general_bundle_cylinder_extension}, which implies that the function 
        \[  \big(\Omega'' \times (0,1] \big) \times [0,1) \rightarrow \overline{\cW}_1 \smallsetminus \cD, \qquad  (w,t,s) \mapsto [0,1,\gamma_{\varphi_s(w)}(t),  \varphi_s(w)] \] 
induces the homeomorphism $Y = \{[0,1,z,w] \in \overline{\cW_1} \smallsetminus \cD : w \in \Omega''\} \cong \H^4$. We can observe that $X$ and $Y$ intersect at their boundary in
    \[ X \cap Y = \{[0,1,z,w] \in \overline{\cW_1} \smallsetminus \cD : \Re(w) = 0\} \cong \H^3,    \]
moreover, $\partial X \smallsetminus X \cap Y \cong \H^3$ and $\partial Y \smallsetminus X \cap Y \cong \H^3$. As $\overline{\cW_1} \smallsetminus \cD = X \cup Y$ can be obtained as two spaces homeomorphic to $\H^4$ glued together by a well placed copy of $\H^3$ at their boundaries, then by Theorem \ref{thm_gluing_upper_spaces}, we conclude $\overline{\cW_1} \smallsetminus \cD \cong \H^4$.

\noindent
\textit{Proof for the set $\overline{\cU_1} \smallsetminus \cD$.} We proceed as in the previous case by taking the function
    \[ F : \big(\Omega \times (0,1] \big) \times (0,1) \rightarrow \overline{\cU_1} \smallsetminus \cD, \qquad  (w,t,s) \mapsto [0,1,\gamma_{\varphi_s(w)}(t),  \varphi_s(w)],
    \]
with $\gamma_z$ the fibers of the $I$-bundle on $\textrm{cl}(\cU_1) \cap h^{-1}(s)$ and $\varphi_s$ an homeomorphism between $\Omega = \{z \in \D^2 : \Im(z) \geq 0\}$ and the base space of the bundle, cf. Figure \ref{Omega_U}. The function $F$ induces the homeomorphism $\textrm{cl}(\cU_1) \smallsetminus \cD \cong \H^4$, but as illustrated in Example \ref{bad_convergence_ex} and Remark \ref{remark_general_bundle_cylinder_extension}, such function is not continuous when approaching to the boundary at $4$-segments. We overcome this difficulty by restricting the function to specific regions were we do have continuity (which is everything we need to obtain cylinder structures) and proving that $\overline{\cU_1} \smallsetminus \cD$ is a manifold with boundary. 

Just as it happened in the previous case, when the fourth vertex $z_4$ is bounded away from the real line, we have uniqueness of convergence of fibers, so the function $F$ induces a cylinder structure when we impose the condition $\Im(z_4) > 0$ and thus, we have $\overline{\cU_1} \smallsetminus (\overline{\cV_1}\cup \cD) \cong \H^4$ (this set is obtained by extracting the front face in each level set $\overline{\cU_1}\cap h^{-1}(s)$, such faces correspond to trapezoid quadrilaterals for $\overline{\cU_1}\cap h^{-1}(s)$ with $0<s<1$, and to 4-segments in the boundary $\widehat{\cU}_1\cap h^{-1}(0)$). We can also observe that if we restrict the convergence  to $4$-segments with only trapezoid quadrilaterals, the phenomenon of bad convergence in Example \ref{bad_convergence_ex} don't appear, so that we have again uniqueness of convergence of fibers and thus, the function $F$ induces the homeomorphism $(\overline{\cU_1} \cap \overline{\cV_1}) \smallsetminus \cD \cong \H^3$. We have thus the decomposition
    \[  \big(\overline{\cU_1} \smallsetminus \cD \big) = \big( \overline{\cU_1} \smallsetminus (\overline{\cV_1} \cup \cD) \big) \sqcup \big( (\overline{\cU_1} \cap \overline{\cV_1} ) \smallsetminus \cD \big)\cong \H^4 \sqcup \H^3,  \]
where the inclusion is $(\overline{\cU_1} \cap \overline{\cV_1} ) \smallsetminus \cD \hookrightarrow \partial \big(\overline{\cU_1} \smallsetminus \cD\big)$. From Lemmas \ref{nivelconos} and \ref{lema:boundary_complement_of_wedges}, it follows that
    \[   \partial (\overline{\cU_1} \smallsetminus \cD) = (\widehat{\cU}_1\smallsetminus \cD_0) \cup \partial (\textrm{cl}(\cU_1) \smallsetminus \cD) \cong \H^3 \cup \big(\H^2 \times (0,1) \big) \cong \H^3 \cup \H^3, \]
and as the sets $(\widehat{\cU}_1 \smallsetminus \cD_0)$ and $\partial (\textrm{cl}(\cU_1) \smallsetminus \cD)$ are glued on $\partial (\widehat{\cU}_1 \smallsetminus \cD_0) \cong \H^2$, then we have that $\partial (\overline{\cU_1} \smallsetminus \cD) \cong \H^3$ by Theorem \ref{thm_gluing_upper_spaces}. 

Moreover, we can proof that $\overline{\cU_1} \smallsetminus \cD$ is manifold with boundary (by all the previous discussion, it only remains to prove this fact around convex $4$-segments): take a $4$-segment $[0,1,x,y] \in \overline{\cU_1}  \smallsetminus \cD$ (we have the conditions $0 < x < 2$, $-1 < y < 1$ and $y < x$), then for $0 < \varepsilon < |x-y| / 2$ small enough, we consider the neighborhood
    \[ A = \{[0,1,z_3,z_4] \in \overline{\cU_1} \smallsetminus \cD : z_3\in S_{\varepsilon}(x), \ z_4 \in S_{\varepsilon}(y) \},  \]
where $S_{\varepsilon}(z) = \{  w \in \C : |\Re(z - w)| < \varepsilon, |\Im(z - w)| < \varepsilon \}$ is a square neighborhood around the complex number $z$. We distinguish to cases:
\begin{enumerate}
    \item When $0 \leq x+\varepsilon \leq 1$, in which case, the neighborhood $A$ is just the two $\varepsilon$-squares $S_{\varepsilon}(x)$ and $S_{\varepsilon}(y)$ with the extra conditions $0 \leq \Im(z_4)$ and $0 \leq \Im(z_3)\leq \Im(z_4)$, see Figure \ref{neighborhood_segment}.$a)$. More precisely, we have a map
    \[  \{ (a,b) \in \R^2 : |a|, |b| < \varepsilon \} \times \{ (c,d) \in \R^2 : 0 \leq c \leq d < \varepsilon \} \rightarrow A, \]
given by $\big((a,b), (c,d)\big) \mapsto [0,1,x + a + i c, y + b + i d]$, which is a homeomorphism, so that $A \cong \R^2 \times \H^2 \cong \H^4$.

\item The general case, where the neighborhood $A$ is just the two $\varepsilon$-squares $S_{\varepsilon}(x)$ and $S_{\varepsilon}(y)$ with the extra conditions $0 \leq \Im(z_4)$, $\Im(z_3) \leq \Im(z_4)$ and $z_3$ is above a piecewise linear curve (above $\{tz_4 + (1-t) : t \in \R\}$ when $\Re(z_3) \geq 1$ and $\Im(z_3) \geq 0$ when $\Re(z_3) \leq 1$), see Figure \ref{neighborhood_segment}.$b)$. This is topologically equivalent to the previous case, so we have once again $A \cong \H^4$.
\end{enumerate}

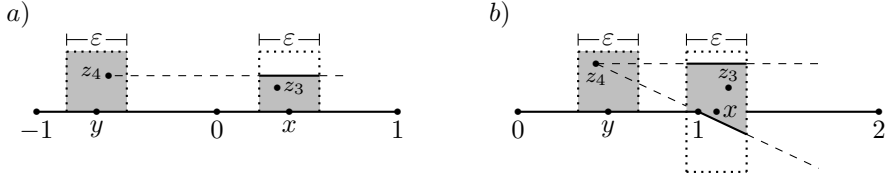
\begin{figure}[!ht]
\begin{center}
\begin{tikzpicture}[scale=.8]
\draw (-7.3,2) node [below] {$a)$};

\fill [gray!50] (-3.3,0)--(-3.3,0.6)--(-2.3,0.6)--(-2.3,0);
\fill [gray!50] (-6.5,0)--(-6.5,1)--(-5.5,1)--(-5.5,0);

\draw [thick] (-7,0)--(-1,0);

\draw (-7,0) node {\tiny{$\bullet$}};
\draw (-7,0) node [below] {$-1$};
\draw (-4,0) node {\tiny{$\bullet$}};
\draw (-4,0) node [below] {$0$};
\draw (-1,0) node {\tiny{$\bullet$}};
\draw (-1,0) node [below] {$1$};

\draw (-6,0) node {\tiny{$\bullet$}};
\draw (-6,0) node [below] {$y$};
\draw [thick,dotted] (-6.5,0)--(-6.5,1);
\draw [thick,dotted] (-6.5,1)--(-5.5,1);
\draw [thick,dotted] (-5.5,0)--(-5.5,1);
\draw (-6.5,1.1)--(-6.5,1.3);
\draw (-6.5,1.2)--(-6.15,1.2);
\draw (-5.5,1.1)--(-5.5,1.3);
\draw (-5.85,1.2)--(-5.5,1.2);
\draw (-6,1.2) node {$\varepsilon$};

\draw (-5.8,0.6) node {\tiny{$\bullet$}};
\draw (-5.75,0.65) node [left] {\footnotesize{$z_4$}};
\draw [dashed] (-5.7,0.6)--(-3.3,0.6);
\draw [thick] (-3.3,0.6)--(-2.3,0.6);
\draw [dashed] (-2.3,0.6)--(-1.8,0.6);
\draw (-3,0.4) node {\tiny{$\bullet$}};
\draw (-3.05,0.4) node [right] {\footnotesize{$z_3$}};

\draw (-2.8,0) node {\tiny{$\bullet$}};
\draw (-2.8,0) node [below] {$x$};
\draw [thick,dotted] (-3.3,0)--(-3.3,1);
\draw [thick,dotted] (-3.3,1)--(-2.3,1);
\draw [thick,dotted] (-2.3,0)--(-2.3,1);
\draw (-3.3,1.1)--(-3.3,1.3);
\draw (-3.3,1.2)--(-2.95,1.2);
\draw (-2.3,1.1)--(-2.3,1.3);
\draw (-2.65,1.2)--(-2.3,1.2);
\draw (-2.8,1.2) node {$\varepsilon$};

\draw (0.7,2) node [below] {$b)$};

\fill [gray!50] (2,0)--(2,1)--(3,1)--(3,0);
\fill [gray!50] (3.8,0.8)--(4.8,0.8)--(4.8,-0.3764)--(4,0)--(3.8,0);

\draw [thick] (1,0)--(4,0);
\draw [thick] (4.8,0)--(7,0);

\draw (1,0) node {\tiny{$\bullet$}};
\draw (1,0) node [below] {$0$};
\draw (2.5,0) node {\tiny{$\bullet$}};
\draw (2.5,0) node [below] {$y$};
\draw (4,0) node {\tiny{$\bullet$}};
\draw (4,0) node [below] {$1$};
\draw (4.3,0) node {\tiny{$\bullet$}};
\draw (4.25,0) node [right] {$x$};
\draw (7,0) node {\tiny{$\bullet$}};
\draw (7,0) node [below] {$2$};

\draw [thick,dotted] (2,0)--(2,1);
\draw [thick,dotted] (2,1)--(3,1);
\draw [thick,dotted] (3,0)--(3,1);
\draw (2,1.1)--(2,1.3);
\draw (2,1.2)--(2.35,1.2);
\draw (3,1.1)--(3,1.3);
\draw (2.65,1.2)--(3,1.2);
\draw (2.5,1.2) node {$\varepsilon$};

\draw (2.3,0.8) node {\tiny{$\bullet$}};
\draw (2.3,0.8) node [below] {\footnotesize{$z_4$}};
\draw [dashed] (2.3,0.8)--(3.8,0.094);
\draw [thick] (4,0)--(4.8,-0.3764);
\draw [dashed] (4.8,-0.3764)--(6,-0.9411);
\draw [dashed] (2.3,0.8)--(3.8,0.8);
\draw [thick] (3.8,0.8)--(4.8,0.8);
\draw [dashed] (4.8,0.8)--(6,0.8);
\draw (4.5,0.4) node {\tiny{$\bullet$}};
\draw (4.5,0.36) node [above] {\footnotesize{$z_3$}};

\draw [thick,dotted] (3.8,-1)--(3.8,1);
\draw [thick,dotted] (3.8,1)--(4.8,1);
\draw [thick,dotted] (4.8,-1)--(4.8,1);
\draw [thick,dotted] (3.8,-1)--(4.8,-1);
\draw (3.8,1.1)--(3.8,1.3);
\draw (3.8,1.2)--(4.15,1.2);
\draw (4.8,1.1)--(4.8,1.3);
\draw (4.45,1.2)--(4.8,1.2);
\draw (4.3,1.2) node {$\varepsilon$};
\end{tikzpicture}
\caption{Square neighborhoods around non-convex 4-segments in $\overline{\cU_1} \smallsetminus \cD$.}
\label{neighborhood_segment}
\end{center}
\end{figure}

By Brown's theorem \cite{brown2}, there is a collar neighborhood around the set of convex $4$-segments $\H^4 \cong (\widehat{\cU}_1 \cap\widehat{\cV}_1 \smallsetminus \cD_0) \times [0,1]^2 \hookrightarrow \overline{\cU_1} \smallsetminus \cD$. This collar neighborhood can be described as follows: there is a continuous function $f : [-1,1] \rightarrow \D^2$ such that $f(-1) = f(1) = 0$, $0 < f(t) < 1$ for every $t \in (-1,1)$, so that the collar neighborhood is thus just the set
    \[   \mathcal{C} = \{ [0,1,z_3,z_4] \in \overline{\cU_1} \smallsetminus \cD : \Im(z_4) \leq f(\Re(z_4))  \} \cong \H^4.   \]
This neighborhood has three boundary components defined by the identities $\Im(z_3) = \Im(z_4)$ (corresponding to $\overline{\cU_1} \cap \overline{\cV_1}$),  $\Im(z_3) = 0$ (corresponding to $\widehat{\cU}_1$), and $\Im(z_4) = f(t)$. The complement of the collar neighborhood can be uniformized with $F$ giving us
    \[  \overline{\cU_1} \smallsetminus (\cD \cup \mathcal{C})= \{ [0,1,z_3,z_4] \in \overline{\cU_1} \smallsetminus \cD : \Im(z_4) \geq f(\Re(z_4))  \} \cong \H^4.  \]
The collar neighborhood and its complement inside $\overline{\cU_1} \smallsetminus \cD$ intersect at
    \[   \{ [0,1,z_3,t + i f(t)] \in \overline{\cU_1} \smallsetminus \cD : -1 < t < 1  \} \cong \H^3,   \]
and thus $\overline{\cU_1} \smallsetminus \cD \cong \H^4$ by Theorem \ref{thm_gluing_upper_spaces}.
\end{proof}

Now we are able to give a description of the closure $\overline{S(4)}$ and of the boundary $\partial S(4)$ as an asymptotic decomposition.

\begin{thm}\label{thm:boundary_decomposition}
The set $\overline{S(4)} \smallsetminus\cD$ has two connected components, both homeomorphic to $\H^4$.
\end{thm}

\begin{proof}
We use Lemma \ref{conos_homeo_h4} and Theorem \ref{thm_gluing_upper_spaces} to glue pairwise copies of $\H^4$. The proof follows by comparing the combinatorics of each level as done in Corollary \ref{boundary_balls}.
\end{proof}

The combination of Theorem \ref{thm:boundary_decomposition} and Proposition \ref{prop_limit_spheres}, gives the conclusion of Theorem \ref{thm_principal_two}.

\subsection{Compactification of the moduli space $\cM(4)$}\label{subsection_compactification_moduli}

We describe the quotient $\partial\cM(4):=\partial S(4)/\langle\sigma,\tau\rangle$ and some properties of the compactification $\overline{\cM(4)}$.

\begin{lem}\label{lema_cociente_disco}
If we denote by $\cD' = \cD /\langle\sigma,\tau\rangle$, then $\cD' \cap \cM(4) \cong C(\D^2)$. Moreover, the boundary $\cD' \cap \partial \cM(4)$ is also homeomorphic to $\D^2$.
\end{lem}

\begin{proof}
For every $0 < s < 1$, the set $h^{-1}(s) \cap \mathrm{cl}(\cU_1) \cap \mathrm{cl}(\cW_1)$ is a fundamental domain for the action of $\langle \sigma, \tau \rangle$ on $\cD_s$, which is a topological bigon whose boundary is decompose as a union of $h^{-1}(s) \cap \mathrm{cl}(\cU_1) \cap \mathrm{cl}(\cV_1) \cap \mathrm{cl}(\cW_1)$ with $h^{-1}(s) \cap \mathrm{cl}(\cU_1) \cap \mathrm{cl}(\cW_1) \cap \mathrm{cl}(\cV_4) \cap \mathrm{cl}(\cW_4)$ (see Lemma \ref{pegados}). As the only non-trivial identification on the boundary is given by $\tau$, acting as a folding on the curve $h^{-1}(s) \cap \mathrm{cl}(\cU_1) \cap \mathrm{cl}(\cV_1) \cap \mathrm{cl}(\cW_1)$, the resulting quotient is homeomorphic to $\D^2$ (c.f. proof of Lemma \ref{Lema:level_sphere_moduli}). As the fundamental domain is constructed in $\textrm{cl}(\cU_1)$, then the cone structure on $\cD' \cap \cM(4)$ follows from Lemma \ref{eucliconos}. For the boundary part $\cD' \cap \partial \cM(4)$, we proceed by analyzing separately the actions of $\langle\sigma,\tau\rangle$ in the sets of non-convex 4-segments $\cN$ and wedges $\cK$.

\noindent
\textit{Non-convex 4-segments.} Since $\cN=\bigcup_{j=0}^3\sigma^j T$ with $T=\{ [0,1,t,x]\colon 0 \leq t\leq x\leq 1 \}$ (see Proposition \ref{prop_topology_wedges_segments} and Figure \ref{wedge_square}), the triangle $T$ is a fundamental domain for action of the subgroup $\langle\sigma\rangle$. The element $\tau$ acts as the reflection of $T$ through the line $\{[0,1,1-t,t]\colon \frac{1}{2} \leq t \leq 1\}$. Then the quotient $\cN/\langle\sigma,\tau\rangle$ is identified with the triangle 
    \[  \{ [0,1,t,x]\colon 0 \leq t\leq x\leq 1-t \},    \]
whose sides are
\[\{[0,1,0,t]\colon 0\leq t\leq 1\}\cup\left\{[0,1,t,t]\colon 0 \leq t \leq\frac{1}{2}\right\}\cup\left\{[0,1,t,1-t]\colon \frac{1}{2} \leq t \leq 1 \right\}.
\]

\noindent
\textit{Wedges.} The action of $\langle \sigma \rangle$ in an element $Z=[0,1,0,z_4]\in\cK$ is as follows:
\[\sigma(Z)=[0,1,1-z_4,1],~~\sigma^2(Z)=\left[0,1,0,\frac{1}{z_4}\right],~~\sigma^3(Z)=\left[0,1,\frac{z_4-1}{z_4},1\right].\]
Notice that $\sigma$ and $\sigma^3$ switch the spheres in $\cK\cong\S^2\vee\S^2$, and $\sigma^2$ acts as the function $z\mapsto\frac{1}{z}$ in the sphere $\{[0,1,0,z_4]\colon z_4\in\C\P^1\}$ (Proposition \ref{prop_topology_wedges_segments}). Moreover, the map $\tau\sigma[0,1,0,z_4]=[0,1,0,\overline{z_4}]$ is the complex conjugate in the same sphere. Then, we can identify the quotient $\cK/\langle\sigma,\tau\rangle$ with the half disc $\{[0,1,0,z]\colon z\in\D,~ Im(z)\geq 0\}$. We conclude that the quotients $\cN/\langle\sigma,\tau\rangle$ and $\cK/\langle\sigma,\tau\rangle$ are identified through the segment $\{[0,1,0,t]\colon 0\leq t\leq 1\}$ at their boundaries, resulting in a space homeomorphic to $\D^2$, see the bottom in Figure \ref{cone_discs}$.b)$.
\end{proof}

\begin{obs}\label{remark_compactification_cone_over_disc}
The set $\cD'$, as a compactification of the open cone $C(\D^2)$ with a disc $\D^2$ at its boundary is somewhat involved. To illustrate the asymptotic behaviour to the boundary, we make the following identification: if we take $\Omega = \{ (z,s) \in \C \times \R : |z| \leq 1, \ \Im(z) \geq s > 0 \}$, then the projection 
\[\textrm{cl}(\cW_1) \cap \textrm{cl}(\cU_1) \rightarrow \Omega \times (0,1], \qquad Z = [0,1,w,z] \mapsto (z,h(Z))   \]
induces the homeomorphism $\big(\cD' \cap \cM(4) \big) \cong \Omega/{\sim} \cong C(\D^2)$ stated in Lemma \ref{lema_cociente_disco}, where the equivalence relation $(z,s) \sim (z',s)$, when $\Im(z) = s$, is the one induced by $\tau$ in $\overline{\cU_1} \cap \overline{\cV_1} \cap \overline{\cW_1}$ (which is analogous to a reflection), see Figure \ref{cone_discs}. Under this identification, a sequence $\{[z_n,s_n]\}_{n \in \N} \subset \Omega/{\sim}$, satisfying $s_n \rightarrow 0$ and $z_n \rightarrow z \in \{w \in \C : |w| \leq 1, \ \Im(w) \geq 0\}$, converges to a quadrilateral in the quotient $\cD_0/\langle\sigma,\tau\rangle$, according to the following three possibilities:
\begin{enumerate}
    \item If $\Im(z) > 0$ or $Im(z)=0$ and $z\leq 0$, then the sequence converges uniquely to the wedge $[0,1,0,z]$.

    \item If $z \in \R$ and $z >0$, then the sequence converges to a $4$-segment $[0,1,x,z]$ with $0\leq x\leq z$. We can obtain every non-convex $4$-segment like this by controlling the slope $s_n /\Im(z_n)$, as we have the relation
    \[  \lim_{n \rightarrow \infty} \frac{s_n}{\Im(z_n)} = \frac{x}{2z}.  \]
    Such sequences were constructed in the proof of Proposition \ref{prop_limit_spheres}.

    \item If $z = 1$, then the sequence converges to a wedge $[0,1,x,1]$ (a convex 4-segment), with $1 \leq x \leq 2$, which satisfies $\tau[0,1,x,1]=[0,1,0,1-x]$. We can obtain every wedge satisfying this condition using parallelograms $[0,1,x + is_n, z_n = a_n + i s_n] \in h^{-1}(s_n) \cap \overline{\cU_1} \cap \overline{\cV_1} \cap \overline{\cW_1}$. 
\end{enumerate}

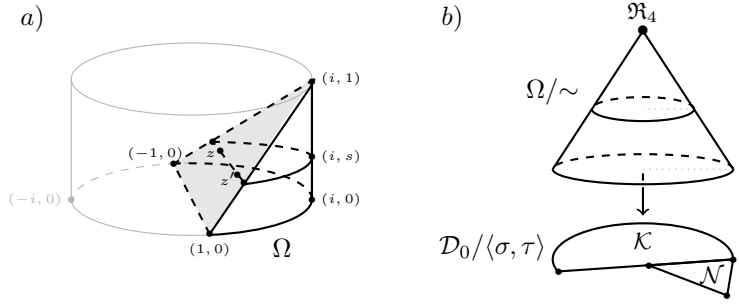
\begin{figure}[h]
\begin{center}
\begin{tikzpicture}[scale=0.8]
\draw (-7,3) node [right] {$a)$};
\fill [gray!20] (-2,2)--(-3.7,-0.58)--(-4.3,0.59);

\draw [black!30,dashed] (-4,0.6) arc (91:180:2 and 0.6);
\draw [black!30] (-6,0) arc (180:278:2 and 0.6);
\draw [black!30] (-6,0)--(-6,2);
\draw [black!30] (-2,2) arc (0:360:2 and 0.6);

\draw [thick,dashed] (-2,0) arc (0:100:2 and 0.6);
\draw [thick] (-2,0) arc (0:-82:2 and 0.6);
\draw [thick] (-2,0)--(-2,2);

\draw [black] (-3.7,-0.58) node {\tiny{$\bullet$}};
\draw (-3.7,-0.58) node [below] {\tiny{$(1,0)$}};

\draw [black!30] (-6,0) node {\tiny{$\bullet$}};
\draw [black!30] (-6,0) node [left] {\tiny{$(-i,0)$}};

\draw [black] (-4.3,0.59) node {\tiny{$\bullet$}};
\draw (-4.6,0.5) node [above] {\tiny{$(-1,0)$}};

\draw [black] (-2,0) node {\tiny{$\bullet$}};
\draw (-2,0) node [right] {\tiny{$(i,0)$}};

\draw (-2,2) node [right] {\tiny{$(i,1)$}};

\draw [black] (-2,1.95) node {\tiny{$\bullet$}};
\draw [thick] (-2,0)--(-2,2);
\draw [thick,dashed] (-3.7,-0.58)--(-4.3,0.59);
\draw [thick] (-2,1.9)--(-3.7,-0.58);
\draw [thick,dashed] (-1.95,2)--(-4.3,0.59);

\draw [thick,dashed] (-2,0.7) arc (0:80:2 and 0.26);
\draw [thick] (-2,0.7) arc (0:-64.5:2 and 0.5);
\draw [black] (-2,0.7) node {\tiny{$\bullet$}};
\draw (-2,0.7) node [right] {\tiny{$(i,s)$}};

\draw [black] (-3.13,0.28) node {\tiny{$\bullet$}};
\draw [black] (-3.65,0.96) node {\tiny{$\bullet$}};
\draw [thick,dashed] (-3.13,0.28)--(-3.65,0.96);

\draw [black] (-3.52,0.8) node {\tiny{$\bullet$}};
\draw (-3.45,0.75) node [left] {\tiny{$z$}};
\draw [black] (-3.24,0.41) node {\tiny{$\bullet$}};
\draw (-3.1,0.35) node [left] {\tiny{$z'$}};

\draw (-2.5,-0.5) node [below] {$\Omega$};

\draw (0,3) node [right] {$b)$};
\draw [black!40,dotted] (5,0.5)--(3.5,0.5);
\draw [black!40,dotted] (4.35,1.5)--(3.5,1.5);
\draw [thick,dashed] (5,0.5) arc (0:180:1.5 and 0.25);
\draw [thick] (2,0.5) arc (180:360:1.5 and 0.25);
\draw [thick,dashed] (4.35,1.5) arc (0:180:0.85 and 0.2);
\draw [thick] (4.35,1.5) arc (0:-180:0.85 and 0.2);
\draw [thick] (2,0.5)--(3.5,2.8);
\draw [thick] (5,0.5)--(3.5,2.8);

\draw [thick] (3.5,0.43)--(3.5,0.31);
\draw [thick][->] (3.5,0.2)--(3.5,-0.27);

\draw [black] (3.5,2.8) node {\small{$\bullet$}};
\draw (3.5,2.8) node [above] {\small{$\mathfrak{R}_4$}};

\draw [thick] (5,-1) arc (0:200:1.5 and 0.6);
\draw [thick] (5,-1)--(2.1,-1.2);
\draw [black] (5,-1) node {\tiny{$\bullet$}};
\draw [black] (2.1,-1.2) node {\tiny{$\bullet$}};
\draw [black] (3.6,-1.1) node {\tiny{$\bullet$}};
\draw [black] (4.9,-1.6) node {\tiny{$\bullet$}};

\draw [thick] (3.6,-1.1)--(5,-1)--(4.9,-1.6)--cycle;
\draw (2,1.8) node {$\Omega/{\sim}$};
\draw (1,-0.8) node {$\cD_0/\langle\sigma,\tau\rangle$};
\draw (3.5,-0.7) node {\small{$\cK$}};
\draw (4.65,-1.25) node {\small{$\cN$}};

\end{tikzpicture}
\caption{$\cD'$ as a compactification of the open cone $\Omega /{\sim} \cong C(\D^2)$.}
\label{cone_discs}
\end{center}
\end{figure}
\end{obs}

\begin{lem}
For every $s \in [0,1)$, the quotient $\big(h^{-1}(s)\smallsetminus \cD_s\big)/\langle\sigma,\tau\rangle$ is homeomorphic to $\R^3$.
\end{lem}

\begin{proof}
Recall from the proof of Lemma \ref{Lema:level_sphere_moduli} that for every $s \in (0,1)$, we have a fundamental domain of the action of $\langle \sigma, \tau, \rangle$ on $h^{-1}(s) \smallsetminus \cD_s$ given by $\big(\overline{\cU_1} \smallsetminus \cD\big) \cap h^{-1}(s)$, together with the region
\[ B_s =\big\{[0,1,z_3,z_4] \in \big(\overline{\cW_1} \smallsetminus \cD)\big) \cap h^{-1}(s) : \Im(z_3) \leq \Im(z_4)\big\}. \]
If we take the identifications on the sphere $\partial \big(\overline{\cU_1} \smallsetminus \cD\big) \cap h^{-1}(s) \cong \H^2$ determined by the actions of elements of the dihedral group $\langle\sigma,\tau\rangle$, we obtain a space homeomorphic to $\H^3$, whose boundary is obtained from the intersection $\big(\overline{\cU_1} \smallsetminus \cD\big) \cap \overline{\cW_4} \cap h^{-1}(s) \cong \R^2$. If we do the analogous identifications on $\partial B_s$, we obtain a space homeomorphic to $\H^3$, whose boundary is obtained from the intersection $\big(\overline{\cW_1} \smallsetminus \cD\big) \cap \overline{\cV_4} \cap h^{-1}(s) \cong \R^2$. As $\sigma^3 \tau (\overline{\cW_1} \cap \overline{\cV_4}) = \overline{\cW_4} \cap \overline{\cU_1}$, we conclude that after performing the identifications by elements of the dihedral group $\langle \sigma , \tau \rangle$ in the fundamental domain, we obtain two copies of $\H^3$ glued through their boundary which result in a space homeomorphic to $\R^3$. By making a limit process, and using Lemma \ref{lema:boundary_complement_of_wedges} and Lemma \ref{lemma_pegados_frontera}, we can observe that the same proof works for the case $h^{-1}(0) = \partial S(4)$ and the result follows.
\end{proof}

\begin{thm}\label{thm:boundary_decomposition_moduli}
The space $\cD' = \cD / \langle \sigma, \tau\rangle$ is homeomorphic to an cone $C(\D^2)$ compactified with the disc $\cD_0 / \langle \sigma, \tau\rangle$. Moreover, the space $\overline{\cM(4)} \smallsetminus \cD'$ is homeomorphic to $\H^4$.
\end{thm}

\begin{proof}
This is a consequence of taking the homeomorphisms of the previous Lemma, together with the cylinder structure given by Lemma \ref{conos_homeo_h4} and Theorem \ref{thm_gluing_upper_spaces} which tells us that $\overline{\cM(4)} \smallsetminus \cD' \cong \R^3 \times [0,1) \cong \H^4$.
\end{proof}

Theorem \ref{thm:boundary_decomposition_moduli} together with the description of $\cD'$ given in Lemma \ref{lema_cociente_disco} and Remark \ref{remark_compactification_cone_over_disc} completes the proof of Theorem \ref{thm_principal_three}.

\section{Conclusions and open questions}\label{section:open_questions}

In this section we make some general conclusions about the technical steps we had to face on the paper and include some questions which we encountered along the way.
\begin{enumerate}
    \item The definition of the height function $h$, given here, generalizes immediately for every $n \geq 3$. We can see however, that when $n \geq 5$, this function will have multiple points were it attains its maximum and thus, it doesn't behave as well as when $n \leq 4$. The authors believe that the following modification to the height function will work: 
    \[  h'([z_1,z_2,z_3,z_4]) = \frac{\min_{1\leq j\leq 4}\{r_j(Z)\}}{\max_{1\leq j,k\leq 4}|z_j-z_k|},    \]
    where $r_j(Z):=d(z_j,\cL_j)$ with $\cL_j:=\overline{z_{j+1} z_{j+2}} \cup \overline{z_{j+2} z_{j+3}} \cup \cdots \cup \overline{z_{j-2} z_{j-1}}$, compare with Definition \ref{distances}. We note that, as we are adding the length of the diagonals $|z_j - z_{j+k}|$ in the denominator of $h'$, the amount of subsets on which we must decompose $S(n)$ so that the function $h$ behaves well, increases considerably, compare with Definition \ref{conos}. 
    
    \item As explained in Lemma \ref{pegados}, the decomosition of $S(4)$ into the subsets $\{\overline{\cU_j}, \overline{\cV_j}, \overline{\cW_j}\}$ and the combinatorics of their intersections, are codified into a graph as a flag structure, dual to the graph. A combinatorically equivalent flag structure is constructed using cones, dual to the regular euclidean cuboctahedron in Theorem \ref{levelspheres}. Is there an intrinsic relation of this combinatorial structure to quadrilaterals? or the appearance of such a structure is a mere coincidence.
    
    \item The height function $h$ is differentiable on an open, dense subset of $S(4)$ and we have a well define gradient flow there, however, this flow arrives to the singular points in finite time. Is it possible to smooth the function, so that it is a differentiable Morse function with a unique maximum at the regular square? such function would induce a smooth regularizing flow on $S(4)$. It is worth noting that Theorem \ref{theo_cone_structure_simple_quadr} gives a continuous regularizing flow on $S(4)$, however in a non-constructive way. The question is thus, to obtain the results given in this paper, but in the differentiable category.
    
    \item As an open subspace of $\C\P^2$, $S(4)$ inherits the Fubini-Study metric which is no longer complete. Moreover, as both $S(4)$ and $\cM(4)$ are homeomorphic to cones over $\S^3$, we can give the spherical metric to each level set of the height function and the corresponding cone metrics on $S(4)$ and $\cM(4)$. Is there a more natural metric on those spaces whose geodesics capture better the deformations, such as the Weil-Petersson metric on Teichm\"uller spaces? This is a geometric version of the suavization problem of $h$ in the previous point.

    \item Is it true that a topological space $X$, with boundary homeomorphic to $\R^{n-1}$ and interior homeomorphic to $\R^n$, must be homeomorphic to $\H^n$? This assertion would have made the proof of Lemma \ref{conos_homeo_h4} shorter. What if we add the condition on $X$ to be a manifold with boundary?

    \item As exemplified in Example \ref{bad_convergence_ex} and explained in Remark \ref{remark_general_bundle_cylinder_extension}, the $I$-bundle structures described in the cones $\{\overline{\cU_j}, \overline{\cV_j}, \overline{\cW_j}\}$ on each level set $h^{-1}(s)$ don't converge good to the $I$-bundle structure described in the cones at the boundary $h^{-1}(0)$. We overcame this difficulty in the proof of Lemma \ref{conos_homeo_h4} by restricting the $I$-bundle structures to subsets where we had good convergence and then gluing the structures with Theorem \ref{thm_gluing_upper_spaces}. The conclusion of Lemma \ref{conos_homeo_h4}, tells us indirectly that there exist $I$-bundle structures which do have good convergence phenomena. How can we describe such $I$-bundle structures?

    \item If we consider the normalized side lengths of a quadrilateral
        \[   \ell_i : S(4) \rightarrow \R, \qquad \ell_i(Z) = \frac{2|z_{i} - z_{i+1}|}{\sum_{j = 1}^4 |z_{j} - z_{j+1}|},  \]
    then the function $f : S(4) \rightarrow \R^4$, $f(Z) = (\ell_1(Z), \ell_2(Z), \ell_3(Z), \ell_4(Z))$ was considered in \cite{KaMil1} and they proved that the fibers $f^{-1}(p)$ are one dimensional curves. It is possible to see that the level sets of the height function $h^{-1}(s)$ are transversal to those curves, and thus, the fibers of $f$ define a singular foliation of curves on $S(4)$. As an example of this, the curve $f^{-1}(1/2,1/2,1/2,1/2)$ is the fiber which contains the regular square and is defined by the intersection (see Remark \ref{recapitulacion_fibrado1}.III)
        \[   \bigcap_{j = 1}^4 \big(\textrm{cl}(\cU_j) \cap \textrm{cl}(\cV_j) \cap \textrm{cl}(\cW_j)\big).  \]

  \item The cone structure $S(4)$ give us the homeomorphisms $h^{-1}(0,1] \cong h^{-1}(\varepsilon,1] \cong S(4)$, for some $\varepsilon > 0$ and thus, it allows us to consider the compactification $\overline{S(4)}^{eucl} \cong  h^{-1}[\varepsilon,1]$, which is homeomorphic to a closed cone over $\S^3$. This compactification is invariant under the action of the dihedral group, so it gives a compactification of the moduli space $\overline{\cM(4)}^{eucl}$ as a closed cone over $\S^3$ as well and in particular $\overline{\cM(4)}^{eucl} \cong \D^4$. A reason to prefer the compactifications of $S(4)$ and $\cM(4)$ considered in this work is because it is compatible with the spaces considered in \cite{KaMil1}, but there is still the question on what properties capture each compactification and why should we prefer one or the other.

    \item How much can we weaken the hypotheses of Proposition \ref{prop:A.1}? As a more concrete example of this question, we can ask the following: If $G = N \rtimes H$ is a group acting by homeomorphisms on a topological space $X$, such that the action of $N$ descends to a well defined action on $(X/H)$, is it true that $(X/H) / N$ is homeomorphic to $X/G$? See Remark \ref{remark_quotients_appendix} and Example \ref{ex:A1}. We should emphasize that we do not have a reference for Remark \ref{remark_quotients_appendix} as it is considered a folklore result.

\end{enumerate}

\appendix
\section{Gluing closed balls and closed half-spaces}\label{section_gluing_balls}

We prove two topological results that are used many times along the paper. The ideas of the proofs of the two results in this appendix are extracted from \cite{rourke}. Recall the notations of the closed ball $\D^n=\{X\in\R^n\colon |X|\leq 1\}$ and the upper half-space $\H^n=\{(x_1,x_2,\dots,x_n)\in\R^n\colon x_n\geq0\}$.

\begin{thm}\label{thm_gluing_upper_spaces}
Let $H_i^n \cong \mathbb{H}^n$ be two copies of the upper half $n$-space, and $H_i^{n-1}\subset \partial H_i^n$, such that $H_i^{n-1} \cong \H^{n-1}$ and $\overline{\partial H_i^n \smallsetminus H_i^{n-1}} \cong \H^{n-1}$. For every homeomorphism $\varphi\colon H_1^{n-1}\to H_2^{n-1}$, the quotient $(H_1^n \cup H_2^n)\big/\varphi$ is homeomorphic to $\H^n$.
\end{thm}

\begin{proof}
We will use that $\H^n=\R^{n-1}\times[0,\infty)$ and the embedding $\H^{n-1}\hookrightarrow\H^n$ given by $(x_1,\dots,x_{n-1})\mapsto(x_1,\dots,x_{n-1},0)$. The key property of the proof is that every homeomorphism $g\colon\R^{n-1}\to\R^{n-1}$ can be extended to a homeomorphism $G\colon\H^n\to\H^n$ as follows
    \[  G(p,t)=(g(p),t)\qquad\forall~~p\in\R^{n-1}=\partial\H^n,~~\text{and}~~t\in[0,\infty).   \]
By applying the key property twice to a homeomorphism $f\colon H_1^{n-1}\to\H^{n-1}$, first extending to $\partial H_1^n$ and then to the whole $H_1^n$, we obtain a homeomorphism $F\colon H_1^n\to\H^n$ whose restriction to $H_1^{n-1}$ is $f$. The homeomorphism $F\sqcup F\colon H_1^n\sqcup H_1^n \to \H^n\sqcup\H^n$, induces a homeomorphism
    \[  \frac{H_1^n\sqcup H_1^n}{id_{H_1^{n-1}}}\cong \frac{\H^n\sqcup\H^n}{id_{\H^{n-1}}}\cong\H^n.  \]
By applying now, the key property twice to $\varphi\colon H_1^{n-1}\to H_2^{n-1}$, we obtain a homeomorphism $\Phi\colon H_1^n\to H_2^n$ that extends $\varphi$. The homeomorphism $Id \sqcup\Phi\colon H_1^n\sqcup H_1^n \to H_1^n \sqcup H_2^n$, induces a homeomorphism from $(H_1^n\sqcup H_1^n)\big/id_{H_1^{n-1}}$ to $(H_1^n \sqcup H_2^n)\big/\varphi$ and the result follows.
\end{proof}

\begin{thm}\label{thm_gluing_balls}
Let $D_i^n \cong \mathbb{D}^n$ be two copies of the closed $n$-ball, and $D_i^{n-1} \subset \partial D_i^n$ such that $D_i^{n-1}\cong \mathbb{D}^{n-1}$ and $\overline{\partial D_i^n \smallsetminus D_i^{n-1}} \cong \mathbb{D}^{n-1}$. For every homeomorphism $\varphi\colon D_1^{n-1}\to D_2^{n-1}$, the quotient $(D_1^n \cup D_2^n)\big/\varphi$ is homeomorphic to $\mathbb{D}^n$.
\end{thm}

\begin{proof}
The proof is analogous to the proof of Theorem \ref{thm_gluing_upper_spaces}. We just have to use the cone structure $\D^n\cong C(\S^{n-1})$ and the embedding $\D^{n-1}\hookrightarrow\D^n$ given by
\[(x_1,\dots,x_{n-1})\mapsto\left(x_1,\dots,x_{n-1},\sqrt{1-x_1^2-\cdots-x_{n-1}^2~}\right),\]
to get the property that every homeomorphism $g\colon \partial \mathbb{D}^n \rightarrow \partial \mathbb{D}^n$ extends to a homeomorphism of $\D^n$ as 
\[ G\colon\mathbb{D}^n \rightarrow \mathbb{D}^n, \quad  G(tp)=t g(p), \quad \forall~~ p \in \S^{n-1}= \partial \mathbb{D}^n.\]
\end{proof}

\begin{obs}
One way to obtain the hypothesis in Theorem \ref{thm_gluing_balls} that the complement of $D_i^{n-1}$ in $\partial D_i^n$ is again a ball is by adding PL-structures to the balls, as it is done for example in \cite[Corollary 3.13]{rourke}. Another way to check the hypothesis is to ask for the balls $D_i^{n-1}$ to be bi-collared in $\partial D_i^n$, so that the complement ball hypothesis follows from the solution of the Schoenflies conjecture by Brown in \cite{brown}.
\end{obs}

\section{Intermediate quotients by normal subgroups}\label{section_intermediate_quotients}
In this appendix, we give a proof of the following result which the authors consider is a well known result, but couldn't find a reference for it. 

\begin{prop}\label{prop:A.1}
If $G$ is a group acting on the right by homeomorphisms on a topological space $X$ and $N \leq G$ is normal subgroup, then there is a homeomorphism between the double quotient $(X / N) / (G/N)$ and $X / G$.
\end{prop}

\begin{proof}
Observe first that the projection map $\pi : X \rightarrow X/G$ descends to the quotient map $\widehat{\pi}: X/N \rightarrow X/G$. As $N$ is a normal subgroup of $G$, the $G$-action on $X$ induces a $G$-action on $X/N$ with $N$ acting trivially, so that it induces a $G/N$-action on $X/N$. As $\widehat{\pi}$ is $G/N$-invariant, it induces a continuous and bijective map $\widehat{\widehat{\pi}}: (X/N)/(G/N) \rightarrow X/G$. If $\widehat{\widehat{A}} \subset (X/N)/(G/N)$ is any subset, then there is a subset $A \subset X$ such that $\widehat{\widehat{A}} = \{ (x N) (G/N) : x \in A\}$, which means that $\widehat{\widehat{A}} = \{ \widehat{x} (G/N) : \widehat{x} \in \widehat{A}\}$, where $\widehat{A} = \{ y N : y \in A\}$. By definition, $\widehat{\widehat{A}}$ is open in $(X/N)/(G/N)$ if and only if its inverse image under the successive projection maps is open, that is,  if the set $ \{  (x g) : x \in A, g \in G \} = A G$ is open. As the projection $\pi$ is an open map, $ \pi(AG) = \pi(A) = \widehat{\widehat{\pi}}(\widehat{\widehat{A}})$ is open in $X/G$ and $\widehat{\widehat{\pi}}$ is an homeomorphism.
\end{proof}

\begin{ejem}[non-example]\label{ex:A1}
Consider the dihedral group $D_4 = \langle \sigma \rangle \rtimes \langle \tau \rangle$ acting in the circle $\S^1 \subset \C$, as $\sigma(z) = e^{\pi i/2 } z$ and $\tau(z) = \overline{z}$. The quotient $\S^1 / \langle \tau \rangle$ is homeomorphic to $\{z \in \S^1 : \Im(z) \geq 0 \}$ and the homeomorphism $\sigma$ descends to a continuous map on $\S^1 /  \langle \tau \rangle$, defined under this identification by
    \[  \Tilde{\sigma}(z) = \left\{\begin{array}{lcr}
         e^{\pi i/2 }z, & & \Re(z) \geq 0,  \\
         e^{3 \pi i/2 }\overline{z}, & & \Re(z) \leq 0.
    \end{array} \right.    \]
Under the equivalence class $z \sim \Tilde{\sigma}(z)$, we have an homeomorphism 
    \[  (\S^1 / \langle \tau \rangle) / \Tilde{\sigma} = (\S^1 / \langle \tau \rangle) / \sim \cong \S^1 / D_4. \]
In this case however, the map $\Tilde{\sigma}^k$ is non-injective, and thus, the action of the cyclic group $ \langle \sigma \rangle$ on $\S^1$ doen't descend to a well defined action on $\S^1 / \langle \tau \rangle$.
\end{ejem}

\begin{obs}\label{remark_quotients_appendix}
As we saw in Example \ref{ex:A1}, there are cases where the condition of the subgroup on Proposition \ref{prop:A.1} to be normal is no longer true and nevertheless the result holds. This is not always the case as there examples of two non-commuting homeomorphisms $f_1, f_2 : X \rightarrow X$, such that for $X_j = X / \langle f_j \rangle$ the quotient $X_1 / \langle f_2 \rangle$ is not homeomorphic to $X_2 / \langle f_1 \rangle$ and both spaces are not homeomorphic to $X / \langle f_1 , f_2 \rangle$.
\end{obs}


\

\noindent
{\small Ahtziri Gonz\'alez, Facultad de Ingeniería Eléctrica, Universidad Michoacana de San Nicol\'as de Hidalgo, Edificio Omega, Ciudad Universitaria, C.P. 58040, Morelia, Michoac\'an, M\'exico.
\noindent
e-mail: ahtziri.lemus@umich.mx}    \\

\noindent
{\small Manuel Sedano-Mendoza,
\noindent
e-mail: msedano@matmor.unam.mx}

\begin{thebibliography}{20}

\bibitem{BEH} \textsc{K Behrend}, \textit{Introduction to Algebraic Stacks}, Lecture Note Series  411, London Mathematical Society, (2014), pp. 1-131.

\bibitem{AHT} \textsc{A Gonz\'alez}, ``\textit{Fundamental group of spaces of Simple Polygons}'' Rocky Mountain J. Math. 48 (6), 2018, pp. 1871-1886.

\bibitem{AHT-JO} \textsc{A Gonz\'alez \& J L\'opez-L\'opez}, \textit{Compactness of Spaces of Convex and Simple Quadrilaterals}, Bulletin of the Australian Mathematical Society, 94(3), 2016, pp. 507-521. doi:10.1017/S0004972716000368

\bibitem{AHT-JO1} \textsc{A Gonz\'alez \& J L\'opez-L\'opez}, \textit{Spaces of Special Quadrilaterals}, Bull. Aust. Math. Soc., vol. 100 (2019), pp. 155-167.

\bibitem{HAT}  \textsc{A. Hatcher}, \textit{Algebraic Topology}, Cambridge Univ. Press, (2002).

\bibitem{HUS} \textsc{D Husemoller}, Fibre Bundles, Graduate Text in Mathematics No. 20. Springer-Verlag ,Third Edition, 1993.

\bibitem{SIN-JON} \textsc{G. A. Jones \& D. Singerman}, \textit{Complex Functions: An Algebraic and Geometric Viewpoint}, Cambridge Univ. Press, (1987).

\bibitem{THU} \textsc{W. Thurston}, \textit{Shapes of polyhedra and triangulations of the sphere, in: The Epstein Birthday Schrift}, Geom. Topol. Monogr., vol. 1 (1998), pp. 511-549.


\bibitem{Cant} \textsc{John Cantwell}, \textit{Topological non-degenerate functions}, T\^ohoku Math Journal, 20(1968), 120-125.

\bibitem{Sch} \textsc{Paul Schmutz Schaller}, \textit{Systoles and topological Morse functions for Riemann surfaces}, J Differential Geometry 52 (1999) 407--452.

\bibitem{Morse} \textsc{M. Morse}, {\it Topologically non-degenerate functions on a compact $n$-manifold M}, Journ. d'Analyse Math., 7(1959), 189-208

\bibitem{Kui} \textsc{N. H. Kuiper}, \textit{A continuous function with two critical points}, Bull. Amer. Math. Soc., 67(1961), 281-285.

\bibitem{ChowGli} \textsc{B. Chow and D. Glickenstein}, \textit{Semidiscrete geometric flows of polygons}, Amer. Math. Monthly 114 (2007), no. 4, 316-328.

\bibitem{brown} \textsc{Morton Brown}, {\it A proof of the generalized Schoenflies theorem}. Bull. Amer. Math. Soc. 66 (1960), 74–76. 

\bibitem{brown2} \textsc{Morton Brown}, {\it Locally flat imbeddings of topological manifolds}, Annals of Mathematics, Vol. 75 (1962), 331-341.
 
\bibitem{rourke} \textsc{ Rourke, Colin Patrick and Sanderson, Brian Joseph}, {\it Introduction to piecewise-linear topology}, Reprint Springer Study Edition Springer-Verlag, Berlin-New York, 1982.

\bibitem{fill} \textsc{Fillastre, François}, {\it From spaces of polygons to spaces of polyhedra following Bavard, Ghys and Thurston} Enseign. Math. (2) 57 (2011), no. 1-2, 23–56. 

\bibitem{Haus} \textsc{J.-C.Hausmann}, {\it Sur la topologie des bras articules}, (Algebraic Topology, Poznan), Lecture Notes in Math. Vol. 1474, Springer, Berlin, 1989, 146- 159.

\bibitem{KaMil1} \textsc{M. Kapovich and J. Millson}, \textit{On the moduli space of polygons in the Euclidean plane}, J Differential Geometry, Vol. 42, No.1 July, 1995.

\bibitem{KaMil2} \textsc{M. Kapovich and J. Millson}, {\it The symplectic geometry of polygons in Euclidean space}, J. Differential Geom. 44 (1996), no. 3, 479–513. 

\bibitem{KaMil3} \textsc{M. Kapovich and J. Millson}, {\it Quantization of bending deformations of polygons in E3, hypergeometric integrals and the Gassner representation}, Canad. Math. Bull. 44 (2001), no. 1, 36–60. 

\bibitem{DeMo} \textsc{P. Deligne and G. D. Mostow}, {\it  Monodromy of hypergeometric functions and nonlattice integral monodromy}, Inst. Hautes Études Sci. Publ. Math. No. 63 (1986), 5–89. 

\bibitem{BaGhys} \textsc{C. Bavard and E. Ghys}, {\it Polygones du plan et polyèdres hyperboliques}, Geom. Dedicata 43 (1992), no. 2, 207–224.

\bibitem{KoYam} \textsc{S. Kojima and Y. Yamashita}, {\it Shapes of stars}, Proc. Amer. Math. Soc. 117 (1993), no. 3, 845–851. 

\bibitem{Loz} \textsc{Lozano-Pérez, Tomás}, {\it Spatial planning: a configuration space approach}, IEEE Trans. Comput. 32 (1983), no. 2, 108–120.

\bibitem{Gott} \textsc{D. H. Gottlieb}, {\it Robots and fibre bundles}, Bull. Soc. Math. Belg. Sér. A 38 (1986), 219–223 (1987). 

\bibitem{HauKnu} \textsc{J-C. Hausmann and A. Knutson}, {\it The cohomology ring of polygon spaces}, Ann. Inst. Fourier (Grenoble) 48 (1998), no. 1, 281–321. 

\bibitem{SmBrouFra} \textsc{S.L. Smith, E. M. Broucke and B. A. Francis}, {\it Curve shortening and the rendezvous problem for mobile autonomous robots}, IEEE Trans. Automat. Control 52 (2007), no. 6, 1154–1159.

\bibitem{CoDeRo} \textsc{R. Connelly, E. D. Demaine and G. Rote}, {\it Straightening polygonal arcs and convexifying polygonal cycles}, Discrete Comput. Geom. 30 (2003), no. 2, 205–239. 

\bibitem{Okoshita} \textsc{K. Ohshika}, {\it Compactifications of Teichm\"uller spaces}, Handbook of Teichm\"uller Theory, Volume IV, (2014).

\bibitem{bestvina} \textsc{M. Bestvina}, {\it The topology of $Out(Fn)$}, Proc. of the ICM, Vol. II (Beijing, 2002), 373–384, Higher Ed. Press, Beijing, 2002. 

\end{thebibliography}
\end{document}